\documentclass[11pt]{article}
\usepackage[english]{babel}
\usepackage[latin1]{inputenc}
\usepackage{amsmath}
\usepackage{amssymb}
\usepackage{amsfonts}
\usepackage{amsthm}
\usepackage{upgreek}
\usepackage{stmaryrd}

\newcommand{\C}{\mathbb{C}}
\newcommand{\N}{\mathbb{N}}

\newcommand{\R}{\mathbb{R}}
\renewcommand{\S}{\mathbb{S}}

\newcommand{\boD}{\mathcal{D}}

\newcommand{\boS}{\mathcal{S}}

\newcommand{\eps}{\varepsilon}

\newtheorem{cor}{Corollary}
\newtheorem{lemma}{Lemma}
\newtheorem{prop}{Proposition}
\newtheorem{step}{Step}
\newtheorem{theorem}{Theorem}
\newtheorem*{theoremn}{Theorem}
\theoremstyle{definition}

\newtheorem{remark}{Remark}

\begin{document}

\title{Asymptotics of the solitary waves for the generalised Kadomtsev-Petviashvili equations}
\author{
\renewcommand{\thefootnote}{\arabic{footnote}} Philippe Gravejat \footnotemark[1]}
\footnotetext[1]{Centre de Recherche en Math\'ematiques de la D\'ecision, Universit\'e Paris Dauphine, Place du Mar\'echal De Lattre De Tassigny, 75775 Paris Cedex 16, France. E-mail: gravejat@ceremade.dauphine.fr}
\date{}
\maketitle

\begin{abstract}
We investigate the asymptotic behaviour of the localised solitary waves for the generalised Kadomtsev-Petviashvili equations. In particular, we compute their first order asymptotics in any dimension $N \geq 2$.
\end{abstract}

\section{Introduction}

\subsection{Motivation and main results}

The present paper deals with the solitary waves for the generalised Kadomtsev-Petviashvili equations
\begin{equation}\label{E1}
\left\{\begin{array}{ll} \partial_t u + u^p \partial_1 u + \partial^3_1 u - \underset{j=2}{\overset{N}{\sum}} \partial_j u_j = 0, \\ \forall j \in \{2,\ldots,N\}, \partial_1 u_j = \partial_j u. \end{array} \right.
\end{equation}
The exponent $p = \frac{m}{n}$ is rational (where $m$ and $n$ are relatively prime, and $n$ is odd). The function $u \mapsto u^p$ is defined by the standard convention
$$\forall u \in \R, u^p = Sign(u)^m |u|^p.$$

Two cases at least are physically relevant. First, the case $p=1$ corresponds to the standard Kadomtsev-Petviashvili equation: a universal model for dispersive, weakly nonlinear long waves, essentially unidimensional in the direction of propagation $x_1$ (see e.g. \cite{KadoPet1}). Second, the case $p = 2$ appears as a model for the evolution of sound waves in antiferromagnetics (see e.g. \cite{FalkTur1}).

The solitary waves for the generalised Kadomtsev-Petviashvili equations are the solutions $u$ of \eqref{E1} of the form
$$u(t,x) = v(x_1 - ct,x_\perp), \ x_\perp = (x_2,\ldots,x_N),$$
which belong to the closure $Y$ of the space $\partial_1 C_0^\infty(\R^N)$ for the norm
$$\forall \phi \in C_0^\infty(\R^N), \| \partial_1 \phi \|_Y = \bigg(\| \nabla \phi\|^2_{L^2(\R^N)} + \| \partial_{1,1}^2 \phi\|^2_{L^2(\R^N)}\bigg)^\frac{1}{2}.$$
They are (at least formally) critical points on $Y$ of the action $S$ defined by
\begin{equation}\label{E2}
\forall v \in Y, S(v) = E(v) + \frac{c}{2} \int_{\R^N} v^2(x) dx,
\end{equation}
where $E$ denotes the energy associated to the generalised Kadomtsev-Petviashvili equations
\begin{equation}\label{E3}
E(v) = \frac{1}{2} \int_{\R^N} \Big( \partial_1 v(x)^2 + \sum_{j=2}^N v_j(x)^2 \Big) dx -\frac{1}{(p+1)(p+2)} \int_{\R^N} v(x)^{p+2} dx.
\end{equation}

The parameter $c > 0$ is the speed of the solitary wave, which moves in the direction $x_1$. We can always make the additional assumption $c = 1$. Indeed, if $v$ is a solitary wave with speed $c$, the equation for $v$ is written as
\begin{equation}\label{E4}
\left\{\begin{array}{ll} - c \partial_1 v + v^p \partial_1 v + \partial^3_1 v - \underset{j=2}{\overset{N}{\sum}} \partial_j v_j = 0, \\ \forall j \in \{2,\ldots,N\}, \partial_1 v_j = \partial_j v. \end{array} \right.
\end{equation}
Therefore, the function $\tilde{v}$, given by the scale changes,
\begin{equation}\label{E5}
\forall x \in \R^N, \tilde{v}(x_1,x_\perp) = c^{-\frac{1}{p}} v \bigg( \frac{x_1}{\sqrt{c}}, \frac{x_\perp}{c} \bigg),
\end{equation}
is a solitary wave with speed $1$. In order to simplify the notation, we will assume from now on that $c$ is equal to $1$. We will recover the arbitrary case by the scale change \eqref{E5}. In particular, with this additional hypothesis, the solitary wave $v$ solves the equation
\begin{equation}\label{E6}
- \Delta v + \partial^4_1 v + \frac{1}{p+1} \partial^2_1 (v^{p+1}) = 0,
\end{equation}
which is the starting point of our analysis.

A. de Bouard and J.-C. Saut first studied the existence and qualitative properties of the solitary waves for the generalised Kadomtsev-Petviashvili equations. In \cite{deBoSau1}, they completely solved the issue of their existence in dimensions two and three: they proved that there exist non-trivial solutions of equation \eqref{E4} in $Y$ if and only if
$$0 < p < \frac{4}{2N - 3}.$$
We believe that their proof extends to any dimension $N \geq 4$
\footnote{In this article, we will derive from integral identities (which are of independent interest) that there are no non-trivial solitary-wave solutions of equation \eqref{E1} in $Y$ if $p \geq \frac{4}{2N - 3}$ (see Corollary \ref{C2}). However, our goal is not to obtain existence results, so we will not consider these existence problems any further.}.
Therefore, we will only consider the exponents $p$ for which there presumably exist non-trivial solutions of equation \eqref{E4} in $Y$, that is,
$$0 < p < \frac{4}{2N - 3}.$$
Moreover, A. de Bouard and J.-C. Saut \cite{deBoSau3} also addressed the issue of the dynamical stability of the family of solitary waves. In dimension two, they proved their orbital stability for $p < \frac{4}{3}$, and their instability for $p > \frac{4}{3}$. Likewise, J.L. Bona and Yue Liu \cite{BonaLiu1} showed the instability of the solitary waves in dimension three.\\
Finally, A. de Bouard and J.-C. Saut \cite{deBoSau2} proved the axisymmetry around the axis $x_1$ of the ground states (the solitary waves which minimise the action $S$ on the space $Y$), and computed the algebraic decay of any solitary wave in dimensions two and three.

\begin{theoremn}[\cite{deBoSau2}]
In dimension two, any solitary wave $v$ of equation \eqref{E1}
satisfies
$$r^2 v \in L^\infty(\R^2), r^2 = x_1^2 + x_2^2.$$
In dimension three, any solitary wave $v$ of equation \eqref{E1}
satisfies
$$\forall 0 \leq \delta < \frac{3}{2}, r^\delta v \in L^2(\R^3), r^2 =
x_1^2 + x_2^2 + x_3^2.$$
\end{theoremn}

\begin{remark}
Their theorem is sharp in dimension two for $p = 1$. Indeed, there exists an explicit solution of equation \eqref{E4} in dimension two, the so-called lump solution $v_c$ given by
$$\forall (x_1,x_2) \in \R^2, v_c(x_1,x_2) = 24c\frac{3 - c x_1^2 +
 c^2 x_2^2}{(3 + c x_1^2 + c^2 x_2^2)^2}.$$
In particular, we cannot expect a decay rate better than $r^{-2}$ in dimension two for $p=1$.
\end{remark}

The goal of this paper is precisely to improve their description of the asymptotics of a solitary wave in any dimension $N \geq 2$.
 
\begin{theorem}\label{T1}
Let $v \in Y$ be a solitary-wave solution of speed $1$ of equation \eqref{E1}. Assume that
$$0 < p < \frac{4}{2N - 3},$$
and consider the function $v_\infty \in C^\infty(\S^{N-1})$ given by
\begin{equation}\label{E7}
\forall \sigma = (\sigma_1,\ldots,\sigma_N) \in \S^{N-1}, v_\infty(\sigma) = \frac{\Gamma(\frac{N}{2})}{2 \pi^\frac{N}{2} (p+1)}(1 - N \sigma_1^2) \int_{\R^N} v(x)^{p+1} dx.
\end{equation}
Then, the function $x \mapsto |x|^N v(x)$ is bounded on $\R^N$, and
\begin{equation}\label{E8}
\forall \sigma \in \S^{N-1}, R^N v(R \sigma) \underset{R \to + \infty}{\to} v_\infty(\sigma).
\end{equation}
Moreover, if $\frac{1}{N} \leq p < \frac{4}{2N - 3}$, this convergence is uniform: it holds in $L^\infty(\S^{N-1})$.
\end{theorem}

\begin{remark}
The function $v_\infty$ is well-defined on $\S^{N-1}$. Indeed, by Theorem \ref{T8}, the integral $\int_{\R^N} v(x)^{p+1} dx$ is finite.
\end{remark}

Theorem \ref{T1} gives the sharp decay rate at infinity, which is exactly equal to $r^{-N}$, of any non-trivial solitary wave $v$ in any dimension $N \geq 2$ for all the exponents $p = \frac{m}{n}$ such that $m$ is an odd number. Indeed, let us consider some non-trivial solitary wave $v$ in $Y$ such that its decay rate is strictly more than $r^{-N}$. Then, there exists some function $d: \R_+ \to \R_+$ such that
$$\forall x \in \R^N, |v(x)| \leq d(|x|),$$
and
$$r^N d(r) \underset{r \to + \infty}{\to} 0.$$
By Theorem \ref{T1}, the function $v_\infty$ is identically equal
to $0$ on $\S^{N-1}$, which gives
$$\int_{\R^N} v(x)^{p+1} dx = 0.$$
Then, $v$ is trivial since $m$ is odd, which leads to a contradiction. Thus, Theorem \ref{T1} is optimal for any non-trivial solitary wave when $m$ is an odd number (which holds for the standard Kadomtsev-Petviashvili equation).

On the other hand, Theorem \ref{T1} may not be sharp if $m$ is even. There may be non-trivial solitary waves whose decay rate is higher than $r^{-N}$. This may happen if the function $v_\infty$ is identically equal to $0$, that is if
$$\int_{\R^N} v(x)^{p+1} dx = 0.$$
We do not actually know of any non-trivial solitary waves which verify such assumption, but we believe that they may exist. Indeed, L. Paumond \cite{Paumond1} proved their existence in dimension $N = 5$ for an equation very similar to equation \eqref{E1}, namely
$$\left\{\begin{array}{ll} \partial_t u + u^p \partial_1 u + \partial^7_1 u - \underset{j=2}{\overset{5}{\sum}} \partial_j u_j = 0, \\ \forall j \in \{2,\ldots,5\}, \partial_1 u_j = \partial_j u. \end{array} \right.$$
More precisely, when $m$ is even, he proved the existence of non-trivial solitary waves $v$ for this equation which verify
$$\forall x \in \R^5, v(x_1,x_2,x_3,x_4,x_5) = - v(x_1,x_4,x_5,x_2,x_3).$$
In particular, such solutions satisfy the condition
$$\int_{\R^5} v(x)^{p+1} dx = 0$$
above. Therefore, Theorem \ref{T1} is not necessarily sharp when $m$ is even.

Before sketching the proof of Theorem \ref{T1}, we must notice another improvement for the standard Kadomtsev-Petviashvili equation.

\begin{theorem}\label{T2}
Let $v \in Y$ be a solitary-wave solution of speed $1$ of equation \eqref{E1}. Assume that $N=2$ or $N=3$, and $p=1$. Then, the function $v_\infty$ is given for every $\sigma \in \S^{N-1}$ by
\begin{equation}\label{E9}
v_\infty(\sigma) = \frac{(7 - 2N) \Gamma(\frac{N}{2})}{2 (2N - 5) \pi^\frac{N}{2}} (1 - N \sigma_1^2) \ E(v) = \frac{(7 - 2N) \Gamma(\frac{N}{2})}{4 \pi^\frac{N}{2}} (1 - N \sigma_1^2) \ S(v).
\end{equation}
\end{theorem}

Theorem \ref{T2} results from Pohozaev identities, which were already derived by A. de Bouard and J.-C. Saut in \cite{deBoSau1}. It links the asymptotics of a solitary wave $v$ to its energy $E(v)$ or its action $S(v)$. In particular, there is only one possible asymptotics for all the solitary waves with the same energy. This seems to be a further evidence of the uniqueness of non-trivial solitary waves (up to translations) in the case of the standard Kadomtsev-Petviashvili equation: we believe that this new evidence could be a useful step towards the resolution of this problem, which is still open to our knowledge.

\subsection{Sketch of the proof of Theorem \ref{T1}}

The goal of Theorem \ref{T1} is to compute the algebraic decay at infinity of a solitary wave $v$, then its first order asymptotics. Our proof is reminiscent of a series of articles by J.L. Bona and Yi A. Li \cite{BonaLi2}, A. de Bouard and J.-C. Saut \cite{deBoSau2}, M. Maris \cite{Maris1,Maris2}, and \cite{Graveja3,Graveja5}. It relies on the use of convolution equations, in particular on a precise analysis of the kernels they involve. The analysis provides the improvements mentioned above. Moreover, it relies on general arguments which could prove fruitful for other equations (with other kernels and nonlinearities) in any dimension. That is the reason why we first explain the main arguments of the proof of Theorem \ref{T1}. We hope that it will also help to clarify this proof.

\subsubsection{Convolution equations}

By equation \eqref{E6}, the solitary wave $v$ satisfies, at least formally, both the convolution equations
\begin{equation}\label{E10}
v = i H_0*(v^p \partial_1 v)
\end{equation}
and
\begin{equation}\label{E11}
v = \frac{1}{p+1} K_0*v^{p+1}.
\end{equation}
Here, $H_0$ and $K_0$ are the kernels whose Fourier transforms are
\begin{equation}\label{E12}
\widehat{H_0}(\xi) = \frac{\xi_1}{|\xi|^2 + \xi_1^4}
\end{equation}
and
\begin{equation}\label{E13}
\widehat{K_0}(\xi) = \frac{\xi_1^2}{|\xi|^2 + \xi_1^4}.
\end{equation}
Equations \eqref{E10} and \eqref{E11} link the asymptotic properties of $v$ to the behaviour of $H_0$ and $K_0$ at infinity. This requires a careful analysis of such kernels, first for deriving equations \eqref{E10} and \eqref{E11} rigorously, second for computing the asymptotics of $v$.

\subsubsection{Main properties of the kernels}

This section is devoted to the study of the kernels $H_0$, $K_0$ and $K_k = - i \partial_k K_0$, given by
\begin{equation}\label{E14}
\forall k \in \{ 1, \ldots, N \}, \widehat{K_k}(\xi) = \frac{\xi_k
 \xi_1^2}{|\xi|^2 + \xi_1^4}.
\end{equation}
In view of the comment above, we first describe their asymptotic properties and their singularities near the origin in order to deduce the asymptotic properties of $v$ later on.

\subsubsection*{Algebraic decay at infinity and singularities near the origin}

Let us consider the spaces of functions $M^\infty_\alpha(\Omega)$ defined by
$$M^\infty_\alpha(\Omega) = \{ u:\Omega \mapsto \C, \| u \|_{M^\infty_\alpha(\Omega)} = \sup \{|x|^\alpha |u(x)|, x \in \Omega \} < + \infty \},$$
for any $\alpha > 0$ and any open subset $\Omega$ of $\R^N$. We say that a function $f$ presents some algebraic decay at infinity if it belongs to some space $M^\infty_\alpha(B(0,1)^c)$ for some $\alpha > 0$. Likewise, $f$ presents some algebraic explosion near the origin if it belongs to some space $M^\infty_\alpha(B(0,1))$ for some $\alpha > 0$.

One goal of Theorem \ref{T1} is to derive the algebraic decay of the solitary wave $v$ at infinity. The first step consists in studying the algebraic decay of $H_0$, $K_0$ and $K_k$ at infinity and their algebraic explosion near the origin. This study first provides their algebraic decay at infinity.

\begin{theorem}\label{T3}
Let $k \in \{1, \ldots, N \}$. The kernels $H_0$, $K_0$ and $K_k$ are continuous on $B(0,1)^c$, and belong to $M^\infty_{N-1}(B(0,1)^c)$, $M^\infty_N(B(0,1)^c)$ and $M^\infty_{N+1}(B(0,1)^c)$ respectively.
\end{theorem}

\begin{remark}
We believe that Theorem \ref{T3} is sharp in the sense that the kernels $H_0$, $K_0$ and $K_k$ do not belong to $M^\infty_{\alpha-1}(B(0,1)^c)$, $M^\infty_\alpha(B(0,1)^c)$ and $M^\infty_{\alpha+1}(B(0,1)^c)$ for any $\alpha > N$.
\end{remark}

This study also yields a description of their singularities near the origin.

\begin{theorem}\label{T4}
Let $1 \leq k \leq N$. There exists some $A > 0$ such that for every $x \in B(0,1)$,
$$\left\{ \begin{array}{lll} (x_1^2 + |x_\perp|)^{N-2} \ |H_0(x)| \leq A(1 + \delta_{N,2} |\ln(|x|)|), \\ (x_1^2 + |x_\perp|)^{N-\frac{3}{2}} \ |K_0(x)| \leq A, \\ (x_1^2 + |x_\perp|)^{N - \frac{1 + \delta_{k,1}}{2}} \ |K_k(x)| \leq A. \end{array} \right.$$
\end{theorem}

\begin{remark}
By Theorem \ref{T4}, $H_0$, $K_0$ and $K_k$ are respectively in $M^\infty_{2N-4}(B(0,1))$, $M^\infty_{2N-3}(B(0,1))$ and $M^\infty_{2N-1-\delta_{k,1}}(B(0,1))$ in any dimension $N > 2$. However, such spaces are not suitable to describe their singularities near the origin. Indeed, their singularities are anisotropic because of the anisotropy of their Fourier transforms given by formulae \eqref{E12}, \eqref{E13} and \eqref{E14}.
\end{remark}

\begin{remark}
We believe that Theorem \ref{T4} is sharp in the sense that it provides the right exponents of the singularities of $H_0$, $K_0$ and $K_k$ near the origin.
\end{remark}

This description is also available in terms of $L^q$-spaces.

\begin{cor}\label{C1}
Let $1 \leq j,k \leq N$ and $q \in [1, + \infty[$. The functions $H_0$, $K_0$ and $x_j^{1-\delta_{k,1}} K_k$ respectively belong to $L^q(B(0,1))$ if $q < \frac{2N - 1}{2N -4}$, $q < \frac{2N - 1}{2N - 3}$ and $q < \frac{2N - 1}{2N - 2}$.
\end{cor}

The proofs of Theorems \ref{T3} and \ref{T4} rely on the form of the Fourier transforms of $H_0$, $K_0$ and $K_k$. By formulae \eqref{E12}, \eqref{E13} and \eqref{E14}, they are rational fractions given by
\begin{equation}\label{E15}
\forall \xi \in \R^N \setminus \{ 0 \}, R(\xi) = \frac{P(\xi)}{|\xi|^2 + \xi_1^4},
\end{equation}
where $P$ is some polynomial function on $\R^N$ of the form
\begin{equation}\label{E16}
\forall \xi \in \R^N, P(\xi) = \underset{j=1}{\overset{N}{\Pi}} \xi_j^{d_j}.
\end{equation}
Theorems \ref{T3} and \ref{T4}, and Corollary \ref{C1} follow from more general results for tempered distributions $f$ whose Fourier transforms $\widehat{f} = R$ are rational fractions of the form \eqref{E15}-\eqref{E16}. Indeed, we can compute explicitly the algebraic decay of such distributions.

\begin{prop}\label{Paumond1}
Let $f$ be a tempered distribution on $\R^N$ whose Fourier transform is of the form \eqref{E15}-\eqref{E16}, and denote
\begin{equation}\label{E17}
d = \underset{j = 1}{\overset{N}{\sum}} d_j = d_1 + d_\perp.
\end{equation}
Assume moreover that $d \neq 0$ if $N = 2$ and $d_1 + 2 d_\perp \leq 4$. Then, $f$ belongs to the space $M^\infty_{N-2+d}(B(0,1)^c)$.
\end{prop}

Likewise, we can describe their singularities near the origin.

\begin{prop}\label{P2}
Let $f$ be a tempered distribution on $\R^N$ whose Fourier transform is of the form \eqref{E15}-\eqref{E16}. Assume moreover that $d \neq 0$ if $N=2$, and $d_1 + 2 d_\perp \leq 4$ where $d$ and $d_\perp$ are defined by \eqref{E17}. Then, there exists some $A > 0$ such that for every $x \in B(0,1) \setminus \{ 0 \}$,
\begin{equation}\label{E18}
(x_1^2 + |x_\perp|)^{N - \frac{5}{2} + \frac{d_1}{2} + d_\perp} |f(x)| \leq A (1 + \delta_{N,2} \delta_{d_1,1} \delta_{d_\perp,0} |\ln(|x|)|).
\end{equation}
In particular, if $d_1 + 2d_\perp < 4$, $f$ belongs to $L^q(B(0,1))$ for every
\begin{equation}\label{E19}
1 \leq q < \frac{2N - 1}{2N - 5 + d_1 + 2d_\perp}.
\end{equation}
Likewise, if $(d_1, d_\perp) = (2,1)$ or $(d_1, d_\perp) = (4,0)$, the functions $x \mapsto x_j f(x)$ belong to $L^q(B(0,1))$ for every
\begin{equation}\label{E20}
1 \leq q < \frac{2N - 1}{2N - 6 + d_1 + 2d_\perp}.
\end{equation}
\end{prop}

\begin{remark}
When $(d_1, d_\perp) = (2,1)$ or $(d_1, d_\perp) = (4,0)$, the distribution $f$ is not a function in $L^1_{loc}(B(0,1))$. The singularities of $f$ near the origin can present some principal values at the origin or some Dirac masses (see Lemma \ref{L3}). However, the distributions $x_j f$ are in $L^1_{loc}(B(0,1))$, so we can study their $L^q$-integrability.
\end{remark}

The first step of the proof of Propositions \ref{Paumond1} and \ref{P2} consists in inductively describing the derivatives of $\widehat{f} = R$, in particular their singularities near the origin and their integrability at infinity. Then, by standard integral expressions, we deduce the behaviour of $f$ near the origin and at infinity. More precisely, we first derive the form of the derivatives of $R$.

\begin{prop}\label{P3}
Let $1 \leq j \leq N$ and $p \in \N$. Let us consider a rational fraction $R$ on $\R^N$ which satisfies formulae \eqref{E15} and \eqref{E16}. Then, the partial derivative $\partial_j^p R$ is given by
\begin{equation}\label{E21}
\forall \xi \in \R^N \setminus \{ 0 \}, \partial_j^p R(\xi) = \frac{P_{j,p}(\xi)}{(|\xi|^2 + \xi_1^4)^{p+1}},
\end{equation}
where $P_{j,p}$ is a polynomial function on $\R^N$. Moreover, there exist some $A_p > 0$ such that the function $P_{1,p}$ satisfies
\begin{align}
\forall \xi \in B(0,1)^c, |P_{1,p}(\xi)| & \leq A_p \sum_{0 \leq k \leq \frac{d_1 + 3p}{4}} \max \{ |\xi_1|, 1 \}^{d_1 + 3p - 4k} |\xi_\perp|^{2k + d_\perp}, \label{E22} \\ \forall \xi \in B(0,1), |P_{1,p}(\xi)| & \leq A_p |\xi|^{p + d_1 + d_\perp}, \label{E23}
\end{align}
where $d_\perp$ is defined by \eqref{E17}. Likewise, if $j \geq 2$, the function $P_{j,p}$ verifies
\begin{align}
\forall \xi \in B(0,1)^c, |P_{j,p}(\xi)| \leq & A_p \sum_{0 \leq k \leq \frac{d_\perp + p}{2}} \max \{ |\xi_\perp|, 1 \}^{d_\perp + p - 2k} \max \{ |\xi_1|, 1 \}^{d_1 + 4k}, \label{E24} \\ \forall \xi \in B(0,1), |P_{1,p}(\xi)| \leq & A_p |\xi|^{p + d_1 + d_\perp}. \label{E25}
\end{align}
\end{prop}

\begin{remark}
By linearity, similar estimates hold for any rational fractions of the form \eqref{E15} (where $P$ is any polynomial function on $\R^N$). However, all the rational fractions considered here will satisfy \eqref{E15} and \eqref{E16}, so, in order to simplify some computations, we will not investigate this point any further.
\end{remark}

The analysis of Proposition \ref{P3} is sufficient to describe the singularities of the derivatives of $R$ near the origin and their integrability at infinity.

\begin{prop}\label{P4}
Let $1 \leq j \leq N$, $p \in \N$ and $q \in [1, + \infty]$. Let us consider a rational fraction $R$ on $\R^N$ which is of the form \eqref{E15}-\eqref{E16}, and denote $d = \underset{j=1}{\overset{N}{\sum}} d_j = d_1 + d_\perp$. Then, the partial derivative $\partial_j^p R$ belongs to $M^\infty_{p + 2 - d}(B(0,1))$. Moreover, if
\begin{equation}\label{E26}
p > d_1 + 2 d_\perp - 4,
\end{equation}
the partial derivative $\partial_1^p R$ belongs to $L^q(B(0,1)^c)$ for
\begin{equation}\label{E27}
q > \frac{2N-1}{p + 4 - d_1 - 2d_\perp},
\end{equation}
while it belongs to $L^\infty(B(0,1)^c)$ for
\begin{equation}\label{E28}
p \geq d_1 + 2 d_\perp - 4.
\end{equation}
Likewise, if $j \geq 2$ and
\begin{equation}\label{E29}
p > \frac{d_1}{2} + d_\perp - 2,
\end{equation}
the partial derivative $\partial_j^p R$ belongs to $L^q(B(0,1)^c)$ for
\begin{equation}\label{E30}
q > \frac{2N - 1}{2p + 4 - d_1 - 2 d_\perp},
\end{equation}
while it belongs to $L^\infty(B(0,1)^c)$ for
\begin{equation}\label{E31}
p \geq \frac{d_1}{2} + d_\perp - 2.
\end{equation}
\end{prop}

\begin{remark}
By linearity, similar estimates hold for any rational fractions of the form \eqref{E15} (where $P$ is any polynomial function on $\R^N$).
\end{remark}

\begin{remark}
We believe that Proposition \ref{P4} sharply describes the singularities of $\partial_j^p R$ near the origin and its integrability at infinity. We believe that this function does not belong to $M^\infty_\alpha(B(0,1))$ for any $\alpha > p + 2 - d$, while $\partial_1^p R$ and $\partial_j^p R$ do not belong to $L^\frac{2N - 1}{p + 4 - d_1 - 2d_\perp}(B(0,1)^c)$ and $L^\frac{2N - 1}{2p + 4 - d_1 - 2d_\perp}(B(0,1)^c)$ (except for $p = d_1 + 2 d_\perp -4$, respectively $p = \frac{d_1}{2} + d_\perp - 2$).
\end{remark}

Then, our argument links the behaviour of the tempered distribution $f$ to the behaviour of the derivatives of its Fourier transform $R$ given by Propositions \ref{P3} and \ref{P4}. This relies on explicit integral expressions which already appeared in \cite{Graveja3} and \cite{Graveja5}, and which are described in the next lemma.

\begin{lemma}\label{L1}
Let $f$ be a tempered distribution on $\R^N$ such that its Fourier transform belongs to $C^\infty(\R^N \setminus \{ 0 \})$. Let us assume moreover that there are some $1 \leq j \leq N$ and $p \in \N^*$ such that
\begin{itemize} 
\item[$(i)$] $\partial_j^p \widehat{f} \in L^1(B(0,1)^c)$,
\item[$(ii)$] $\partial_j^{p-1} \widehat{f} \in L^1(B(0,1))$,
\item[$(iii)$] $|.| \partial_j^p \widehat{f} \in L^1(B(0,1))$.
\end{itemize}
Then, the function $x \mapsto x_j^p f(x)$ is continuous on $\R^N$ and satisfies for every $\lambda > 0$,
\begin{equation}\label{E32}
\begin{split}
\forall x \in \R^N, x_j^p f(x) & = \frac{i^p}{(2 \pi)^N} \bigg( \int_{B(0,\lambda)^c} \partial_j^p \widehat{f}(\xi) e^{ix.\xi} d\xi + \frac{1}{\lambda} \int_{S(0,\lambda)} \xi_j \partial_j^{p-1} \widehat{f}(\xi) d\xi \\ & + \int_{B(0,\lambda)} \partial_j^p \widehat{f}(\xi) (e^{ix.\xi} - 1) d\xi \bigg).
\end{split}
\end{equation}
\end{lemma} 

Lemma \ref{L1} links the algebraic decay of a distribution $f$ at infinity, or its algebraic explosion near the origin, to the integrability of some derivatives of its Fourier transform $\widehat{f}$. Indeed, it is sufficient to prove that the right-hand side of equation \eqref{E32} is uniformly bounded on $B(0,1)$ (or $B(0,1)^c$) to infer that $f$ belongs to $M^\infty_p(B(0,1))$ (or $M^\infty_p(B(0,1)^c)$). In particular, Lemma \ref{L1} seems relevant to study the algebraic decay of the tempered distributions $f$ whose Fourier transforms $R$ satisfy \eqref{E15}-\eqref{E16}, or their singularity near the origin. Indeed, Proposition \ref{P4} yields a lot of information on their integrability. However, some derivatives of $R$ are not sufficiently integrable at infinity to satisfy the assumptions of Lemma \ref{L1}. In order to overcome this difficulty, we prove the next alternative of Lemma \ref{L1}.

\begin{lemma}\label{L2}
Let $1 \leq j \leq N$, $\lambda > 0$, and let $f$ be a tempered distribution on $\R^N$ such that its Fourier transform belongs to $C^\infty(\R^N \setminus \{ 0 \})$. Assume moreover that there exist some $1 \leq p \leq m$ and $A > 0$ such that
\begin{itemize} 
\item[$(i)$] $\forall \xi \in \R^N \setminus \{ 0 \}, |\widehat{f}(\xi)| \leq A (|\xi|^{-r} + |\xi|^s)$,
\item[$(ii)$] $\forall (k,\xi) \in \{ 0,\ldots,p \} \times B(0,1), |\xi|^{N-p+k} |\partial_j^k \widehat{f}(\xi)| \leq A$,
\item[$(iii)$] $\partial_j^m \widehat{f} \in L^1(B(0,1)^c)$,
\item[$(iv)$] $\forall k \in \{ 0,\ldots,m-1 \}, \partial_j^k \widehat{f} \in L^{q_{m-k}}(B(0,1)^c)$,
\end{itemize}
where $r < N$, $s \geq 0$, $1 < q_k < \frac{N}{N-k}$ if $1 \leq k \leq N-1$, and $1 < q_k \leq + \infty$ if $k > N$. Then, the function $x \mapsto x_j^p f(x)$ is continuous on $\Omega_j = \{x \in \R^N, x_j \neq 0 \}$, and satisfies for every $x \in \Omega_j$,
\begin{equation}\label{E33}
\begin{split}
x_j^p f(x) = & \frac{i^p}{(2 \pi)^N} \bigg( (-ix_j)^{p-m} \int_{B(0,\lambda)^c} \partial_j^m \widehat{f}(\xi) e^{ix.\xi} d\xi + \sum_{k=p}^{m-1} \frac{(-ix_j)^{p-k-1}}{\lambda} \int_{S(0,\lambda)} \xi_j \\ & \partial_j^k \widehat{f}(\xi) e^{ix.\xi} d\xi + \frac{1}{\lambda} \int_{S(0,\lambda)} \xi_j \partial_j^{p-1}
\widehat{f}(\xi) d\xi + \int_{B(0,\lambda)} \partial_j^p \widehat{f}(\xi) (e^{ix.\xi} - 1) d\xi \bigg).
\end{split}
\end{equation}
\end{lemma}

The assumptions of Lemma \ref{L2} are tailored for tempered distributions whose Fourier transforms are rational fractions of the form \eqref{E15}-\eqref{E16}. However, the proof of Lemma \ref{L2}, which results from Lemma \ref{L1} and some integrations by parts, can be adapted to other distributions which do not necessarily satisfy all the hypothesis of Lemma \ref{L2}.

By Proposition \ref{P4}, we now apply Lemma \ref{L2} to some
tempered distributions $f$ whose Fourier transforms are rational
fractions of the form \eqref{E15}-\eqref{E16}.

\begin{prop}\label{P5}
Let $j \in \{ 1,\ldots,N \}$, $\lambda > 0$, $\Omega_j = \{ x \in \R^N, x_j \neq 0\}$, and let $f$ be a tempered distribution on $\R^N$ whose Fourier transform $R$ is of the form \eqref{E15}-\eqref{E16}. Assume moreover that $d \neq 0$ if $N = 2$, and $d_1 + 2 d_\perp \leq 4$, where $d$ and $d_\perp$ are defined by \eqref{E17}. Then, the function $x \mapsto x_j^{p_j} f(x)$ is continuous on $\Omega_j$, and is given for every $x \in \Omega_j$ by
\begin{equation}\label{E34}
\begin{split}
x_j^{p_j} f(x) = & \frac{i^{p_j}}{(2 \pi)^N} \bigg( (-ix_j)^{p_j - m_j} \int_{B(0,\lambda)^c} \partial_j^{m_j} R(\xi) e^{ix.\xi} d\xi + \sum_{k=p_j}^{m_j-1} \frac{(-ix_j)^{p_j - k - 1}}{\lambda} \\ & \int_{S(0,\lambda)} \partial_j^k R(\xi) \xi_j e^{ix.\xi} d\xi + \frac{1}{\lambda} \int_{S(0,\lambda)} \xi_j \partial_j^{p_j - 1} R(\xi) d\xi + \int_{B(0,\lambda)} \partial_j^{p_j} R(\xi) \\ & (e^{ix.\xi} - 1) d\xi \bigg),
\end{split}
\end{equation}
where $p_j = N - 2 + d$, $m_1 = 2N - 4 + d_1 + 2d_\perp$, and if $j \geq 2$, $m_j = N - 2 + d$.
\end{prop}

\begin{remark}
Here, we make two additional assumptions. Indeed, if $d = 0$ and $N = 2$, all the derivatives of $R$ are not integrable near the origin. Therefore, we cannot expect to prove some formula like \eqref{E34} for $f$. On the other hand, the assumption $d_1 + 2 d_\perp \leq 4$ is only technical. Lemma \ref{L2} requires some integrability at infinity for the derivatives of $R$. This is not possible for the derivatives of low order if $d_1 + 2 d_\perp > 4$. However, we believe that Lemma \ref{L2} can be improved to compute a formula like \eqref{E34} even in the case $d_1 + 2 d_\perp > 4$. Since it is not useful in our context, we will not investigate this point any further.
\end{remark}

Formula \eqref{E34} links the algebraic decay of $f$ at infinity (or its algebraic explosion near the origin) to the integrability properties of $R$. Indeed, it suffices to choose $\lambda = \frac{1}{|x|}$, and to bound the right-hand side of formula \eqref{E34} by Proposition \ref{P4} to obtain the algebraic decay of $f$ stated by Proposition \ref{Paumond1}. Likewise, we choose $\lambda = \frac{1}{|x_j|^{1 + \delta_{j,1}}}$, and bound the right-hand side of formula \eqref{E34} by Proposition \ref{P4} to establish the algebraic explosion of $f$ near the origin. Thus, the main interest of formula \eqref{E34} is the possibility of choosing $\lambda$ in the most fruitful way. In particular, we can adapt the value of $\lambda$ to the anisotropy of our problem to obtain the anisotropic estimates \eqref{E18} of Proposition \ref{P2}.

Finally, Theorems \ref{T3} and \ref{T4}, and Corollary \ref{C1} are direct consequences of Propositions \ref{Paumond1} and \ref{P2}. They correspond to the cases $(d_1,d_\perp) = (1,0)$ for $H_0$, $(d_1,d_\perp) = (2,0)$ for $K_0$ and $(d_1,d_\perp) = (2+\delta_{k,1}, 1-\delta_{k,1})$ for $K_k$.

\subsubsection*{Pointwise limit of the kernel $K_0$ at infinity}

The integral expressions of Lemmas \ref{L1} and \ref{L2} have another important application: the computation of the pointwise limit at infinity of $K_0$, that is the limit when $R$ tends to $+ \infty$ of the function
$$R \mapsto R^N K_0(R \sigma - y),$$
for every $\sigma \in \S^{N-1}$ and $y \in \R^N$.

\begin{theorem}\label{T5}
Let $\sigma \in \S^{N-1}$ and $y \in \R^N$. Then,
\begin{equation}\label{E35}
R^N K_0(R \sigma - y) \underset{R \to + \infty}{\to}
\frac{\Gamma(\frac{N}{2})}{2 \pi^\frac{N}{2}} (1 - N \sigma_1^2).
\end{equation}
\end{theorem}

Theorem \ref{T5} results from formula \eqref{E34}. Indeed, by
Proposition \ref{P5}, the kernel $K_0$ satisfies for every $j \in \{
1,\ldots,N\}$ and $x \in \Omega_j = \{ x \in \R^N, x_j \neq 0\}$,
\begin{equation}\label{E36}
\begin{split}
x_j^N K_0(x) = & \frac{i^N}{(2 \pi)^N} \bigg( (-ix_j)^{N - m_j}
\int_{B(0,\lambda)^c} \partial_j^{m_j} \widehat{K_0}(\xi) e^{ix.\xi}
d\xi + \sum_{k=N}^{m_j-1} \frac{(-ix_j)^{N - k - 1}}{\lambda}
\\ & \int_{S(0,\lambda)} \xi_j \partial_j^k \widehat{K_0}(\xi)
e^{ix.\xi} d\xi + \frac{1}{\lambda} \int_{S(0,\lambda)} \xi_j
\partial_j^{N-1} \widehat{K_0}(\xi) d\xi + \int_{B(0,\lambda)}
\partial_j^N \widehat{K_0}(\xi) \\ & (e^{ix.\xi} - 1) d\xi \bigg),
\end{split}
\end{equation}
where $m_1 = 2N - 2$ and $m_j = N$ if $j \geq 2$. By applying the
dominated convergence theorem to this formula, we obtain for any $j
\in \{ 1,\ldots,N \}$ such that $\sigma_j \neq 0$,
\begin{equation}\label{E37}
\begin{split}
R^N K_0(R \sigma - y) \underset{R \to + \infty}{\to} & \frac{i^N}{(2
 \pi \sigma_j)^N} \bigg( \frac{i}{\sigma_j} \int_{B(0,1)^c}
\partial_j^{N+1} \widehat{R_{1,1}}(\xi) e^{i \sigma. \xi} d\xi +
\frac{i}{\sigma_j} \int_{\S^{N-1}} \xi_j e^{i \sigma. \xi} \\ & \partial_j^N
\widehat{R_{1,1}}(\xi)  d\xi + \int_{\S^{N-1}}
\xi_j \partial_j^{N-1} \widehat{R_{1,1}}(\xi) d\xi + \int_{B(0,1)}
\partial_j^N \widehat{R_{1,1}}(\xi) \\ & (e^{i \sigma. \xi} - 1) d\xi
\bigg).
\end{split}
\end{equation}
Here, the distribution $R_{1,1}$ is the so-called composed Riesz kernel
given by
\begin{equation}\label{E38}
\widehat{R_{1,1}}(\xi) = \frac{\xi_1^2}{|\xi|^2}.
\end{equation}
To complete the proof of Theorem \ref{T5}, it now remains to prove
that the right-hand sides of equations \eqref{E35} and \eqref{E37} are equal.

\begin{theorem}\label{T6}
Let $1 \leq j \leq N$, $y \in \R^N$ and $\sigma \in \S^{N-1}$ such
that $\sigma_j \neq 0$. Then,
\begin{equation}\label{E39}
\begin{split}
\frac{\Gamma(\frac{N}{2})}{2 \pi^\frac{N}{2}} (1 - N \sigma_1^2) = & \frac{i^N}{(2 \pi \sigma_j)^N} \bigg( \frac{i}{\sigma_j} \int_{B(0,1)^c} \partial_j^{N+1} \widehat{R_{1,1}}(\xi) e^{i \sigma. \xi} d\xi + \frac{i}{\sigma_j} \int_{\S^{N-1}} \xi_j e^{i \sigma. \xi} \\ & \partial_j^N \widehat{R_{1,1}}(\xi) d\xi + \int_{\S^{N-1}} \xi_j \partial_j^{N-1} \widehat{R_{1,1}}(\xi) d\xi + \int_{B(0,1)} \partial_j^N \widehat{R_{1,1}}(\xi) \\ & (e^{i \sigma. \xi} - 1) d\xi \bigg).
\end{split}
\end{equation}
\end{theorem}

Theorem \ref{T6} follows once more from Lemma \ref{L2}, and from an
explicit formula for $R_{1,1}$. Indeed, by standard Riesz operator
theory, the composed Riesz kernel $R_{1,1}$ is equal to
\begin{equation}\label{E40}
R_{1,1}(x) = \frac{\Gamma(\frac{N}{2})}{2 \pi^\frac{N}{2}} \bigg( PV
\bigg( \frac{|x|^2 - N x_1^2}{|x|^{N+2}} 1_{B(0,1)} \bigg) +
\frac{|x|^2 - N x_1^2}{|x|^{N+2}} 1_{B(0,1)^c} \bigg) +
\frac{\delta_0}{N},
\end{equation}
where $\delta_0$ is the Dirac mass at the origin, and $PV \bigg( \frac{|x|^2 - N x_1^2}{|x|^{N+2}} 1_{B(0,1)} \bigg)$ is the principal value of the function $x \mapsto \frac{|x|^2 - N x_1^2}{|x|^{N+2}} 1_{B(0,1)}(x)$ at the origin, which is defined by
\begin{equation}\label{E41}
\forall \phi \in C_0^\infty(\R^N), \Big< PV \bigg( \frac{|x|^2 - N
 x_1^2}{|x|^{N+2}} 1_{B(0,1)} \bigg), \phi \Big> = \int_{B(0,1)}
\frac{|x|^2 - N x_1^2}{|x|^{N+2}} (\phi(x) - \phi(0)) dx.
\end{equation}
Then, formula \eqref{E39} is a consequence of expression \eqref{E40}.

The main interest of Theorem \ref{T5} lies in the computation of the
first order term $v_\infty$ of the asymptotics of $v$. Indeed, the
computation of $v_\infty$ will result in subsection 2.5 from the
computation of the pointwise limit of $K_0$. However, formulae \eqref{E7} and \eqref{E35} already emphasise the link between the pointwise limit of $K_0$ and $v_\infty$.

\subsubsection*{Rigorous formulation of equations \eqref{E10} and \eqref{E11}}

We conclude our study of $H_0$, $K_0$ and $K_k$ by giving a rigorous sense to the convolution equations \eqref{E10} and \eqref{E11}.

Indeed, let us consider equation \eqref{E10}. By Theorem \ref{T3} and Corollary \ref{C1}, the kernel $H_0$ belongs to $L^q(\R^N)$ for every $\frac{N}{N-1} < q < \frac{2N - 1}{2N - 4}$. However, we will prove in Theorem \ref{T8} that the functions $v$ and $\nabla v$ belong to $L^q(\R^N)$ for every $1 < q < + \infty$. Therefore, the function $v^p \partial_1 v$ belongs to $L^1(\R^N)$. Thus, by Young's inequalities, equation \eqref{E10} makes sense in $L^q(\R^N)$ for every $\frac{N}{N-1} < q < \frac{2N - 1}{2N - 4}$. In particular, it makes sense almost everywhere, which will be sufficient in the following.

Likewise, by Theorem \ref{T3} and Corollary \ref{C1}, the kernel $K_0$ belongs to $L^q(\R^N)$ for every $1 < q < \frac{2N - 1}{2N - 3}$, and by Theorem \ref{T8}, the function $v^{p+1}$ belongs to $L^1(\R^N)$. Thus, equation \eqref{E11} makes sense in $L^q(\R^N)$ for every $1 < q < \frac{2N - 1}{2N - 3}$, and consequently, almost everywhere.

However, we will also consider the gradient of $v$. In particular, we will take the gradient of equation \eqref{E11}, which leads to another difficulty. Indeed, by Theorem \ref{T4}, the first order derivatives of $K_0$, given up to a multiplicative coefficient by $K_k$, have non-integrable singularities near the origin. We are not allowed to differentiate the convolution equation \eqref{E11} without additional care: we cannot write
$$\partial_k v = \frac{1}{p+1} \partial_k K_0 * v^{p+1}.$$
Our method to overcome this difficulty is reminiscent of some classical arguments in distribution theory, using integral formulae. Indeed, by Theorem \ref{T3} and Corollary \ref{C1}, the kernel $K_0$ has first order partial derivatives in the sense of distributions, which are equal to
$$\partial_k K_0 = i \ K_k 1_{B(0,1)^c} + i \ PV(K_k 1_{B(0,1)}) + \bigg( \int_{\S^{N-1}} K_0(y) y_k dy \bigg) \delta_0,$$
where $PV(K_k 1_{B(0,1)})$ denotes (as above for the composed Riesz kernels) the principal value at the origin of $K_k$
$$\forall \phi \in C_0^\infty(B(0,1)), <PV(K_k 1_{B(0,1)}), \phi> = \int_{B(0,1)} K_k(x) (\phi(x) - \phi(0)) dx.$$
Thus, provided we know sufficient smoothness for $v^{p+1}$, we will be able to take the derivative of equation \eqref{E11} and to obtain an explicit integral expression for its derivative.

\begin{lemma}\label{L3}
Let us consider a function $f \in C^0(\R^N)$ such that
\begin{itemize}
\item[(i)] $f \in L^\infty(\R^N) \cap M^\infty_{N(p+1)}(\R^N)$,
\item[(ii)] $\nabla f \in L^\infty(\R^N)^N$,
\end{itemize}
and let $g = K_0*f$. Then, $g$ is of class $C^1$ on $\R^N$. Moreover, its partial derivative $\partial_k g$ is given for every $x \in \R^N$ by
\begin{equation}\label{E42}
\begin{split}
\partial_k g(x) & = i \int_{B(0,1)^c} K_k(y) f(x-y) dy + i \int_{B(0,1)} K_k(y) (f(x-y) - f(x)) dy \\ & + \bigg( \int_{\S^{N-1}}  K_0(y) y_k dy \bigg) f(x).
\end{split}
\end{equation}
\end{lemma}

The function $v^{p+1}$ will satisfy the assumptions of Lemma \ref{L3}, so we will be able to give a sense to the gradient of equation \eqref{E11} and to complete our asymptotic analysis of $v$.

\subsubsection{Decay properties of the solitary waves}

In this section, we establish the algebraic decay of a solitary wave $v$ and of its gradient. In particular, we extend to any dimension $N \geq 2$ the decay properties stated by A. de Bouard and J.-C. Saut \cite{deBoSau2} in dimension $N = 2$.

\begin{theorem}\label{T7}
Let $v \in Y$ be a solitary-wave solution of speed $1$ of equation \eqref{E1}. Assume that $0 < p < \frac{4}{2N - 3}$. Then, $v$ belongs to $M^\infty_N(\R^N)$, while its gradient $\nabla v$ belongs to $M^\infty_{\min \{(p+1)N,N+1\} }(\R^N)$.
\end{theorem}

Theorem \ref{T7} relies on a standard argument which is reminiscent of a series of papers by J.L. Bona and Yi A. Li \cite{BonaLi2}, A. de Bouard and J.-C. Saut \cite{deBoSau2}, and M. Maris \cite{Maris1} \cite{Maris2} (see also \cite{Graveja3} \cite{Graveja5}). This argument links the algebraic decay of the solitary waves to the algebraic decay of the kernels of the convolution equations they satisfy. These algebraic decays are identical, because the nonlinearity $v^{p+1}$ is superlinear.

To get a feeling for this claim, let us consider the simplified model
$$f = K*f^q,$$
where we assume that $q > 1$, $f$ and $K$ are smooth functions such that $f$ belongs to $L^q(\R^N)$ for every $1 < q \leq + \infty$, and the kernel $K$ belongs to $L^1(\R^N)$ and to $M^\infty_{\alpha_K}(\R^N)$ for some $\alpha_K > 0$. In order to study the algebraic decay of $f$, we write for every $x \in \R^N$ and $\alpha > 0$,
\begin{equation}\label{E43}
\begin{split}
|x|^\alpha |f(x)| \leq & A \bigg( \int_{\R^N} |x-y|^\alpha |K(x-y)| |f(y)|^q dy + \int_{\R^N} |K(x-y)| |y|^\alpha |f(y)|^q dy \bigg) \\ \leq & A \big( \| K \|_{M^\infty_\alpha(\R^N)} \| f \|^q_{L^q(\R^N)} + \| K \|_{L^1(\R^N)} \| f \|^q_{M^\infty_{\frac{\alpha}{q}}(\R^N)} \big).
\end{split}
\end{equation}
The function $f$ now belongs to $L^q(\R^N)$, while the kernel $K$ belongs to $L^1(\R^N)$ and to $M^\infty_\alpha(\R^N)$, provided $0 \leq \alpha \leq \alpha_K$. Therefore, if $0 \leq \alpha \leq \alpha_K$, equation \eqref{E43} reduces to
\begin{equation}\label{E44}
\| f \|_{M^\infty_\alpha(\R^N)} \leq A + A \| f \|^q_{M^\infty_{\frac{\alpha}{q}}(\R^N)}.
\end{equation}
Equation \eqref{E44} links the algebraic decay with exponent $\alpha$ of $f$ to its algebraic decay with exponent $\frac{\alpha}{q}$. In particular, if we know some algebraic decay with a small exponent $\alpha_0 > 0$, a bootstrap argument yields that $f$ belongs to $M^\infty_\alpha(\R^N)$ for $\alpha = q \alpha_0$, $\alpha = q^2 \alpha_0$, $\ldots$, that is for every $\alpha \in [0, \alpha_K]$. This provides a striking optimal decay property for superlinear equations. Indeed, assuming $f$ possesses some algebraic decay, then, $f$ decays as fast as the kernel. However, some decay of $f$ must be established first in order to initiate the inductive argument.

The situation is more involved for the function $v$ and for the kernels $H_0$, $K_0$ and $K_k$. The main difficulties come from the singularities near the origin of the kernels. In particular, in the case of $K_k$, we must adapt our argument to equation \eqref{E42}, which is no longer a convolution equation.

However, in order to perform the argument above, we first need to determine some integrability and initial decay for $v$. Indeed, we already know the integrability and the algebraic decay of $H_0$, $K_0$ and $K_k$ by Theorems \ref{T3} and \ref{T4}, and Corollary \ref{C1}. The integrability properties of $v$ follow from the following theorem due to A. de Bouard and J.-C. Saut \cite{deBoSau1, deBoSau2} in dimensions $N=2$ and $N=3$
\footnote{The proof of A. de Bouard and J.-C. Saut is still relevant in any dimension $N \geq 4$. Therefore, we will omit the proof of Theorem \ref{T8}, and refer to \cite{Graveja0} for a detailed proof in any dimension $N \geq 4$.}.

\begin{theorem}[\cite{deBoSau1,deBoSau2}]\label{T8}
Let $v \in Y$ be a solitary-wave solution of speed $1$ of equation \eqref{E1}. Assume that $0 < p < \frac{4}{2N - 3}$. Then, $v$ is bounded and continuous on $\R^N$. Moreover, the functions $v$, $\nabla v$ and $\partial_1^2 v$ belong to $L^q(\R^N)$ for every $q \in ]1, + \infty[$.
\end{theorem}

\begin{remark}
Another interest of Theorem \ref{T8} is to give a rigorous sense to equations \eqref{E10}, \eqref{E11} and \eqref{E42} (which was already mentioned in subsection 1.2.2).
\end{remark}

It now only remains to compute some initial decay for $v$. We deduce it from the next proposition due to A. de Bouard and J.-C. Saut \cite{deBoSau2}
\footnote{A. de Bouard and J.-C. Saut \cite{deBoSau2} proved Proposition \ref{P6} in dimensions two and three. However, their proof is still relevant in any dimension $N \geq 4$. Therefore, we will omit the proof of Proposition \ref{P6}, and refer to \cite{deBoSau2} for more details.}.

\begin{prop}[\cite{deBoSau2}]\label{P6}
Let $v \in Y$ be a solitary-wave solution of speed $1$ of equation \eqref{E1}. Assume that $0 < p < \frac{4}{2N - 3}$. Then,
\begin{equation}\label{E45}
\int_{\R^N} |x|^2 (|\nabla v(x)|^2 + |\partial^2_1 v(x)|^2) dx < + \infty.
\end{equation}
\end{prop}

Proposition \ref{P6} gives very weak decay for the gradient of $v$ and its second order partial derivative $\partial_1^2 v$. In particular, it does not provide local algebraic decay, but integral algebraic decay. However, it is sufficient to apply the inductive argument above and to prove Theorem \ref{T7}.

\subsubsection{Asymptotics of the solitary waves}

Theorem \ref{T7} establishes the decay properties of a solitary wave $v$ and of its gradient. In order to complete the proof of Theorem \ref{T1}, we now compute the first order asymptotics of $v$, i.e. the limit when $|x|$ tends to $+ \infty$ of the function $x \mapsto |x|^N v(x)$. Our argument is reminiscent of \cite{Graveja5} and relies once more on the convolution equation \eqref{E11}.

Indeed, by equation \eqref{E11}, this computation reduces to compute the limit when $R$ tends to $+ \infty$ of the functions $v_R$ defined by
\begin{equation}\label{E46}
\forall \sigma \in \S^{N-1}, v_R(\sigma) = R^N v(R \sigma) = \frac{R^N}{p+1} \int_{\R^N} K_0(R\sigma - y) v^{p+1}(y) dy,
\end{equation}
and to prove that this convergence is uniform on the sphere $\S^{N-1}$, provided $p \geq \frac{1}{N}$. Here, we first invoke the dominated convergence theorem for formula \eqref{E46} to compute the pointwise limit of $v_R$.

\begin{prop}\label{P7}
Let $v \in Y$ be a solitary-wave solution of speed $1$ of equation \eqref{E1}. Assume that $0 < p < \frac{4}{2N - 3}$. Then,
\begin{equation}\label{E47}
\forall \sigma \in \S^{N-1}, v_R(\sigma) \underset{R \to +  \infty}{\to} v_\infty(\sigma) = \frac{\Gamma(\frac{N}{2})}{2 \pi^\frac{N}{2} (p+1)} (1 - N \sigma_1^2) \int_{\R^N} v^{p+1}(y) dy.
\end{equation}
\end{prop}

Indeed, by Theorem \ref{T5}, the integrand of formula \eqref{E46} satisfies for every $y \in \R^N$ and $\sigma \in \S^{N-1}$,
$$R^N K_0(R\sigma - y) v^{p+1}(y) \underset{R \to + \infty}{\to} \frac{\Gamma(\frac{N}{2})}{2 \pi^\frac{N}{2} (p+1)} (1 - N \sigma_1^2) v^{p+1}(y).$$
Then, it remains to dominate this integrand by means of Theorems \ref{T3}, \ref{T4} and \ref{T8} to obtain Proposition \ref{P7}.

Then, we apply Ascoli-Arzela's theorem to prove the uniformity of the convergence. 

\begin{prop}\label{P8}
Let $v \in Y$ be a solitary-wave solution of speed $1$ of equation \eqref{E1}. Assume that $\frac{1}{N} \leq p < \frac{4}{2N - 3}$. Then,
\begin{equation}\label{E48}
\| v_R - v_\infty \|_{L^\infty(\S^{N-1})} \underset{R \to +  \infty}{\to} 0.
\end{equation}
\end{prop}

Indeed, we already know the existence of a pointwise limit at infinity, so Ascoli-Arzela's theorem gives the uniformity of the convergence. However, this theorem requires some compactness: we deduce it from the algebraic decay of the gradient of $v$. Provided $p \geq \frac{1}{N}$, this gradient belongs to
$M^\infty_{N+1}(\R^N)$ by Theorem \ref{T7}. Therefore, the gradients on the sphere $\S^{N-1}$ of the functions $v_R$ are uniformly bounded on $\S^{N-1}$. This yields the compactness result necessary to apply Ascoli-Arzela's theorem, and to prove assertion \eqref{E48}.

\begin{remark}
Our argument based on Ascoli-Arzela's theorem fails to prove the uniformity of the convergence when $0 < p < \frac{1}{N}$. However, the nonlinearity $v^{p+1}$ is less and less smooth at the origin when $p$ tends to $0$. Therefore, the convergence may be no longer uniform when $p$ is too small.
\end{remark}

This concludes the proof of Theorem \ref{T1}, which is a direct consequence of Theorem \ref{T7}, and Propositions \ref{P7} and \ref{P8}.

\subsection{Plan of the paper}

The paper splits into two parts. The first part is devoted to the analysis of $H_0$, $K_0$ and $K_k$. In the first section, we establish some properties of the rational fractions of the form \eqref{E15}-\eqref{E16} stated in Propositions \ref{P3} and \ref{P4}. In the second one, we prove Lemmas \ref{L1} and \ref{L2} to obtain the integral formula \eqref{E34}: it is the starting point of the proofs of Propositions \ref{Paumond1} and \ref{P2}. The third section deals with the algebraic decay of $H_0$, $K_0$ and $K_k$ at infinity, and their explosion near the origin, stated in Theorems \ref{T3} and \ref{T4}, and Corollary \ref{C1}. The fourth section is concerned with the pointwise limit of $K_0$ obtained by Theorems \ref{T5} and \ref{T6}, while in the last section, we prove Lemma \ref{L3} to give a rigorous sense to the derivative of the convolution equation \eqref{E11}.

The proof of Theorem \ref{T1} forms the core of the second part. The first section is devoted to the proof of Theorem \ref{T7}, which gives the optimal algebraic decay of a solitary wave and of its gradient. In the second one, we conclude the proof of Theorem \ref{T1} by computing the asymptotics of a solitary wave in Propositions \ref{P7} and \ref{P8}. In the last section, we focus on the standard Kadomtsev-Petviashvili equation. We link the asymptotics of a solitary wave to its energy and its action by proving Theorem \ref{T2}. As mentioned above, Theorem \ref{T2} follows from the standard Pohozaev identities derived by A. de Bouard and J.-C. Saut in \cite{deBoSau1}
\footnote{A. de Bouard and J.-C. Saut \cite{deBoSau1} actually proved Lemma \ref{L4} in dimensions two and three. However, their proof is still relevant in dimensions $N \geq 4$, so we refer to \cite{deBoSau1} for the proof of Lemma \ref{L4}.}.

\begin{lemma}[\cite{deBoSau1}]\label{L4}
Let us consider some positive real number $p$ and a solitary-wave solution $v \in Y$ of speed $1$ of equation \eqref{E1}. Then, the following identities hold for every $k \in \{ 2,\ldots,N \}$:
\begin{align}
& \int_{\R^N} \Big( - v(x)^2 + \frac{2}{p+2} v(x)^{p+2} - 3 \partial_1
 v(x)^2 + \sum_{j=2}^N v_j(x)^2 \Big) dx = 0, \label{E49} \\ &
\int_{\R^N} \Big( v(x)^2 - \frac{2}{(p+1)(p+2)} v(x)^{p+2} +
 \partial_1 v(x)^2 - 2 v_k(x)^2 + \sum_{j=2}^N v_j(x)^2 \Big) dx =
0, \label{E50} \\ & \int_{\R^N} \Big( v(x)^2 - \frac{1}{p+1}
 v(x)^{p+2} + \partial_1 v(x)^2 + \sum_{j=2}^N v_j(x)^2 \Big) dx =
0. \label{E51}
\end{align}
\end{lemma}

Finally, we mention another consequence of Lemma \ref{L4}: the non-existence of non-trivial solutions of equation \eqref{E4} in $Y$ when $p \geq \frac{4}{2N - 3}$ in any dimension $N \geq 4$
\footnote{A. de Bouard and J.-C. Saut \cite{deBoSau1} already proved their non-existence in dimensions two and three.}.

\begin{cor}\label{C2}
Let us consider a solitary-wave solution $v \in Y$ of speed $1$ of equation \eqref{E1} and, assume that $N \geq 4$ and $p \geq \frac{4}{2N - 3}$. Then, $v$ is constant.
\end{cor}

\section{Main properties of the kernels $H_0$, $K_0$ and $K_k$}

In this part, we state some properties of $H_0$, $K_0$ and $K_k$. We first study their algebraic decay at infinity and their explosion near the origin. These follow from the integrability properties of their Fourier transforms using some integral expressions. Then, we compute the pointwise limit at infinity of $K_0$, before deriving rigorously equations \eqref{E10}, \eqref{E11} and \eqref{E42}.

\subsection{Properties of the Fourier transforms of the kernels $H_0$, $K_0$ and $K_k$}

This first section is devoted to the analysis of $\widehat{H_0}$, $\widehat{K_0}$ and $\widehat{K_k}$. By formulae \eqref{E12}, \eqref{E13} and \eqref{E14}, they are rational fractions of the form \eqref{E15}-\eqref{E16}. In Proposition \ref{P3}, we describe the derivatives of all the rational fractions of this form. Here, the main difficulty comes from their anisotropy in the direction $\xi_1$: we must distinguish the terms including the variable $\xi_1$ and the other terms. More precisely, Proposition \ref{P3} results from the next inductive argument.

\def\proof{\par{\it Proof of Proposition \ref{P3}}. \ignorespaces}
\begin{proof}
By induction, the derivatives $\partial_j^p R$ satisfy equation \eqref{E21}. Here, $P_{j,p}$ are polynomial functions on $\R^N$ given by
\begin{equation}\label{E52}
\begin{split}
P_{j,0}(\xi) = & P(\xi) = \underset{j=1}{\overset{N}{\Pi}} \xi_j^{d_j}, \\ P_{j,p+1}(\xi) = & (|\xi|^2 + \xi_1^4) \partial_j P_{j,p}(\xi) - 2 (p+1) (\xi_j + 2 \delta_{j,1} \xi_1^3) P_{j,p}(\xi).
\end{split}
\end{equation}
In particular, the inductive definitions of $P_{j,p}$ are different according to $j=1$ or $j \geq 2$. Thus, we split our analysis into two cases depending on the value of $j$.

{\it Case $j=1$.} The polynomial function $P_{1,p}$ is given by
\begin{equation}\label{E53}
P_{1,p}(\xi) = \sum_{k=0}^{+ \infty} a_{k,p}(\xi_1) |\xi_\perp|^{2k}
\underset{j=2}{\overset{N}{\Pi}} \xi_j^{d_j},
\end{equation}
where the functions $a_{k,p}$ are polynomial functions on $\R$. Indeed, by formulae \eqref{E52}, if the function $P_{1,p}$ is of the
form above, the function $P_{1,p+1}$ is given by
$$P_{1,p+1}(\xi) = \sum_{k=0}^{+ \infty} \Big( (|\xi_\perp|^2 +
\xi_1^2 + \xi_1^4) a_{k,p}'(\xi_1) - 2 (p+1) (\xi_1 + 2 \xi_1^3)
a_{k,p}(\xi_1) \Big) |\xi_\perp|^{2k} \underset{j=2}{\overset{N}{\Pi}}
\xi_j^{d_j}.$$
Therefore, with the usual convention $a_{-1,p} = 0$, the functions
$a_{k,p}$ are inductively given on $k \in \N$ by
\begin{align*}
 a_{k,0}(\xi_1) = & \delta_{k,0} \xi_1^{d_1}, \\ a_{k,p+1}(\xi_1) = &
 a_{k-1,p}'(\xi_1) + (\xi_1^2 + \xi_1^4) a_{k,p}'(\xi_1) - 2 (p+1)
 \xi_1 (1 + 2 \xi_1^2) a_{k,p}(\xi_1).
\end{align*}
Thus, they are polynomial functions on $\R$. Moreover, letting $d_{k,p} = {\rm deg}(a_{k,p}) \in \N \cup \{ - \infty \}$ be their degree (with the usual convention ${\rm deg}(0) = - \infty$), we obtain
$$d_{k,p+1} \leq \max \{ d_{k-1,p} - 1, d_{k,p} + 3 \}.$$
It follows by induction that
$$d_{k,p} \leq d_1 + 3p - 4k.$$
Therefore, we deduce from equation \eqref{E53} that
$$\forall \xi \in B(0,1)^c, |P_{1,p}(\xi)| \leq A_p \sum_{0 \leq k
 \leq \frac{d_1 + 3p}{4}} \max \{ |\xi_1|, 1 \}^{d_1 + 3p - 4k}
|\xi_\perp|^{2k + d_\perp},$$
which leads to inequality \eqref{E22}.

On the other hand, by equations \eqref{E52} and a straightforward
inductive argument, either the function $P_{1,p}$ is identically equal to $0$, or its terms of lowest degree are of degree $v_{1,p} \geq d_1 + d_\perp + p$. In particular, it follows that
$$\forall \xi \in B(0,1), |P_{1,p}(\xi)| \leq A_p |\xi|^{v_{1,p}} \leq A_p |\xi|^{d_1 + d_\perp + p},$$
which is exactly inequality \eqref{E23}.

{\it Case $j \geq 2$.} The polynomial function $P_{j,p}$ is given by
\begin{equation}\label{E54}
P_{j,p}(\xi) = \sum_{k=0}^{+ \infty} b_{k,p}(\xi_\perp)
\xi_1^{2k+d_1},
\end{equation}
where the functions $b_{k,p}$ are polynomial functions on
$\R^{N-1}$. Indeed, by formulae \eqref{E52}, if the function $P_{j,p}$ is of the form above, the function $P_{j,p+1}$ is given by
$$P_{j,p+1}(\xi) = \sum_{k=0}^{+ \infty} \bigg( \partial_j
b_{k,p}(\xi_\perp) \Big( \xi_1^4 + \xi_1^2 + |\xi_\perp|^2 \Big) - 2
(p+1) \xi_j b_{k,p}(\xi_\perp) \bigg) \xi_1^{2k + d_1}.$$
Hence, with the usual convention $b_{-1,p} = 0$, the functions $b_{k,p}$ are inductively given on $k \in \N$ by
\begin{align*}
b_{k,0}(\xi_\perp) = & \delta_{k,0} \underset{j=2}{\overset{N}{\Pi}}
\xi_j^{d_j}, \\ b_{k,p+1}(\xi_\perp) = & \partial_j
b_{k-2,p}(\xi_\perp) + \partial_j b_{k-1,p}(\xi_\perp) + |\xi_\perp|^2
\partial_j b_{k,p}(\xi_\perp) - 2 (p+1) \xi_j b_{k,p}(\xi_\perp).
\end{align*}
Thus, they are polynomial functions on $\R^{N-1}$. Moreover, letting $d_{k,p}' = {\rm deg}(b_{k,p}) \in \N \cup \{ - \infty
\}$ be their degree, we obtain
$$d_{k,p+1}' \leq \max \{ d_{k-2,p}' - 1, d_{k-1,p}' - 1, d_{k,p}' + 1 \}.$$
It follows by induction that
$$d_{k,p}' \leq d_\perp + p - k - \nu(k),$$
where $\nu(k) = 0$ if $k$ is even, and $\nu(k) = 1$ if $k$ is
odd. Therefore, equation \eqref{E54} leads to
\begin{align*}
|P_{j,p}(\xi)| \leq & A_p \bigg( \sum_{0 \leq 2k \leq d_\perp + p}
\max \{ |\xi_\perp|, 1 \}^{d_\perp - p - 2k} |\xi_1|^{4k + d_1} \\ + &
\sum_{0 \leq 2k + 1 \leq d_\perp + p} |\xi_1|^{4k + 2 + d_1} \max \{
|\xi_\perp|, 1 \}^{d_\perp - p - 2k - 2} \bigg) \\ \leq & A_p \sum_{0
 \leq k \leq \frac{d_\perp + p}{2}} \max \{ |\xi_\perp|, 1
\}^{d_\perp - p - 2k} \max\{ |\xi_1|, 1 \}^{4k + d_1},
\end{align*}
for every $\xi \in B(0,1)^c$, which is inequality \eqref{E24}.

On the other hand, by equations \eqref{E52} and a straightforward
inductive argument, either the function $P_{j,p}$ is identically equal to $0$, or its terms of lowest degree are of degree $v_{j,p} \geq d_1 + d_\perp + p$. In particular, it follows that
$$\forall \xi \in B(0,1), |P_{j,p}(\xi)| \leq A_p |\xi|^{v_{j,p}} \leq A_p |\xi|^{d_1 + d_\perp + p},$$
which completes the proof of Proposition \ref{P3}.
\end{proof}

Proposition \ref{P3} yields estimates of the derivatives of the rational fraction $R$, which are sufficiently sharp to describe the singularities of its derivatives near the origin, and their integrability at infinity in Proposition \ref{P4}.

\def\proof{\par{\it Proof of Proposition \ref{P4}}. \ignorespaces}
\begin{proof}
Let us first consider the behaviour of $\partial_j^p R$ near the origin. By equation \eqref{E21}, and estimates \eqref{E23} and \eqref{E25} of Proposition \ref{P3}, we have
$$|\partial_j^p R (\xi)| \leq A_p \frac{|\xi|^{p + d}}{(|\xi|^2 + \xi_1^4)^{p+1}} \leq \frac{A_p}{|\xi|^{p + 2 - d}},$$
for every $\xi \in B(0,1)$. Therefore, the function $\partial_j^p R$ belongs to $M^\infty_{p + 2 - d}(B(0,1))$.

Then, let us consider the integrability properties of $\partial_j^p R$ at infinity. By equation \eqref{E21}, and inequalities \eqref{E22} and \eqref{E24} of Proposition \ref{P3}, we can estimate the $L^q$-norm of $\partial_j^p R$ for every $q \in [1, + \infty]$. However, estimates \eqref{E22} and \eqref{E24} differ because of the anisotropy of $R$. Thus, we split our study into two cases depending on the value of $j$.

{\it Case $j=1$.} We first assume that $q < + \infty$. Then, from
Proposition \ref{P3},
\begin{equation}\label{E55}
\int_{B(0,1)^c} |\partial_1^p R(\xi)|^q d\xi \leq A \sum_{0 \leq k \leq \frac{d_1 + 3p}{4}} \int_1^{+ \infty} J_{q(d_1 + 3p - 4k), q(2k + d_\perp), q(p+1)}(r) dr,
\end{equation}
where we denote for every $\lambda \in [1, + \infty[$ and $(\alpha, \beta, \gamma) \in (\R_+^*)^3$,
\begin{equation}\label{E56}
J_{\alpha, \beta, \gamma}(\lambda) = \int_{S(0, \lambda)} \frac{\max
 \{|\xi_1| , 1 \}^{\alpha} \max \{ 1, |\xi_\perp| \}^\beta}{(|\xi|^2
 + \xi_1^4)^\gamma} d\xi.
\end{equation}
However, by using spherical coordinates $x = (r \cos(\theta_1),
\ldots, r \sin(\theta_1) \ldots \sin(\theta_ {N-1}))$, and the
successive changes of variables $u = \tan(\theta_1)$ and $v =
\frac{u}{\sqrt{\lambda}}$, we compute
\begin{align*}
J_{\alpha, \beta, \gamma}(\lambda) \leq & A \lambda^{N - 1 + \alpha +
 \beta - 2 \gamma} \int_0^{\frac{\pi}{2}} \frac{\max
 \{\cos(\theta_1), \frac{1}{\lambda} \}^\alpha}{(1 + \lambda^2
 \cos(\theta_1)^4)^\gamma} d\theta_1
\\ \leq & A \lambda^{N - 1 + \beta - 2 \gamma} \bigg(
\int_{\cos(\theta_1) \leq \frac{1}{\lambda}} d\theta_1 +
\lambda^\alpha \int_0^{+ \infty} \frac{(1 + u^2)^{2 \gamma -
 \frac{\alpha}{2} - 1}}{((1 + u^2)^2 + \lambda^2)^\gamma} du \bigg)
\\ \leq & A \lambda^{N - 2 + \beta - 2 \gamma} \bigg( 1 + \lambda^{ 1
 + \alpha} \bigg(\int_0^1 \frac{du}{(1 + \lambda^2)^\gamma} +
\lambda^{-\frac{\alpha + 1}{2}} \int_\frac{1}{\sqrt{\lambda}}^{+
 \infty} \frac{v^{4 \gamma - \alpha - 2}}{(v^4 + 1)^\gamma} dv \bigg)
\bigg).
\end{align*}
Thus, for every $\lambda > 1$ and $(\alpha, \beta, \gamma) \in
(\R_+^*)^3$,
\begin{equation}\label{E57}
J_{\alpha, \beta, \gamma}(\lambda) \leq A \bigg( \lambda^{N - 2 +
 \beta - 2 \gamma} + \lambda^{N - 1 + \alpha + \beta - 4 \gamma} +
\lambda^{N - \frac{3}{2} + \frac{\alpha}{2} + \beta - 2 \gamma}
\bigg).
\end{equation}
Hence, by equations \eqref{E55} and \eqref{E57},
$$\int_{B(0,1)^c} |\partial_1^p R(\xi)|^q d\xi \leq A \int_1^{+
 \infty} \bigg( r^{N - 1 + q(d_1 + d^\perp - p - 4)} + r^{N -
 \frac{3}{2} + q(\frac{d_1}{2} + d^\perp - \frac{p}{2} - 2)} \bigg)
dr.$$
In particular, if $N - 1 + q(d_1 + d^\perp - p - 4) < -1$ and $N -
\frac{3}{2} + q(\frac{d_1}{2} + d^\perp - \frac{p}{2} - 2) < - 1$, the
derivative $\partial_1^p R$ belongs to $L^q(B(0,1)^c)$. However, by
assumption \eqref{E26}, this system of inequalities reduces to
$$q > \frac{2N-1}{p + 4 - d_1 - 2d_\perp},$$
which completes the proof of assertion \eqref{E27}.

Let us now consider the case $q = + \infty$. It follows from Proposition \ref{P3} that
\begin{align*}
|\partial_1^p R(\xi)| \leq & A \sum_{0 \leq k \leq \frac{d_1 + 3p}{4}}
\frac{\max \{|\xi_1| , 1 \}^{d_1 + 3p - 4k} |\xi_\perp|^{2k +
 d_\perp}}{(|\xi|^2 + \xi_1^4)^{p+1}} \\ \leq & A \bigg( \max
\{\xi_1^2 , |\xi| \}^{\frac{d_1}{2} - \frac{p}{2} + d_\perp - 2} +
|\xi|^{\frac{d_1}{2} - \frac{p}{2} + d_\perp - 2} \bigg).
\end{align*}
for any $\xi \in \R^N$. Hence, the function $\partial_1^p R$ belongs to $L^\infty(B(0,1)^c)$ if
assumption \eqref{E28} holds.

{\it Case $j \geq 2$.} We first assume $q < + \infty$. Then, from Proposition \ref{P3} and definition \eqref{E56},
$$\int_{B(0,1)^c} |\partial_j^p R(\xi)|^q d\xi \leq A \sum_{0 \leq k
 \leq \frac{d_\perp + p}{2}} \int_1^{+ \infty} J_{q(4k + d_1),
 q(d_\perp + p - 2k), q(p+1)}(r) dr,$$
so, by formula \eqref{E57},
$$\int_{B(0,1)^c} |\partial_j^p R(\xi)|^q d\xi \leq A \int_1^{+
 \infty} r^{N-\frac{3}{2}} \bigg( r^{\frac{1}{2} + q(d_1 + 2 d^\perp
 - 2p - 4)} + r^{q(\frac{d_1}{2} + d^\perp - p - 2)} \bigg) dr.$$
In particular, if assumption \eqref{E29} holds, the function $\partial_j^p R$ belongs to $L^q(B(0,1)^c)$ for
$$q > \frac{2N - 1}{2p + 4 - d_1 - 2 d_\perp},$$
which completes the proof of assertion \eqref{E30}.

Let us finally consider the case $q = + \infty$. We deduce from Proposition \ref{P3} that for every $\xi \in \R^N \setminus \{ 0 \}$,
\begin{align*}
|\partial_j^p R(\xi)| \leq & A \sum_{0 \leq k \leq \frac{d_\perp + p}{2}} \frac{\max \{|\xi_\perp| , 1 \}^{d_\perp + p - 2k} \max \{|\xi_1| , 1 \}^{4k + d_1}}{(|\xi|^2 + \xi_1^4)^{p+1}} \\ \leq & A \Big( \max \{\xi_1^2 , |\xi| \}^{\frac{d_1}{2} + d_\perp - p - 2} + \max \{ |\xi_\perp|, 1 \}^{d_\perp - p - 2} \\ + & \max \{ |\xi_1|, 1 \}^{d_1 + 2 d_\perp - 2p - 4} \Big).
\end{align*}
Hence, the function $\partial_j^p R$ belongs to $L^\infty(B(0,1)^c)$ if assumption \eqref{E31} holds, which concludes the proof of Proposition \ref{P4}.
\end{proof}

\subsection{Integral representations of some classes of tempered
 distributions}

In this section, we derive the integral expressions of Lemmas \ref{L1} and \ref{L2}, and Proposition \ref{P5}. They link the properties of a tempered distribution $f$ with the properties of the derivatives of its Fourier transform (in particular, when this Fourier transform is a rational fraction $R$ of the form \eqref{E15}-\eqref{E16}). For the sake of completeness, we first show Lemma \ref{L1} (which is probably well-known to the experts and which already appeared in a slightly different form in \cite{Graveja3} and \cite{Graveja5}). Indeed, though it cannot be applied to our kernels, Lemma \ref{L1} is the key ingredient of the proof of Lemma \ref{L2}.

\def\proof{\par{\it Proof of Lemma \ref{L1}}. \ignorespaces}
\begin{proof}
Let us consider $\lambda > 0$, and $\psi \in C^\infty_c(\R^N, [0,1])$ such that
\begin{equation}\label{E58}
\left\{ \begin{array}{ll} {\rm supp}(\psi) \subset B(0,2 \lambda), \\ \forall x \in B(0, \lambda), \psi(x) = 1. \end{array} \right.
\end{equation}
By standard duality, we have for every $g \in \boS(\R^N)$,
\begin{align*}
<x_j^p f, \widehat{g}> = & <x_j^p f, \widehat{\psi g}> + <x_j^p f,
\widehat{(1 - \psi)g}> \\ = & i^p \Big( - <\partial_j^{p-1}
\widehat{f}, \partial_j(\psi g)> + <\partial_j^p \widehat{f}, (1 -
\psi) g> \Big).
\end{align*}
However, by assumptions $(i)$ and $(ii)$, and the smoothness of
$\widehat{f}$, the tempered distributions $\partial_j^p \widehat{f}$
and $\partial_j^{p-1} \widehat{f}$ belong to $L^1(B(0,\lambda)^c)$ and $L^1(B(0,2\lambda))$. Hence, by assumption \eqref{E58},
\begin{equation}\label{E59}
<x_j^p f, \widehat{g}> = i^p \bigg( - \int_{B(0, 2 \lambda)}
\partial_j^{p-1} \widehat{f} (\xi) \partial_j(\psi g) (\xi) d\xi +
\int_{B(0, \lambda)^c} \partial_j^p \widehat{f} (\xi) (1 - \psi(\xi))
g(\xi) d\xi \bigg).
\end{equation}
However, by assumptions $(iii)$ and \eqref{E58}, and integrations by parts,
\begin{align*}
\int_{B(0, 2 \lambda)} \partial_j^{p-1} \widehat{f} (\xi)
\partial_j(\psi g) (\xi) d\xi = & - \int_{B(0, \lambda)} \partial_j^p
\widehat{f} (\xi) (g (\xi) - g(0)) d\xi - \frac{g(0)}{\lambda} \int_{S(0,
 \lambda)} \xi_j \\ & \partial_j^{p-1} \widehat{f}(\xi) d\xi -
\int_{\lambda < |\xi| < 2 \lambda} \partial_j^p \widehat{f} (\xi)
\psi(\xi) g (\xi) d\xi,
\end{align*}
while by assumptions \eqref{E58},
\begin{align*}
\int_{B(0, \lambda)^c} \partial_j^p \widehat{f} (\xi) (1 - \psi(\xi))
g(\xi) d\xi = & \int_{\lambda < |\xi| < 2 \lambda} \partial_j^p
\widehat{f} (\xi) (1 - \psi(\xi)) g (\xi) d\xi \\ + & \int_{B(0, 2
 \lambda)^c} \partial_j^p \widehat{f} (\xi) g(\xi) d\xi.
\end{align*}
Hence, equation \eqref{E59} becomes
\begin{align*}
<x_j^p f, \widehat{g}> = & i^p \bigg( \int_{B(0, \lambda)} \partial_j^p
\widehat{f}(\xi) (g (\xi) - g(0)) d\xi + \frac{g(0)}{\lambda} \int_{S(0,
 \lambda)} \partial_j^{p-1} \widehat{f}(\xi) \xi_j d\xi\\ + &
\int_{B(0, \lambda)^c} \partial_j^p \widehat{f} (\xi) g (\xi) d\xi
\bigg).
\end{align*}
However, since $g$ is in $S(\R^N)$, it satisfies
$$\forall \xi \in \R^N, g(\xi) = \frac{1}{(2\pi)^N} \int_{\R^N}
\widehat{g}(x) e^{ix.\xi} dx,$$
which yields
\begin{align*}
<x_j^p f,\widehat{g}> = & \frac{i^p}{(2\pi)^{N}} \int_{\R^N} \widehat{g}(x)
\bigg( \int_{B(0, \lambda)} \partial_j^p \widehat{f}(\xi) (e^{ix.\xi}
- 1) d\xi + \frac{1}{\lambda} \int_{S(0, \lambda)} \xi_j
\partial_j^{p-1} \widehat{f}(\xi) d\xi \\ + & \int_{B(0, \lambda)^c}
\partial_j^p \widehat{f} (\xi) e^{ix.\xi} d\xi \bigg) dx.
\end{align*}
Therefore, by standard duality, the tempered distribution $x_j^p f$ is
equal to the tempered distribution $\Phi$ defined for every $x \in
\R^N$ by
\begin{equation}\label{E60}
\begin{split}
\Phi(x) = & \frac{i^p}{(2 \pi)^N} \bigg( \int_{B(0, \lambda)}
\partial_j^p \widehat{f}(\xi) (e^{ix.\xi} - 1) d\xi +
\frac{1}{\lambda} \int_{S(0, \lambda)} \partial_j^{p-1}
\widehat{f}(\xi) \xi_j d\xi \\ + & \int_{B(0, \lambda)^c} \partial_j^p
\widehat{f} (\xi) e^{ix.\xi} d\xi \bigg).
\end{split}
\end{equation}
Indeed, by assumptions $(i)$ and $(iii)$, $\Phi$ belongs to
$L^1_{loc}(\R^N)$ and satisfies
$$\forall x \in \R^N, |\Phi(x)| \leq A (1 + |x|).$$
Therefore, $\Phi$ is well-defined on $\R^N$ and is a tempered
distribution. Moreover, it is continuous on $\R^N$ by assumptions
$(i)$ and $(iii)$ once more, and a standard application of the
dominated convergence theorem. Thus, the function $x \mapsto x_j^p
f(x)$ is continuous on $\R^N$ and verifies formula \eqref{E32} by
equation \eqref{E60}.
\end{proof}

Then, Lemma \ref{L2} follows from Lemma \ref{L1}.

\def\proof{\par{\it Proof of Lemma \ref{L2}}. \ignorespaces}
\begin{proof}
Let us consider $\eps > 0$, and let $f_\eps$ denote the tempered distribution defined by
\begin{equation}\label{E61}
\forall \xi \in \R^N, \widehat{f_\eps}(\xi) = \widehat{f}(\xi)
e^{- \eps |\xi|^2}.
\end{equation}
The functions $\widehat{f_\eps}$ decay much faster than $\widehat{f}$ at infinity. In particular, we can apply Lemma \ref{L1} to them, which is not always possible for $\widehat{f}$. Thus, Lemma \ref{L2} results from applying Lemma \ref{L1} to $f_\eps$, and after several integrations by parts, taking the limit $\eps \to 0$. In order to do so, we first establish some properties of $f_\eps$.

\begin{step}\label{S11}
Let $(j,p) \in \{ 1,\ldots,N \} \times \N$. The function
$\widehat{f_\eps}$ belongs to $C^\infty(\R^N \setminus \{ 0 \})$
and satisfies
\begin{equation}\label{E62}
\forall \xi \in \R^N \setminus \{ 0 \}, \partial_j^p
\widehat{f_\eps}(\xi) = \sum_{k=0}^p C_p^k \eps^\frac{k}{2}
\partial_j^{p-k} \widehat{f}(\xi) S_k(\sqrt{\eps} \xi_j) e^{-
 \eps |\xi|^2}.
\end{equation}
Here, the functions $S_k$ are polynomial functions on $\R$ of degree
less or equal to $k$. In particular, $S_0$ is identically equal to the
constant function $1$. 
\end{step}

Indeed, let us consider the function $\Psi$ given by
$$\forall x \in \R, \Psi(x) = e^{-x^2}.$$
By induction, there exist polynomial functions $(S_k)_{k \in \N}$
on $\R$ such that
$$\forall (x,k) \in \R \times \N, \Psi^{(k)}(x) = S_k(x) e^{-x^2}.$$
Moreover, the degree of $S_k$ is less or equal to $k$, and $S_0$ is identically equal to $1$.

However, the function $\widehat{f_\eps}$ is in $C^\infty(\R^N
\setminus \{ 0 \})$ and satisfies by formula \eqref{E61},
$$\forall \xi \in \R^N, \widehat{f_\eps}(\xi) = \widehat{f}(\xi)
\underset{k=1}{\overset{N}{\Pi}} \Psi(\sqrt{\eps} \xi_k).$$
It follows from Leibniz' formula that for every $j \in \{ 1,\ldots,N
\}$ and $p \in \N$,
$$\forall \xi \in \R^N \setminus \{ 0 \}, \partial_j^p
\widehat{f_\eps}(\xi) = \sum_{k=0}^p C_p^k \eps^\frac{k}{2}
\partial_j^{p-k} \widehat{f}(\xi) S_k(\sqrt{\eps} \xi_j) e^{-
 \eps |\xi|^2},$$
which completes the proof of Step \ref{S11}.

We now apply Lemma \ref{L1} to $f_\eps$.

\begin{step}\label{S12}
Let $j \in \{ 1,\ldots,N \}$ and $\lambda > 0$. The function $x \mapsto x_j^p f_\eps(x)$ is continuous on $\R^N$. Moreover, it is equal to
\begin{equation}\label{E63}
\begin{split}
\forall x \in \R^N, x_j^p f_\eps(x) & = \frac{i^p}{(2 \pi)^N}
\bigg( \int_{B(0,\lambda)^c} \partial_j^p \widehat{f_\eps}(\xi)
e^{ix.\xi} d\xi + \frac{1}{\lambda} \int_{S(0,\lambda)} \xi_j
\partial_j^{p-1} \widehat{f_\eps}(\xi) d\xi \\ & +
\int_{B(0,\lambda)} \partial_j^p \widehat{f_\eps}(\xi) (e^{ix.\xi}
- 1) d\xi \bigg).
\end{split}
\end{equation}
\end{step}

Indeed, by formula \eqref{E61}, $f_\eps$ is a tempered
distribution on $\R^N$ whose Fourier transform is in $C^\infty(\R^N
\setminus \{ 0 \})$. Moreover,
$$\forall \xi \in B(0,1)^c, |\partial_j^p \widehat{f_\eps}(\xi)|
\leq A_\eps \sum_{k=0}^p |\xi|^k |\partial_j^{p-k}
\widehat{f}(\xi)| e^{- \eps |\xi|^2},$$
by formula \eqref{E62}, and the functions $\partial_j^k \widehat{f}$ belong to $L^{q_{m-k}}(B(0,1)^c)$ for every $k \in \{ 0,\ldots,p \}$ by assumption $(iv)$, so, the function $\partial_j^p \widehat{f_\eps}$ belongs to $L^1(B(0,1)^c)$. On the other hand, by formula
\eqref{E62},
$$\forall \xi \in B(0,1), |\partial_j^{p-1} \widehat{f_\eps}(\xi)|
\leq A_\eps \sum_{k=0}^{p-1} |\partial_j^{p-1-k}
\widehat{f}(\xi)|,$$
so, by assumption $(ii)$,
$$\forall \xi \in B(0,1), |\partial_j^{p-1} \widehat{f_\eps}(\xi)|
\leq A_\eps \sum_{k=0}^{p-1} |\xi|^{1-N+k} \leq
\frac{A_\eps}{|\xi|^{N-1}}.$$
Hence, $\partial_j^{p-1} \widehat{f_\eps}$ is in $L^1(B(0,1))$. Likewise, by formula \eqref{E62} and assumption $(ii)$,
$$\forall \xi \in B(0,1), |\xi| |\partial_j^p
\widehat{f_\eps}(\xi)| \leq \frac{A_\eps}{|\xi|^{N-1}},$$
so $|.| \partial_j^p \widehat{f_\eps}$ also belongs to $L^1(B(0,1))$. Finally, $f_\eps$ satisfies all the assumptions of Lemma \ref{L1}. Thus, the function $x \mapsto x_j^p f_\eps(x)$ is continuous on $\R^N$ and satisfies equation \eqref{E63}.

We now integrate equation \eqref{E63} by parts in order to compute
formula \eqref{E33} for $f_\eps$.

\begin{step}\label{S13}
Let $\Omega_j = \{ x \in \R^N, x_j \neq 0 \}$. Then, for every $x \in
\Omega_j$,
\begin{equation}\label{E64}
\begin{split}
x_j^p f_\eps(x) = & \frac{i^p}{(2 \pi)^N} \bigg( (-ix_j)^{p-m}
\int_{B(0, \lambda)^c} \partial_j^m \widehat{f_\eps}(\xi) e^{ix
 .\xi} d\xi + \sum_{k=p}^{m-1} \frac{(-ix_j)^{p-k-1}}{\lambda}
\int_{S(0, \lambda)} \xi_j \\ & \partial_j^k \widehat{f_\eps}(\xi)
e^{ix.\xi} d\xi + \frac{1}{\lambda} \int_{S(0, \lambda)} \xi_j
\partial_j^{p-1} \widehat{f_\eps}(\xi) d\xi + \int_{B(0, \lambda)}
\partial_j^p \widehat{f_\eps}(\xi) (e^{ix .\xi} - 1) d\xi
\bigg).
\end{split}
\end{equation}
\end{step}

Indeed, let us consider the integrals defined for every $x \in \Omega_j$ and $k \in \{ p,\ldots,m \}$ by
\begin{equation}\label{E65}
I_k(x) = (-ix_j)^{-k} \int_{B(0,\lambda)^c} \partial_j^k
\widehat{f_\eps}(\xi) e^{ix.\xi} d\xi.
\end{equation}
By formula \eqref{E62}, we compute
$$\forall \xi \in B(0,1)^c, |\partial_j^k \widehat{f_\eps}(\xi)|
\leq A_\eps \sum_{l=0}^k |\xi|^l |\partial_j^{k-l}
\widehat{f}(\xi)| e^{- \eps |\xi|^2}.$$
Since the functions $\partial_j^l \widehat{f}$ are in $C^\infty(\R^N
\setminus \{ 0 \})$ and $L^{q_{m-l}}(B(0,1)^c)$ for any $l \in \{
0,\ldots,m \}$ (with $q_m=1$) by assumptions $(iii)$ and $(iv)$, the
derivative $\partial_j^k \widehat{f_\eps}$ belongs to
$L^1(B(0,1)^c)$. Thus, the integral $I_k(x)$ is well-defined for every $x \in \Omega_j$.

Moreover, an integration by parts of $I_k(x)$ leads to
$$\forall k \in \{ p,\ldots,m-1 \}, I_k(x) = I_{k+1}(x) +
\frac{(-ix_j)^{-k-1}}{\lambda} \int_{S(0,\lambda)} \partial_j^k
\widehat{f_\eps}(\xi) \xi_j e^{ix.\xi} d\xi,$$
so, by summing this identity for any $k$ between $p$ and $m-1$,
$$I_p(x) = I_m(x) + \frac{1}{\lambda} \sum_{k=p}^{m-1} (-ix_j)^{-k-1}
\int_{S(0,\lambda)} \partial_j^k \widehat{f_\eps}(\xi) \xi_j
e^{ix.\xi} d\xi,$$
which gives by definition \eqref{E65},
\begin{align*}
\int_{B(0,\lambda)^c} \partial_j^p \widehat{f_\eps}(\xi)
e^{ix.\xi} d\xi = & (-ix_j)^{p-m} \int_{B(0,\lambda)^c} \partial_j^m
\widehat{f_\eps}(\xi) e^{ix.\xi} d\xi \\ + & \sum_{k=p}^{m-1}
\frac{(-ix_j)^{p -k-1}}{\lambda} \int_{S(0,\lambda)} \partial_j^k
\widehat{f_\eps}(\xi) \xi_j e^{ix.\xi} d\xi.
\end{align*}
Then, equation \eqref{E64} follows by replacing the first term in the right-hand side of equation \eqref{E63} by the expression just above.

In order to derive formula \eqref{E33}, it only remains to study the convergences when $\eps$ tends to $0$ of both sides of equation \eqref{E64}. We begin by the convergence of the left-hand side.

\begin{step}\label{S14}
Let $j \in \{ 1,\ldots,N \}$. Then,
\begin{equation}\label{E66}
x_j^p f_\eps \underset{\eps \to 0}{\to} x_j^p f \ {\rm in} \
\boD'(\R^N).
\end{equation}
\end{step}

Indeed, by standard duality, we have
$$\forall \phi \in C_c^\infty(\R^N), <x_j^p f_\eps, \phi> = i^p
<\widehat{f_\eps}, \partial_j^p \widehat{\phi}>.$$
Moreover, by assumption $(i)$ and formula \eqref{E61},
$\widehat{f_\eps}$ is in $L^1(\R^N)$, so
$$<x_j^p f_\eps, \phi> = i^p \int_{\R^N} \widehat{f_\eps}(\xi)
\partial_j^p \widehat{\phi}(\xi) d\xi = i^p \int_{\R^N}
\widehat{f}(\xi) \partial_j^p \widehat{\phi}(\xi) e^{- \eps
 |\xi|^2} d\xi.$$
On the other hand, it follows from assumption $(i)$ that
$$\forall \xi \in \R^N, |\widehat{f}(\xi) \partial_j^p
\widehat{\phi}(\xi) e^{- \eps |\xi|^2}| \leq A |\xi|^{-r}
|\partial_j^p \phi(\xi)|.$$
Therefore, by the dominated convergence theorem,
$$\int_{\R^N} \widehat{f}(\xi) \partial_j^p \widehat{\phi}(\xi) e^{-
 \eps |\xi|^2} d\xi \underset{\eps \to 0}{\to} \int_{\R^N}
\widehat{f}(\xi) \partial_j^p \widehat{\phi}(\xi) d\xi.$$
However, by assumption $(i)$, the function $\widehat{f}$ belongs to
$L_{loc}^1(\R^N)$, so, by standard duality,
$$\int_{\R^N} \widehat{f}(\xi) \partial_j^p \widehat{\phi}(\xi) d\xi =
(-i)^p <x_j^p f, \phi>.$$
Finally,
$$<x_j^p f_\eps, \phi> \underset{\eps \to 0}{\to} <x_j^p f,
\phi>,$$
and the distribution $x_j^p f_\eps$ converges towards $x_j^p f$ in
$\boD'(\R^N)$.

Then, we consider the convergence of the right-hand side of equation
\eqref{E64}.

\begin{step}\label{S15}
Let $g: \Omega_j \mapsto \R$ be the function defined for every $x \in \Omega_j$ by
\begin{equation}\label{E67}
\begin{split}
g(x) = & \frac{i^p}{(2 \pi)^N} \bigg( (-ix_j)^{p-m} \int_{B(0,
 \lambda)^c} \partial_j^m \widehat{f}(\xi) e^{ix .\xi} d\xi +
\sum_{k=p}^{m-1} \frac{(-ix_j)^{p-k-1}}{\lambda} \int_{S(0, \lambda)}
\xi_j \\ &\partial_j^k \widehat{f}(\xi) e^{ix.\xi} d\xi +
\frac{1}{\lambda} \int_{S(0, \lambda)} \xi_j \partial_j^{p-1}
\widehat{f}(\xi) d\xi + \int_{B(0, \lambda)} \partial_j^p
\widehat{f}(\xi) (e^{ix .\xi} - 1) d\xi \bigg).
\end{split}
\end{equation}
Then, $g$ is continuous on $\Omega_j$ and satisfies
\begin{equation}\label{E68}
x_j^p f_\eps \underset{\eps \to 0}{\to} g \ {\rm in} \ L_{loc}^\infty(\Omega_j).
\end{equation}
\end{step}

Indeed, by assumptions $(ii)$ and $(iii)$, the functions $\partial_j^m
\widehat{f}$ and $|.| \partial_j^p \widehat{f}$ respectively belong
to $L^1(B(0,\lambda)^c)$ and $L^1(B(0,\lambda))$. Moreover, the
function $\widehat{f}$ is in $C^\infty(\R^N \setminus \{ 0
\})$. Therefore, by equation \eqref{E67} and a standard application of
the dominated convergence theorem, $g$ is well-defined and continuous
on $\R^N$.

Moreover, equations \eqref{E62}, \eqref{E64} and \eqref{E67} lead to
\begin{equation}\label{E69}
\begin{split}
& \frac{(2 \pi)^N}{i^p} \Big( x_j^p f_\eps(x) - g(x) \Big) = (-ix_j)^{p-m} \int_{B(0, \lambda)^c} \sum_{k=1}^m C_m^k \eps^\frac{k}{2} \partial_j^{m-k} \widehat{f}(\xi) S_k(\sqrt{\eps} \xi_j) e^{- \eps |\xi|^2} \\ & e^{ix.\xi} d\xi + \sum_{k=p}^{m-1} \frac{(-ix_j)^{p-k-1}}{\lambda} \int_{S(0, \lambda)} \sum_{l=1}^k C_k^l \eps^\frac{l}{2} \partial_j^{k-l} \widehat{f}(\xi) S_l(\sqrt{\eps} \xi_j) e^{- \eps |\xi|^2} e^{ix.\xi} \xi_j d\xi + \frac{1}{\lambda} \\ & \int_{S(0, \lambda)} \sum_{k=1}^{p-1} C_{p-1}^k \eps^\frac{k}{2} \xi_j \partial_j^{p-k-1} \widehat{f}(\xi) S_k(\sqrt{\eps} \xi_j) e^{-  \eps |\xi|^2} d\xi + \int_{B(0, \lambda)} \sum_{k=1}^p C_p^k \eps^\frac{k}{2} \partial_j^{p-k} \widehat{f}(\xi) \\ & S_k(\sqrt{\eps} \xi_j) e^{- \eps |\xi|^2} (e^{ix .\xi} - 1) d\xi,
\end{split}
\end{equation}
for every $x \in \Omega_j$. Then, assumption $(iv)$ and Step \ref{S11} give for the first term in the right-hand side
\begin{align*}
& \bigg|\sum_{k=1}^m \int_{B(0, \lambda)^c} C_m^k \eps^\frac{k}{2} \partial_j^{m-k} \widehat{f}(\xi) S_k(\sqrt{\eps} \xi_j) e^{- \eps |\xi|^2 + ix.\xi} d\xi \bigg| \\ \leq & A \sum_{k=1}^m \eps^\frac{k}{2} \|\partial_j^{m-k} \widehat{f} \|_{L^{q_k}(B(0,\lambda)^c)} \bigg( \int_{B(0,\lambda)^c} |S_k(\sqrt{\eps} \xi_j)|^{q_k'} e^{- \eps q_k' |\xi|^2} d\xi \bigg)^\frac{1}{q_k'} \\ \leq & A \sum_{k=1}^m \eps^{\frac{k}{2} - \frac{N}{2 q_k'}} \bigg( \int_{\R^N} |S_k(\xi_j)|^{q_k'} e^{-q_k' |\xi|^2} d\xi \bigg)^\frac{1}{q_k'}
\end{align*}
where $q_k' = \frac{q_k}{q_k - 1}$. Moreover, we know that $$1 < q_k <
\frac{N}{N-k}$$ if $1 \leq k \leq N-1$, and $1 < q_k \leq + \infty$,
otherwise. Therefore, there is some $\delta > 0$ such that for every
$k \in \{ 1, \ldots, m \}$,
$$\frac{k}{2} - \frac{N}{2 q_k'} \geq \delta.$$
Hence, we obtain
\begin{equation}\label{E70}
\bigg| \sum_{k=1}^m \int_{B(0, \lambda)^c} C_m^k \eps^\frac{k}{2}
\partial_j^{m-k} \widehat{f}(\xi) S_k(\sqrt{\eps} \xi_j) e^{-
 \eps |\xi|^2 + ix.\xi} d\xi \bigg| \leq A \eps^\delta.
\end{equation}
Then, by Step \ref{S11} and the smoothness of $\widehat{f}$ on $S(0, \lambda)$, the second term in the right-hand side of equation \eqref{E69} satisfies
\begin{equation}\label{E71}
\begin{split}
& \bigg| \sum_{k=p}^{m-1} \frac{(-ix_j)^{p-k-1}}{\lambda} \int_{S(0, \lambda)} \xi_j \sum_{l=1}^k C_k^l \eps^\frac{l}{2} \partial_j^{k-l} \widehat{f}(\xi) S_l(\sqrt{\eps} \xi_j) e^{-
\eps |\xi|^2 + ix.\xi} d\xi \bigg| \\ \leq & A \eps^\frac{1}{2} \sum_{k=p}^{m-1} |x_j|^{p-k-1}.
\end{split}
\end{equation}
Likewise, by Step \ref{S11} and the smoothness of $\widehat{f}$ on
$S(0, \lambda)$, the third term in the right-hand side of equation \eqref{E69} verifies
\begin{equation}\label{E72}
\bigg| \frac{1}{\lambda} \int_{S(0, \lambda)} \xi_j \sum_{k=1}^{p-1}
C_{p-1}^k \eps^\frac{k}{2} \partial_j^{p-k-1} \widehat{f}(\xi)
S_k(\sqrt{\eps} \xi_j) e^{- \eps |\xi|^2} d\xi \bigg| \leq A
\eps^\frac{1}{2}.
\end{equation}
Finally, by Step \ref{S11} and assumption $(ii)$, the last term in the right-hand side of equation \eqref{E69} verifies 
\begin{equation}\label{E73}
\begin{split}
& \bigg| \int_{B(0, \lambda)} \sum_{k=1}^p C_p^k \eps^\frac{k}{2}
\partial_j^{p-k} \widehat{f}(\xi) S_k(\sqrt{\eps} \xi_j) e^{-
 \eps |\xi|^2} (e^{ix .\xi} - 1) d\xi \bigg| \\ \leq & A |x|
\eps^\frac{1}{2} \int_{B(0, \lambda)} \sum_{k=1}^p |\xi|^{k-N+1}
d\xi \leq  A \eps^\frac{1}{2} |x|.
\end{split}
\end{equation}
Thus, equations \eqref{E69}, \eqref{E70}, \eqref{E71}, \eqref{E72}
and \eqref{E73} give
$$\forall x \in \Omega_j, |x_j^p f_\eps(x) - g(x)| \leq A
\bigg( \eps^\delta + \eps^\frac{1}{2} \Big( 1 + |x| +
\sum_{k=p}^{m-1} |x_j|^{p-k-1} \Big) \bigg).$$
Therefore, the functions $x \mapsto x_j^p f_\eps$ converge towards
$g$ in $L^\infty_{loc}(\Omega_j)$ when $\eps$ tends to $0$.

Finally, by assertions \eqref{E66} and \eqref{E68}, and the uniqueness
of the limit of the functions $x \mapsto x_j^p f_\eps$, the
function $x \mapsto x_j^p f$ is identically equal to $g$ on
$\Omega_j$. In particular, by Step \ref{S15}, it is continuous on
$\Omega_j$ and satisfies equation \eqref{E33}, which concludes the proof of Lemma \ref{L2}.
\end{proof}

Then, Lemma \ref{L2} yields the integral expression \eqref{E34} of
Proposition \ref{P5} for the tempered distributions $f$ whose Fourier transforms are rational fractions $R$ of the form \eqref{E15}-\eqref{E16}.

\def\proof{\par{\it Proof of Proposition \ref{P5}}. \ignorespaces}
\begin{proof}
Let $p_j = N - 2 + d$, $m_1 = 2N - 4 + d_1 + 2d_\perp$ and
$m_j = N - 2 + d$ if $j \geq 2$. In order to prove formula
\eqref{E34}, we apply Lemma \ref{L2} to $f$ with $p = p_j$ and $m =
m_j$.

Indeed, by Proposition \ref{P4}, the tempered distribution $f$ satisfies all the assumptions of Lemma \ref{L2}. By formulae \eqref{E15}-\eqref{E16}, its Fourier transform $R$ belongs to $C^\infty(\R^N \setminus \{ 0 \})$ and satisfies
$$\forall \xi \in \R^N \setminus \{ 0 \}, |R(\xi)| \leq A \bigg( |\xi|^d + |\xi|^{d - 2} \bigg).$$
Thus, since $d \neq 0$ if $N = 2$, assumption $(i)$ is satisfied by $f$.
 
Moreover, by Proposition \ref{P4}, there is some positive real number $A$ such that
$$\forall (k,\xi) \in \{ 0,\ldots,N-1 \} \times B(0,1), |\xi|^{k + 2 - d} |\partial_j^k R(\xi)| \leq A.$$
Therefore, $f$ verifies assumption $(ii)$ with $p_j = N - 2 + d$.

On the other hand, by Proposition \ref{P4}, the function $\partial_j^{m_j} R$ belongs to $L^1(B(0,1))^c$. Moreover, since $d_1 + 2 d_\perp \leq 4$, the function $\partial_j^k R$ belongs to $L^q(B(0,1)^c)$ for every $k \in \{ 0,\ldots,m_j-1 \}$, when $q > q_{m_j-k}^j = \frac{2N-1}{(1 + \delta_{j,1}) k + 4 - d_1 - 2d_\perp}$. In particular, we compute
$$\forall k \in \{ 1,\ldots,N-1 \}, q_k^1 = \frac{2N-1}{2N - k} < \frac{N}{N - k},$$
while for every $j \in \{ 2,\ldots,N \}$,
$$\forall k \in \{ 1,\ldots,N-1 \}, q_k^j = \frac{2N-1}{2N - 2k + d_1} < \frac{N}{N - k}.$$
Therefore, the distribution $f$ also verifies assumptions $(iii)$ and $(iv)$ of Lemma \ref{L2}, and we can apply it with $p = p_j$ and $m = m_j$. Thus, the function $x \mapsto x_j^{p_j} f(x)$ is continuous on $\Omega_j$, and formula \eqref{E34} holds by equation \eqref{E33}.
\end{proof}

\subsection{Algebraic decay and explosion near the origin of the kernels $H_0$, $K_0$ and $K_k$}

This section deals with the proofs of Propositions \ref{Paumond1} and
\ref{P2} (which immediately yield Theorems \ref{T3} and \ref{T4}, and Corollary \ref{C1}). Their proofs follow from expression \eqref{E34} of Proposition \ref{P5}. We estimate it for any value of $\lambda$ before choosing it either to compute the algebraic decay at infinity or the explosion near the origin of the considered distribution. Let us begin by the algebraic decay properties of Proposition \ref{Paumond1} and Theorem \ref{T3}.

\def\proof{\par{\it Proof of Proposition \ref{Paumond1}}. \ignorespaces}
\begin{proof}
By Proposition \ref{P5}, the functions $x \mapsto x_j^{p_j} f(x)$ are continuous on $\Omega_j$ for every $j \in \{ 1,\ldots,N
\}$. Therefore, the restriction of $f$ to $\R^N \setminus \{ 0 \}$ is continuous. Moreover, by choosing $x \in B(0,1)^c$, $j \in \{ 1,\ldots,N \}$ such that $|x_j| \geq \frac{|x|}{\sqrt{N}}$, and $\lambda = \frac{1}{|x|}$, we infer from Proposition \ref{P4} and formula \eqref{E34} that
\begin{align*}
|x|^{N-2+d} |f(x)| \leq & A \frac{|x|^{N-2+d}}{|x_j|^{N-2+d}} \bigg(
|x_j|^{N-2+d-m_j} \bigg( \int_{B(0,1)^c} |\partial_j^{m_j}
\widehat{f}(\xi)| d\xi \\ + & \int_{\frac{1}{|x|} \leq |\xi| \leq 1}
\frac{d\xi}{|\xi|^{m_j+2-d}} \bigg) + \sum_{k=N-2+d}^{m_j-1}
|x_j|^{N-3+d-k} \int_{S(0,\frac{1}{|x|})} \frac{d\xi}{|\xi|^{k+2-d}}
\\ + & \int_{S(0,\frac{1}{|x|})} \frac{d\xi}{|\xi|^{N-1}} +
|x| \int_{B(0,\frac{1}{|x|})} \frac{d\xi}{|\xi|^{N-1}} \bigg).
\end{align*}
Hence, since $|x_j| \geq \frac{|x|}{\sqrt{N}}$ and $m_j \geq N - 2 +
d$,
$$|x|^{N-2+d} |f(x)| \leq A_N (|x|^{N-2+d-m_j} + 1 ) \leq A_N.$$
Thus, the distribution $f$ belongs to $M^\infty_{N-2+d}(B(0,1)^c)$.
\end{proof}

Then, we deduce Theorem \ref{T3}.

\def\proof{\par{\it Proof of Theorem \ref{T3}}. \ignorespaces}
\begin{proof}
Theorem \ref{T3} is a direct application of Proposition \ref{Paumond1} with $(d_1,d_\perp) = (1,0)$ for $H_0$, $(d_1,d_\perp) = (2,0)$ for $K_0$ and $(d_1,d_\perp) = (2+\delta_{k,1}, 1-\delta_{k,1})$ for $K_k$.
\end{proof}

We now turn to the study of the singularities of $f$ near the origin.

\def\proof{\par{\it Proof of Proposition \ref{P2}}. \ignorespaces}
\begin{proof}
We first derive the local estimates \eqref{E18} for $f$. Then, we describe its singularities in terms of $L^q$-spaces.

The proof relies on expression \eqref{E34}. Indeed, the distribution $f$ satisfies the assumptions of Proposition \ref{P5}. Therefore, by the same argument as in the proof of Proposition \ref{Paumond1}, its restriction to $B(0,1) \setminus \{ 0 \}$ is continuous. Moreover, its restriction satisfies equation \eqref{E34}. Thus, it remains to bound the second term of equation \eqref{E34} independently of $x$ to obtain the local estimates \eqref{E18}. However, there are at least two difficulties. The first one results from the form of equation \eqref{E34}. Indeed, this equation is suitable to give the optimal algebraic decay of $f$, but not its optimal algebraic explosion near the origin. More precisely, an ingenuous way to estimate the last term of equation \eqref{E34} consists in writing by Proposition \ref{P4},
\begin{align*}
\bigg| \int_{B(0,\lambda)} \partial_j^{p_j} \widehat{f}(\xi)
 (e^{ix.\xi} - 1) d\xi \bigg| \leq & A |x| \int_{B(0,1)}
|\partial_j^{p_j} \widehat{f}(\xi)| |\xi| d\xi + A \int_{1 \leq |\xi|
 \leq \lambda} |\partial_j^{p_j} \widehat{f}(\xi)| d\xi \\ \leq & A
|x| + A_\lambda.
\end{align*}
By this argument, we cannot expect to prove (for any choice of $\lambda$) that this term is bounded by $A |x_j|^d$ for some $A \geq 0$ and $d \in \N$, which is actually the goal of Proposition \ref{P2}. Thus, in order to precisely describe the singularities of $f$ near the origin, we integrate by parts the last term of equation \eqref{E34}.

Indeed, consider $j \in \{ 1,\ldots,N \}$, $\lambda > 1$ and $x \in
B(0,1) \cap \Omega_j$. By Proposition \ref{P4}, the functions
$\partial_j^k \widehat{f}$ belong to $L^1(B(0,\lambda))$ for every $0
\leq k < p_j$. Therefore, by integrating by parts the integral above, we deduce for any integer $1 \leq q_j \leq p_j$,
\begin{align*}
\int_{B(0,\lambda)} \partial_j^{p_j} \widehat{f}(\xi) (e^{ix.\xi} - 1) d\xi = & \sum_{k=1}^{q_j} \frac{(-ix_j)^{k-1}}{\lambda} \int_{S(0,\lambda)} \xi_j e^{ix. \xi} \partial_j^{p_j - k} \widehat{f}(\xi) d\xi - \int_{S(0,\lambda)} \frac{\xi_j}{\lambda} \\ & \partial_j^{p_j - 1} \widehat{f}(\xi) d\xi + (-i x_j)^{q_j} \int_{B(0,\lambda)} e^{ix. \xi} \partial_j^{p_j-q_j} \widehat{f}(\xi) d\xi.
\end{align*}
Hence, equation \eqref{E34} becomes
\begin{equation}\label{E74}
\begin{split}
x_j^{p_j} f(x) = & \frac{i^{p_j}}{(2 \pi)^N} \bigg( (-ix_j)^{p_j -
 m_j} \int_{B(0,\lambda)^c} \partial_j^{m_j} \widehat{f}(\xi)
e^{ix.\xi} d\xi + \sum_{k=p_j}^{m_j-1} \frac{(-ix_j)^{p_j - k -
 1}}{\lambda} \\ & \int_{S(0,\lambda)} \xi_j \partial_j^k
\widehat{f}(\xi) e^{ix.\xi} d\xi + \sum_{k=1}^{q_j}
\frac{(-ix_j)^{k-1}}{\lambda} \int_{S(0,\lambda)} \xi_j
\partial_j^{p_j - k} \widehat{f}(\xi) e^{ix. \xi} d\xi \\ + & (-i
x_j)^{q_j} \int_{B(0,\lambda)} e^{ix. \xi} \partial_j^{p_j-q_j}
\widehat{f}(\xi) d\xi \bigg).
\end{split}
\end{equation}
In particular, the argument used above to bound the integral on
$B(0,\lambda)$ now gives by Proposition \ref{P4},
\begin{align*}
& \bigg| (-i x_j)^{q_j} \int_{B(0,\lambda)} e^{ix. \xi}
 \partial_j^{p_j-q_j} \widehat{f}(\xi) d\xi \bigg| \\ \leq & A
|x_j|^{q_j} \bigg( \int_{B(0,1)} |\partial_j^{p_j - q_j}
\widehat{f}(\xi)| d\xi + \int_{1 \leq |\xi| \leq \lambda}
|\partial_j^{p_j - q_j} \widehat{f}(\xi)| d\xi \bigg) \leq A_\lambda
|x_j|^{q_j}.
\end{align*}
Therefore, by choosing $q_j$ and $\lambda$ appropriately, we can now
sharply describe the singularities of $f$ near the origin.

However, the anisotropy of the considered distributions is a second difficulty. Their explosion near the origin is not the same in every direction, so we split our analysis into two steps: first, the case of the direction $x_1$, second, the case of the directions $x_j$ with $j \neq 1$. In particular, we choose different values of $q_j$ and $\lambda$ according to the value of $j$.

\setcounter{step}{0}
\begin{step}\label{S21}
Estimates of $f$ near the origin in the direction $x_1$.
\end{step}

Let us consider $x \in B(0,1) \cap \Omega_1$ and $\lambda \geq 1$. By Proposition \ref{P5}, the restriction of $f$ to $B(0,1) \cap \Omega_1$ is continuous and satisfies by equation \eqref{E74},
\begin{equation}\label{E75}
\begin{split}
|x_1|^{m_1} |f(x)| \leq & A \bigg( \int_{B(0,\lambda)^c}
|\partial_1^{m_1} \widehat{f}(\xi)| d\xi + \sum_{k=p_1-q_1}^{m_1-1}
\frac{|x_1|^{m_1 - 1 - k}}{\lambda} \int_{S(0,\lambda)} |\partial_1^k
\widehat{f}(\xi)| |\xi_1| d\xi \\ + & |x_1|^{m_1-p_1+q_1}
\int_{B(0,\lambda)} |\partial_1^{p_1-q_1} \widehat{f}(\xi)| d\xi \bigg),
\end{split}
\end{equation}
where we denote as above $p_1 = N-2+d$ and $m_1 = 2N-4+d_1+2d_\perp$,
and choose $q_1 = \min \{ p_1,2 \}$
\footnote{Indeed, since $d \neq 0$ if $N=2$, $q_1$ satisfies in any
 case $1 \leq q_1 \leq p_1$.}.
Then, the first term in the right-hand side of equation \eqref{E75} verifies by Proposition \ref{P3} and equation \eqref{E56},
$$\int_{B(0,\lambda)^c} |\partial_1^{m_1} \widehat{f}(\xi)| d\xi \leq
A \sum_{0 \leq k \leq \frac{d_1 + 3m_1}{4}} \int_\lambda^{+ \infty}
J_{d_1 + 3m_1 - 4k, 2k + d_\perp, m_1 + 1}(r) dr,$$
so, by equation \eqref{E57},
\begin{equation}\label{E76}
\int_{B(0,\lambda)^c} |\partial_1^{m_1} \widehat{f}(\xi)| d\xi \leq A
\int_\lambda^{+ \infty} r^{-\frac{3}{2}} dr \leq
\frac{A}{\sqrt{\lambda}}.
\end{equation}
Likewise, Proposition \ref{P3}, and equations \eqref{E56} and \eqref{E57} lead to
\begin{equation}\label{E77}
\int_{S(0,\lambda)} |\xi_1| |\partial_1^k \widehat{f}(\xi)| d\xi
\leq A \sum_{0 \leq l \leq \frac{d_1 + 3k}{4}} J_{1 + d_1 + 3k - 4l,
 2l + d_\perp, k+1}(\lambda) \leq A \lambda^{N - 3 + \frac{d_1}{2} +
 d_\perp - \frac{k}{2}}.
\end{equation}
for every $p_1 - q_1 \leq k \leq m_1 - 1$. Finally, by Propositions \ref{P3} and \ref{P4}, and equation \eqref{E56}, the last term in the right-hand side of equation \eqref{E75} satisfies
\begin{align*}
\int_{B(0,\lambda)} |\partial_1^{p_1 - q_1} \widehat{f}(\xi)| d\xi
\leq & \int_{B(0,1)} |\partial_1^{p_1 - q_1} \widehat{f}(\xi)| d\xi + \int_{1 \leq |\xi| \leq \lambda} |\partial_1^{p_1 - q_1}
\widehat{f}(\xi)| d\xi \\ \leq A & + A \sum_{0 \leq k \leq
 \frac{d_1 + 3p_1 - 3q_1}{4}} \int_1^\lambda J_{d_1 + 3p_1 - 3q_1 -
 4k, 2k + d_\perp, p_1-q_1+1}(r) dr,
\end{align*}
which gives
\begin{equation}\label{E78}
\int_{B(0,\lambda)} |\partial_1^{p_1 - q_1} \widehat{f}(\xi)| d\xi
\leq A \Big( 1 + (1 - \delta) \lambda^{\frac{N - 3 + d_\perp +
 q_1}{2}} + \delta \ln(\lambda) \Big),
\end{equation}
by equation \eqref{E57} (with $\delta = \delta_{N,2} \delta_{d_1,1} \delta_{d_\perp,0}$). Then, it follows from equations \eqref{E75}, \eqref{E76}, \eqref{E77} and \eqref{E78} that
\begin{align*}
|x_1|^{m_1} |f(x)| \leq & A \bigg( \lambda^{-\frac{1}{2}} +
\sum_{k=p_1-q_1}^{m_1-1} \Big( |x_1|^{2N - 5 + d_1 + 2d_\perp - k}
\lambda^{N - 4 + d_\perp + \frac{d_1 - k}{2}} \Big) \\ + &
|x_1|^{N-2+d_\perp+q_1} \Big( 1 + (1 - \delta) \lambda^{\frac{N -
 3 + d_\perp + q_1}{2}} + \delta \ln(\lambda) \Big) \bigg).
\end{align*}
However, since $|x| < 1$, we can set $\lambda = \frac{1}{x_1^2} > 1$
to obtain
$$|x_1|^{m_1} |f(x)| \leq A \bigg( |x_1| + |x_1|^3 + |x_1|^{N - 2 +
 d_\perp + q_1} \Big( 1 + \delta |\ln(|x_1|)| \Big) \bigg),$$
and to deduce the estimate of $f$ near the origin in the direction $x_1$:
\begin{equation}\label{E79}
|x_1|^{2N-5+d_1+2d_\perp} |f(x)| \leq A \Big( 1 + \delta_{N,2}
\delta_{d_1,1} \delta_{d_\perp,0} |\ln(|x_1|)| \Big).
\end{equation}

\begin{step}\label{S22}
Estimates of $f$ near the origin in the directions $x_j$, $j \neq 1$.
\end{step}

Now, we consider the singularity of $f$ near the origin in every
direction $x_j$ with $j \in \{ 2,\ldots,N \}$. In order to do so, let $x \in B(0,1) \cap \Omega_j$ and $\lambda \geq 1$. Then, equation \eqref{E74} yields
\begin{equation}\label{E80}
\begin{split}
|x_j|^{p_j} |f(x)| \leq & A \bigg( \int_{B(0,\lambda)^c}
|\partial_j^{p_j} \widehat{f}(\xi)| d\xi + \sum_{k=1}^{q_j}
\frac{|x_j|^{k-1}}{\lambda} \int_{S(0,\lambda)} |\xi_j|
|\partial_j^{p_j-k} \widehat{f}(\xi)| d\xi \\ + & |x_j|^{q_j}
\int_{B(0,\lambda)} |\partial_j^{p_j - q_j} \widehat{f}(\xi)| d\xi
\bigg),
\end{split}
\end{equation}
where we let $p_j = N - 2 + d$ and choose $q_j = \min \{ k \in \N,
k > \frac{d_1 + 1}{2} \} - \delta$ with $\delta = \delta_{N,2}
\delta_{d_1,1} \delta_{d_\perp,0}$ as above
\footnote{Indeed, since $d \neq 0$ if $N=2$, $q_j$ satisfies in any
 case $1 \leq q_j \leq p_j$.}.
On one hand, by Proposition \ref{P3} and equation \eqref{E56},
$$\int_{B(0,\lambda)^c} |\partial_j^{p_j} \widehat{f}(\xi)| d\xi \leq A \sum_{0 \leq k \leq \frac{d_\perp + p_j}{2}} \int_\lambda^{+ \infty} J_{4k + d_1, d_\perp + p_j - 2k, p_j + 1}(r) dr,$$
so, by equation \eqref{E57},
\begin{equation}\label{E81}
\int_{B(0,\lambda)^c} |\partial_j^{p_j} \widehat{f}(\xi)| d\xi \leq A \int_\lambda^{+ \infty} \frac{dr}{r^{\frac{3 + d_1}{2}}} \leq
\frac{A}{\lambda^{\frac{1+d_1}{2}}}.
\end{equation}
On the other hand, Proposition \ref{P3}, and equations \eqref{E56} and \eqref{E57} yield
\begin{equation}\label{E82}
\begin{split}
\int_{S(0,\lambda)} |\xi_j| |\partial_j^{p_j-k} \widehat{f}(\xi)| d\xi \leq & A \sum_{0 \leq l \leq \frac{p_j - k + d_\perp}{2}} J_{d_1 + 4l, d_\perp + 1 + p_j - k - 2l, p_j - k + 1}(\lambda) \\ \leq & A \lambda^{-\frac{d_1+1}{2} + k}.
\end{split}
\end{equation}
for every $1 \leq k \leq q_j$. Finally, we deduce from Proposition \ref{P3} and equation \eqref{E56},
\begin{align*}
\int_{B(0,\lambda)} |\partial_j^{p_j-q_j} \widehat{f}(\xi)| d\xi \leq & \int_{B(0,1)} |\partial_j^{p_j-q_j} \widehat{f}(\xi)| d\xi + \int_{1 \leq |\xi| \leq \lambda} |\partial_j^{p_j-q_j} \widehat{f}(\xi)| d\xi \\ \leq A & \bigg( 1 + \sum_{0 \leq k \leq \frac{p_j-q_j+d_\perp}{2}} \int_1^\lambda J_{4k + d_1, d_\perp + p_j - q_j - 2k,p_j-q_j+1}(r) dr \bigg),
\end{align*}
which yields by equation \eqref{E57},
\begin{equation}\label{E83}
\int_{B(0,\lambda)} |\partial_j^{p_j-q_j} \widehat{f}(\xi)| d\xi \leq
A \Big( 1 + (1 - \delta) \lambda^{\frac{-1-d_1}{2} + q_j} + \delta
\ln(\lambda) \Big).
\end{equation}
Then, it follows from equations \eqref{E80}, \eqref{E81}, \eqref{E82} and \eqref{E83} that
\begin{align*}
|x_j|^{p_j} |f(x)| \leq & A \bigg( \lambda^{-\frac{1+d_1}{2}} +
\sum_{k=1}^{q_j} \Big( |x_j|^{k-1} \lambda^{-\frac{d_1+3}{2} + k}
\Big) + |x_j|^{q_j} \Big( 1 + (1 - \delta) \lambda^{\frac{-1-d_1}{2} +
 q_j} \\ + & \delta \ln(\lambda) \Big) \bigg).
\end{align*}
However, since $|x| < 1$, we can set $\lambda = \frac{1}{|x_j|} > 1$
to obtain the estimate of $f$ near the origin in the direction $x_j$:
\begin{equation}\label{E84}
|x_j|^{N-\frac{5}{2}+\frac{d_1}{2}+d_\perp} |f(x)| \leq A \Big( 1 +
\delta_{N,2} \delta_{d_1,1} \delta_{d_\perp,0} |\ln(|x_j|)| \Big).
\end{equation}

\begin{step}\label{S23}
Proof of the local estimate \eqref{E18}.
\end{step}

Estimate \eqref{E18} follows from Steps \ref{S21} and
\ref{S22}. Indeed, assume first that $N=2$, $d_1=1$ and
$d_\perp=0$. Since the restriction of $f$ to $B(0,1) \setminus \{ 0
\}$ is continuous, equations \eqref{E79} and \eqref{E84} yield for
every $x \in B(0,1) \setminus \{ 0 \}$,
$$|f(x)| \leq A \min \{ |\ln(|x_1|)|, |\ln(|x_2|)| \} \leq A
|\ln(|x|)|,$$
which is exactly estimate \eqref{E18}. Likewise, if $N \neq 2$, $d_1
\neq 1$ or $d_\perp \neq 0$, the restriction of $f$ to $B(0,1)
\setminus \{ 0 \}$ is also continuous, and equations \eqref{E79} and
\eqref{E84} give for every $x \in B(0,1) \setminus \{ 0 \}$,
\begin{align*}
\Big( x_1^2 + |x_\perp| \Big)^{N - \frac{5}{2} + \frac{d_1}{2} + d_\perp} |f(x)| \leq & A \Big( |x_1|^{2N-5+d_1+2d_\perp} + \sum_{j=2}^N |x_j|^{N-\frac{5}{2}+\frac{d_1}{2}+d_\perp} \Big) |f(x)| \\ \leq & A,
\end{align*}
which concludes the proof of estimate \eqref{E18}.

\begin{step}\label{S24}
$L^q$-integrability of $f$ near the origin.
\end{step}

We now turn to the integral estimates of $f$ near the origin. They follow from the local estimates \eqref{E18} by a standard argument of distribution theory. Indeed, let $j \in \{ 1, \ldots, N \}$ and let $f_j$ and $g_j$ denote the tempered distributions defined on $\R^N$ and $\R^N \setminus \{ 0 \}$ by
\begin{equation}\label{E85}
\begin{split}
f_j & = x_j^{\delta'} f, \\ g_j & = x_j^{\delta'} f_{|\R^N \setminus \{
 0 \}}.
\end{split}
\end{equation}
Here, $\delta'$ is equal to $1$ if $(d_1,d_\perp) = (2,1)$ or
$(d_1,d_\perp) = (4,0)$, and $0$, otherwise, while $f_{|\R^N \setminus
 \{ 0 \}}$ denotes the restriction of $f$ to $\R^N \setminus \{ 0
\}$. As mentioned above, $f_{|\R^N \setminus \{ 0 \}}$ is actually
continuous on $\R^N \setminus \{ 0 \}$, so $g_j$ is also continuous
on $\R^N \setminus \{ 0 \}$. Moreover, equation \eqref{E18} gives
\begin{equation}\label{E86}
|g_j(x)| \leq A \frac{(1 + \delta |\ln(|x|)|)
 |x|^{\delta'}}{(|x_1|^2 + |x_\perp|)^s},
\end{equation}
for every $x \in \R^N \setminus \{ 0 \}$ (with $\delta = \delta_{N,2} \delta_{d_1,1} \delta_{d_\perp,0}$ and $s = N - \frac{5}{2} + \frac{d_1}{2} + d_\perp$). Thus, if $N=2$, $d_1=1$ and $d_\perp=0$, we compute for every $q \geq 1$,
$$\int_{B(0,1)} |g_j(x)|^q dx \leq A \int_{B(0,1)} |\ln(|x|)|^q dx < + \infty,$$
so the distribution $g_j$ belongs to $L^q(B(0,1))$ for every $q \geq 1$, which is the desired result.

On the other hand, if $N \neq 2$, $d_1 \neq 1$ or $d_\perp \neq 0$,
equation \eqref{E86} leads to
$$\int_{B(0,1)} |g_j(x)|^q dx \leq A \int_{B(0,1)} \frac{|x|^{q
 \delta'}}{(|x_1|^2 + |x_\perp|)^{q s}} dx,$$
for every $q \geq 1$, which gives
\begin{align*}
\int_{B(0,1)} |g_j(x)|^q dx \leq & A \int_0^1 r^{N-1+q(\delta' - s)}
\bigg( \int_0^{\frac{\pi}{2}} \frac{\sin(\theta_1)^{N-2} d\theta_1}{(r
 \cos^2(\theta_1) + \sin(\theta_1))^{q s}} \bigg) dr \\ \leq & A
\int_0^1 r^{2N-2+q(\delta' - 2s)} \bigg( \int_0^{+ \infty}
\frac{u^{N-2} (1 + r^2 u^2)^{qs - \frac{N}{2}} du}{(1 + u \sqrt{1 +
 r^2 u^2})^{q s}} \bigg) dr.
\end{align*}
by using the spherical coordinates and the change of variables $u = \frac{\tan(\theta_1)}{r}$. However, we compute for every $r \in ]0,1]$,
$$\int_0^{+ \infty} \frac{u^{N-2} (1 + r^2 u^2)^{qs - \frac{N}{2}}
 du}{(1 + u \sqrt{1 + r^2 u^2})^{q s}} \leq A \bigg( 1 +
\int_1^{\frac{1}{r}} \frac{du}{u^{qs - N + 2}} + r^{qs - N}
\int_\frac{1}{r}^{+ \infty} \frac{du}{u^2} \bigg),$$
so
$$\int_0^{+ \infty} \frac{u^\frac{N}{2} (1 + r^2 u^2)^{qs -
 \frac{N}{2}} du}{(1 + u \sqrt{1 + r^2 u^2})^{q s}} \leq \left\{
 \begin{matrix} A & {\rm if} & qs > N-1, \\ A |\ln(r)| & {\rm if} &
 qs = N-1, \\ A r^{qs - N + 1} & {\rm if} & qs < N-1. \end{matrix}
\right.$$
Thus, we deduce
$$\int_{B(0,1)} |g_j(x)|^q dx \leq A \left\{ \begin{matrix} \int_0^1
 r^{2N-2+q\delta'-2qs} dr & {\rm if} & qs > N-1, \\ \int_0^1
 r^{q\delta'} |\ln(r)| dr & {\rm if} & qs = N-1, \\ \int_0^1
 r^{N+q\delta'-qs-1} dr & {\rm if} & qs < N-1. \end{matrix}
\right.$$
Hence, the distribution $g_j$ belongs to $L^q(B(0,1))$ for
$$1 \leq q < \frac{2N-1}{2s-\delta'},$$
which is exactly condition \eqref{E19} if $\delta' = 0$, and condition
\eqref{E20} if $\delta' = 1$.

We now deduce the same result for $f_j$. Indeed, $g_j$ is continuous
on $\R^N \setminus \{ 0 \}$, and belongs to $L^1(B(0,1))$. Therefore, it defines a distribution on the whole space $\R^N$ (and not only on the subset $\R^N \setminus \{ 0 \}$). Then, it follows from definition \eqref{E85} that the support of the distribution $f_j - g_j$ is included in the singleton $\{ 0 \}$. By Schwartz' theorem, there exist some integer $M$ and some real numbers $\lambda_\alpha$ such
that
\begin{equation}\label{E87}
f_j - g_j = \sum_{|\alpha| \leq M} \lambda_\alpha \partial^\alpha \delta_0.
\end{equation}
However, by Proposition \ref{Paumond1}, definition \eqref{E85} and the proof above, the distribution $g_j$ belongs to $L^q(B(0,1))$ for every $1 \leq q < \frac{2N-1}{2N-5+d_1+2d_\perp - \delta'}$, and to
$L^q(B(0,1)^c)$ for every $q > \frac{N}{N-2+d-\delta'}$. Therefore, it belongs to $L^q(\R^N)$ provided
$$\frac{N}{N-2+d-\delta'} < q < \frac{2N-1}{2N-5+d_1+2d_\perp-
 \delta'}.$$
In particular, if $N > 4 - 2d$, $g_j$ belongs to some space
$L^q(\R^N)$ with $1 \leq q \leq 2$. Therefore, its Fourier transform
$\widehat{g_j}$ belongs to $L^{q'}(\R^N)$ (where $q' = \frac{q}{q-1}$). Then, since equation \eqref{E87} gives for almost every $\xi \in \R^N$,
$$i \partial_j^{\delta'} \bigg( \frac{\underset{l=1}{\overset{N}{\Pi}} \xi_l^{d_l}}{|\xi|^2 + \xi_1^4} \bigg) = \widehat{g_j}(\xi) + \sum_{|\alpha| \leq M} \lambda_\alpha i^{|\alpha|} \xi^\alpha,$$
the distribution $g_j$ belongs to $L^{q'}(\R^N)$ if and only if all the real numbers $\lambda_\alpha$ vanish. Finally, by equation \eqref{E87}, the distribution $f_j$ is equal to $g_j$ if $N > 4 - 2d$. Thus, $f_j$ belongs to $L^q(B(0,1))$ when condition \eqref{E19} (if $\delta' = 0$) or condition \eqref{E20} (if $\delta' = 1$) holds.

Assume now that $N \leq 4 - 2d$. In this case, the coefficient
$\delta'$ is necessarily equal to $0$, and we can simplify our notation by letting $g$ denote the distributions $g_j$ above. Moreover, let $f'$ denote the tempered distribution whose Fourier transform is
\begin{equation}\label{E88}
\widehat{f'}(\xi) = \frac{\underset{j=1}{\overset{N}{\Pi}}
 \xi_j^{d_j'}}{|\xi|^2 + \xi_1^4},
\end{equation}
where $d_j' = d_j + 2 \delta_{j,1}$ for every $j \in \{1,\ldots,N \}$, and $g'$, the restriction of $f'$ to the set $\R^N \setminus \{ 0 \}$. On one hand, since $d_1' \geq 2$ and $N \leq 4 - 2d$, we obtain that $d' \neq 0$ if $N =2$, $d_1' + 2 d_\perp' \leq 4$ and $N > 4 - 2d'$ (with the usual notation $d' = d_1' + \ldots + d_N' = d_1' + d_\perp'$). Hence, by the argument above, the distribution $f'$ is equal to $g'$. On the other hand, by equation \eqref{E88}, the distribution $f'$ is equal to $-\partial^2_1 f$. Likewise, by
definition \eqref{E85}, $g'$ is equal to $-\partial^2_1 g$. Thus,
equation \eqref{E87} leads to
$$0 = f' - g' = \sum_{|\alpha| \leq M} \lambda_\alpha \partial_1^2
\partial^\alpha \delta_\alpha.$$
It follows that all the real numbers $\lambda_\alpha$
vanish. Therefore, the distribution $f$ is equal to $g$: it belongs to $L^q(B(0,1))$ when condition \eqref{E19} holds, which completes the analysis of the $L^q$-integrability of $f$ near the origin and the proof of Proposition \ref{P2}.
\end{proof}

Finally, Theorem \ref{T4} and Corollary \ref{C1} follow from
Proposition \ref{P2}.

\def\proof{\par{\it Proof of Theorem \ref{T4}}. \ignorespaces}
\begin{proof}
Theorem \ref{T4} is a direct consequence of equation \eqref{E18} with $(d_1,d_\perp) = (1,0)$ for $H_0$, $(d_1,d_\perp) = (2,0)$ for $K_0$, and $(d_1,d_\perp) = (2+\delta_{k,1}, 1-\delta_{k,1})$ for $K_k$.
\end{proof}

\def\proof{\par{\it Proof of Corollary \ref{C1}}. \ignorespaces}
\begin{proof}
Likewise, Corollary \ref{C1} follows from conditions \eqref{E19} and \eqref{E20} with $(d_1,d_\perp) = (1,0)$ for $H_0$, $(d_1,d_\perp) = (2,0)$ for $K_0$, and $(d_1,d_\perp) = (2+\delta_{k,1}, 1-\delta_{k,1})$ for $K_k$.
\end{proof}

\subsection{Pointwise limit at infinity of the kernel $K_0$}

This section is devoted to the proofs of Theorems \ref{T5} and \ref{T6}. Since the proof of Theorem \ref{T5} requires formula \eqref{E39} of Theorem \ref{T6}, we first show Theorem \ref{T6}. However, their proofs follow from the same argument: the use of explicit integral expressions. Indeed, we apply Proposition \ref{P5} to $K_0$ to compute formula \eqref{E36}, while we apply Lemma \ref{L2} to $R_{1,1}$ to get formula \eqref{E33}. Then, we make use of the dominated convergence theorem to compute the limits at infinity of such expressions, and to prove that they are equal. Finally, we obtain an explicit expression of this limit by formulae \eqref{E40} and \eqref{E41}.

\def\proof{\par{\it Proof of Theorem \ref{T6}}. \ignorespaces}
\begin{proof}
As mentioned above, our argument relies on expression \eqref{E33}. Indeed, the kernel $R_{1,1}$ satisfies all the assumptions of Lemma \ref{L2}. By formula \eqref{E38}, its Fourier transform is bounded on $\R^N \setminus \{ 0 \}$ and belongs to $C^\infty(\R^N \setminus \{ 0 \})$. Therefore, assumption $(i)$ holds with $r=s=0$. Moreover, we compute for every $j \in \{ 1,\ldots,N \}$ and $n \in \N$,
$$\forall \xi \in \R^N \setminus \{ 0 \}, |\partial_j^n
\widehat{R_{1,1}}(\xi)| \leq \frac{A}{|\xi|^n}.$$
Hence, assumptions $(ii)$ and $(iii)$ hold with $p=N$ and $m=N+1$. Moreover, the function $\partial_j^n \widehat{R_{1,1}}$ belongs to $L^q(B(0,1)^c)$ for any $q > q_{N+1-n} = \frac{N}{n}$. In particular, we notice that
$$\forall n \in \{ 1,\ldots,N-1 \}, q_n = \frac{N}{N+1-n} <
\frac{N}{N-n}.$$
Thus, the kernel $R_{1,1}$ also verifies assumption $(iv)$ and we can write by Lemma \ref{L2} for every $\lambda > 0$ and $x \in \Omega_j = \{ x \in \R^N, x_j \neq 0 \}$,
\begin{equation}\label{E89}
\begin{split}
x_j^N R_{1,1}(x) & = \frac{i^N}{(2 \pi)^N} \bigg( \frac{i}{x_j} \int_{B(0,\lambda)^c} \partial_j^{N+1} \widehat{R_{1,1}}(\xi) e^{ix.\xi} d\xi + \frac{i}{\lambda x_j} \int_{S(0,\lambda)} \xi_j \partial_j^N \widehat{R_{1,1}}(\xi) \\ & e^{ix.\xi} d\xi + \frac{1}{\lambda} \int_{S(0,\lambda)} \xi_j \partial_j^{N-1} \widehat{R_{1,1}}(\xi) d\xi + \int_{B(0,\lambda)} \partial_j^N \widehat{R_{1,1}}(\xi) (e^{ix.\xi} - 1) d\xi \bigg).
\end{split}
\end{equation}
On the other hand, by formulae \eqref{E40} and \eqref{E41}, the
restriction of $R_{1,1}$ to $\R^N \setminus \{ 0 \}$ is given by
$$\forall x \in \R^N \setminus \{ 0 \}, R_{1,1}(x) =
\frac{\Gamma(\frac{N}{2})}{2 \pi^\frac{N}{2}} \frac{|x|^2 - N
 x_1^2}{|x|^{N+2}},$$
so equation \eqref{E89} becomes for every $x \in \Omega_j$,
\begin{align*}
& \frac{\Gamma(\frac{N}{2})}{2 \pi^\frac{N}{2}} \frac{|x|^2 - N x_1^2}{|x|^{N+2}} = \frac{i^N}{(2 \pi x_j)^N} \bigg( \frac{i}{x_j} \int_{B(0,\lambda)^c} \partial_j^{N+1} \widehat{R_{1,1}}(\xi) e^{ix.\xi} d\xi + \frac{i}{\lambda x_j} \int_{S(0,\lambda)} \xi_j e^{ix.\xi} \\ & \partial_j^N \widehat{R_{1,1}}(\xi) d\xi + \frac{1}{\lambda} \int_{S(0,\lambda)} \xi_j \partial_j^{N-1} \widehat{R_{1,1}}(\xi) d\xi + \int_{B(0,\lambda)} \partial_j^N \widehat{R_{1,1}}(\xi) (e^{ix.\xi} - 1) d\xi \bigg).
\end{align*}
In particular, if we consider $y \in \R^N$, $\sigma \in \S^{N-1}$, with $\sigma_j \neq 0$), and $R > 0$ (with $R |\sigma_j| > 2 |y_j|$), we compute for $x = R\sigma - y$ and $\lambda = \frac{1}{R}$ after the change of variables $u = R\xi$,
\begin{equation}\label{E90}
\begin{split}
& \frac{\Gamma(\frac{N}{2})}{2 \pi^\frac{N}{2}} \frac{|\sigma - \frac{y}{R}|^2 - N (\sigma_1 - \frac{y_1}{R})^2}{|\sigma - \frac{y}{R}|^{N+2}} = \frac{i^N}{(2 \pi (\sigma_j - \frac{y_j}{R}))^N} \bigg( \frac{i}{\sigma_j - \frac{y_j}{R}} \int_{B(0,1)^c} \frac{\partial_j^{N+1} \widehat{R_{1,1}} (\frac{u}{R})}{R^{N+1}} \\ & e^{i (\sigma - \frac{y}{R}). u} du + \int_{\S^{N-1}} \frac{\partial_j^{N-1} \widehat{R_{1,1}} (\frac{u}{R})}{R^{N-1}} u_j du + \frac{i}{\sigma_j - \frac{y_j}{R}} \int_{\S^{N-1}} \frac{\partial_j^N \widehat{R_{1,1}} (\frac{u}{R})}{R^N} e^{i (\sigma - \frac{y}{R}). u} u_j du \\ & + \int_{B(0,1)} \frac{\partial_j^N \widehat{R_{1,1}} (\frac{u}{R})}{R^N} (e^{i (\sigma - \frac{y}{R}). u} - 1) du \bigg).
\end{split}
\end{equation}
In order to get formula \eqref{E39}, it now remains to compute the
limit when $R$ tends to $+ \infty$ of all the terms of equation
\eqref{E90}. Here, we make use of the homogeneity of
$\widehat{R_{1,1}}$. Indeed, by formula \eqref{E38}, the function
$\widehat{R_{1,1}}$ is a homogeneous rational fraction of degree
$0$. Therefore, the partial derivative $\partial_j^k
\widehat{R_{1,1}}$ is a homogeneous rational fraction of degree
$-k$. Thus, equation \eqref{E90} becomes
\begin{equation}\label{E91}
\begin{split}
& \frac{\Gamma(\frac{N}{2})}{2 \pi^\frac{N}{2}} \frac{|\sigma - \frac{y}{R}|^2 - N (\sigma_1 - \frac{y_1}{R})^2}{|\sigma - \frac{y}{R}|^{N+2}} = \frac{i^N}{(2 \pi (\sigma_j - \frac{y_j}{R}))^N} \bigg( \frac{i}{\sigma_j - \frac{y_j}{R}} \int_{B(0,1)^c} \partial_j^{N+1} \widehat{R_{1,1}} (u) \\ & e^{i (\sigma - \frac{y}{R}). u} du + \int_{\S^{N-1}} \partial_j^{N-1} \widehat{R_{1,1}} (u) u_j du + \frac{i}{\sigma_j - \frac{y_j}{R}} \int_{\S^{N-1}} \partial_j^N \widehat{R_{1,1}} (u) e^{i (\sigma - \frac{y}{R}). u} u_j du \\ & + \int_{B(0,1)} \partial_j^N \widehat{R_{1,1}} (u) (e^{i (\sigma - \frac{y}{R}). u} - 1) du \bigg).
\end{split}
\end{equation}
We now invoke the dominated convergence theorem to compute the
limit of the right-hand side of equation \eqref{E91}. Indeed, by
homogeneity of the partial derivatives of $\widehat{R_{1,1}}$, we
compute for the first term in the right-hand side,
$$\forall u \in B(0,1)^c, \Big| \partial_j^{N+1} \widehat{R_{1,1}}(u) e^{i (\sigma - \frac{y}{R}). u} \Big| \leq \frac{A}{|u|^{N+1}}.$$ 
Likewise, the third term in the right-hand side satisfies
$$\forall u \in \S^{N-1}, \Big| u_j \partial_j^N \widehat{R_{1,1}}(u) e^{i (\sigma - \frac{y}{R}). u} \Big| \leq A,$$
while the fourth term verifies
$$\forall u \in B(0,1), \Big| \partial_j^N \widehat{R_{1,1}}(u) (e^{i (\sigma - \frac{y}{R}). u} - 1) \Big| \leq \frac{A}{|u|^{N-1}}
\Big|\sigma - \frac{y}{R} \Big| \leq \frac{A}{|u|^{N-1}},$$
provided $R \geq 2 |y|$. Therefore, by taking the limit $R \to + \infty$ in equation \eqref{E91}, the dominated convergence theorem
yields
\begin{align*}
& \frac{\Gamma(\frac{N}{2})}{2 \pi^\frac{N}{2}} (1 - N \sigma_1^2) = \frac{i^N}{(2 \pi \sigma_j)^N} \bigg( \frac{i}{\sigma_j} \bigg( \int_{B(0,1)^c} \partial_j^{N+1} \widehat{R_{1,1}}(u) e^{i \sigma. u} du + \int_{\S^{N-1}} \partial_j^N \widehat{R_{1,1}}(u) \\ & e^{i \sigma. u} u_j du \bigg) + \int_{\S^{N-1}} u_j \partial_j^{N-1} \widehat{R_{1,1}}(u) du + \int_{B(0,1)} \partial_j^N \widehat{R_{1,1}}(u) (e^{i \sigma. u} - 1) du \bigg),
\end{align*}
which is exactly formula \eqref{E39}.
\end{proof}

Then, we deduce the pointwise limit of $K_0$ given by Theorem \ref{T5}.

\def\proof{\par{\it Proof of Theorem \ref{T5}}. \ignorespaces}
\begin{proof}
Let $\sigma \in \S^{N-1}$ and $y \in \R^N$, and consider some $j \in \{ 1,\ldots,N \}$ such that $\sigma_j \neq 0$ and some $R > \max \{ 2|y|, \frac{2 |y_j|}{|\sigma_j|} \}$. The kernel $K_0$ fulfils all the assumptions of Proposition \ref{P5} with $d_1=2$ and $d_\perp=0$. Therefore, by equation \eqref{E34}, formula \eqref{E36} holds with $m_1 = 2N-2$ and $m_j=N$ if $j \geq 2$. In particular, after the change of variables $u = R\xi$, this formula becomes for $x = R\sigma - y$ and $\lambda = \frac{1}{R}$,
\begin{equation}\label{E92}
\begin{split}
& R^N K_0(R\sigma - y) = \frac{i^N}{(2 \pi (\sigma_j - \frac{y_j}{R}))^N} \bigg( \Big( -i \Big( \sigma_j - \frac{y_j}{R} \Big) \Big)^{N - m_j} \int_{B(0,1)^c} \frac{\partial_j^{m_j} \widehat{K_0} (\frac{u}{R})}{R^{m_j}} \\ & e^{i (\sigma - \frac{y}{R}). u} du + \sum_{k=N}^{m_j-1} \Big( -i \Big( \sigma_j - \frac{y_j}{R} \Big) \Big)^{N - k - 1} \int_{\S^{N-1}} \frac{\partial_j^k \widehat{K_0} (\frac{u}{R})}{R^k} e^{i (\sigma - \frac{y}{R}). u} u_j du \\ & + \int_{\S^{N-1}} \frac{\partial_j^{N-1} \widehat{K_0} (\frac{u}{R})}{R^{N-1}} u_j du + \int_{B(0,1)} \frac{\partial_j^N \widehat{K_0} (\frac{u}{R})}{R^N} (e^{i (\sigma - \frac{y}{R}). u} - 1) du \bigg).
\end{split}
\end{equation}
Then, denoting $N_j=\max\{ N+1, m_j \}$, we claim that
\begin{equation}\label{E93}
\begin{split}
& R^N K_0(R\sigma - y) = \frac{i^N}{(2 \pi (\sigma_j - \frac{y_j}{R}))^N} \bigg( \Big( -i \Big( \sigma_j - \frac{y_j}{R} \Big) \Big)^{N - N_j} \int_{B(0,1)^c} \frac{\partial_j^{N_j} \widehat{K_0} (\frac{u}{R})}{R^{N_j}} \\ & e^{i (\sigma - \frac{y}{R}). u} du + \sum_{k=N}^{N_j-1} \Big( -i \Big( \sigma_j - \frac{y_j}{R} \Big) \Big)^{N - k - 1} \int_{\S^{N-1}} \frac{\partial_j^k \widehat{K_0} (\frac{u}{R} )}{R^k} e^{i (\sigma - \frac{y}{R}). u} u_j du \\ & + \int_{\S^{N-1}} \frac{\partial_j^{N-1} \widehat{K_0} (\frac{u}{R})}{R^{N-1}} u_j du + \int_{B(0,1)} \frac{\partial_j^N \widehat{K_0} (\frac{u}{R})}{R^N} (e^{i (\sigma - \frac{y}{R}). u} - 1) du \bigg).
\end{split}
\end{equation}
Indeed, by Proposition \ref{P4} (with $d_1 = 2$ and $d_\perp = 0$),
the function $\partial_j^k\widehat{K_0}$ belongs to $L^1(B(0,1)^c)$ for
every $k \geq m_j$. Thus, equation \eqref{E93} results from equation
\eqref{E92} by several integrations by parts of the integral on
$B(0,1)^c$.

In order to get formula \eqref{E35}, it now remains to compute the limit when $R$ tends to $+ \infty$ of equation \eqref{E93}. In particular, we must compute the limit when $R \to + \infty$ of the
functions $u \mapsto R^{-k} \partial_j^{k} \widehat{K_0} (\frac{u}{R})$ for any $k \in \{ N-1,\ldots,N_j \}$. Thus, we must
describe more precisely $\partial_j^{k} \widehat{K_0}$. Indeed, the
Fourier transform $\widehat{K_0}$ satisfies all the assumptions of
Proposition \ref{P3} with $d_1=2$ and $d_\perp=0$. Therefore, by
equations \eqref{E21}, \eqref{E23} and \eqref{E25}, the partial
derivative $\partial_j^k \widehat{K_0}$ is equal to
\begin{equation}\label{E94}
\forall k \in \N, \forall \xi \in \R^N \setminus \{ 0 \}, \partial_j^k
\widehat{K_0}(\xi) = \frac{P_{j,k}(\xi)}{(|\xi|^2 + \xi_1^4)^{k+1}},
\end{equation}
where $P_{j,k}$ is a polynomial function on $\R^N$ which satisfies,
\begin{equation}\label{E95}
\forall \xi \in B(0,1), |P_{j,k}(\xi)| \leq A_k |\xi|^{k + 2}.
\end{equation}
Moreover, denoting
\begin{equation}\label{E96}
\forall \xi \in \R^N, P_{j,k}(\xi) = \underset{i =
 0}{\overset{d_k}{\sum}} Q_i(\xi),
\end{equation}
where $Q_i$ is a homogeneous polynomial function either equal to $0$
or of degree $i$, we claim that the non-vanishing function $Q_{i_0}$
of lowest degree is equal to
\begin{equation}\label{E97}
\forall \xi \in \R^N, Q_{i_0}(\xi) = |\xi|^{2k+2} \partial_j^k
\widehat{R_{1,1}}(\xi).
\end{equation}
Indeed, the partial derivative $\partial_j^k \widehat{R_{1,1}}(\xi)$
satisfies by induction,
$$\forall \xi \in \R^N \setminus \{ 0 \}, \partial_j^k
\widehat{R_{1,1}}(\xi) = \frac{S_k(\xi)}{|\xi|^{2k+2}},$$
where the functions $S_k$ are polynomial functions given by
\begin{align*}
S_0(\xi) = & \xi_1^2, \\ S_{k+1}(\xi) = & |\xi|^2 \partial_j S_k(\xi)
- 2(k+1) \xi_j S_k(\xi).
\end{align*}
In particular, the polynomial function $S_k$ is of degree $k +
2$. Likewise, by equations \eqref{E52}, the functions $P_{j,k}$ are
given by
\begin{align*}
P_{j,0}(\xi) = & \xi_1^2, \\ P_{j,k+1}(\xi) = & (|\xi|^2 + \xi_1^4)
\partial_j P_{j,k}(\xi) - 2 (k+1) (\xi_j + 2 \delta_{j,1} \xi_1^3)
P_{j,k}(\xi).
\end{align*}
Therefore, the homogeneous term of lowest degree $Q_{i_0}$ of
$P_{j,k}$ satisfies exactly the same equation as $S_k$. Hence, the
function $Q_{i_0}$ is equal to $S_k$, which leads to formula
\eqref{E97}. In particular, the degree $i_0$ of $Q_{i_0}$ is equal to
$k+2$.

We now compute the limit of the function $u \mapsto R^{-k} \partial_j^{k} \widehat{K_0} (\frac{u}{R})$ when $R$ tends to $+
\infty$. Indeed, by equations \eqref{E94} and \eqref{E96}, we compute for every $k \in \N$ and $u \in \R^N \setminus \{ 0 \}$,
$$R^{-k} \partial^k_j \widehat{K_0} \Big( \frac{u}{R} \Big) =
\frac{\underset{i = k+2}{\overset{d_k}{\sum}} R^{k+2-i} Q_i(u)}{(|u|^2 + R^{-2} u_1^4)^{k+1}},$$
which gives by equation \eqref{E97},
\begin{equation}\label{E98}
R^{-k} \partial^k_j \widehat{K_0} \Big( \frac{u}{R} \Big)
\underset{R \to + \infty}{\to} \frac{Q_{k+2}(u)}{|u|^{2k+2}} =
\partial_j^k \widehat{R_{1,1}}(u).
\end{equation}

Finally, we invoke the dominated convergence theorem to compute the
limit of the right-hand side of equation \eqref{E93}. Indeed, the first term in the right-hand side of equation \eqref{E93} becomes by the change of variables $u = R\xi$,
\begin{align*}
\int_{B(0,1)^c} R^{-N_j} \partial_j^{N_j} \widehat{K_0}
\Big(\frac{u}{R} \Big) e^{i (\sigma - \frac{y}{R}). u} du = & \int_{1
 < |u| < R} R^{-N_j} \partial_j^{N_j} \widehat{K_0} \Big(\frac{u}{R}
\Big) e^{i (\sigma - \frac{y}{R}). u} du \\ + & \int_{B(0,1)^c}
R^{N-N_j} \partial_j^{N_j} \widehat{K_0} (\xi) e^{i (R\sigma -
 y). \xi} d\xi.
\end{align*}
On one hand, by Proposition \ref{P5}, the function $\partial_j^{N_j}
\widehat{K_0}$ belongs to $L^1(B(0,1)^c)$, so
$$\int_{B(0,1)^c} R^{N-N_j} \partial_j^{N_j} \widehat{K_0} (\xi) e^{i (R\sigma - \xi). u} du \leq R^{N - N_j} \int_{B(0,1)^c}
|\partial_j^{N_j} \widehat{K_0}(\xi)| d\xi \underset{R \to +
 \infty}{\to} 0.$$
On the other hand, by equations \eqref{E94} and \eqref{E95},
$$\forall 1 < |u| < R, \Big| R^{-N_j} \partial_j^{N_j} \widehat{K_0}
\Big( \frac{u}{R} \Big) e^{i (\sigma - \frac{y}{R}). u} \Big| \leq
\frac{A}{|u|^{N_j}},$$
so, since $N_j \geq N+1$, by the dominated convergence theorem and
assertion \eqref{E98},
$$\int_{1 < |u| < R} R^{-N_j} \partial_j^{N_j} \widehat{K_0}
\Big(\frac{u}{R} \Big) e^{i (\sigma - \frac{y}{R}). u} du \underset{R
 \to + \infty}{\to} \int_{B(0,1)^c} \partial_j^{N_j}
\widehat{R_{1,1}} (u) e^{i \sigma. u} du.$$
Hence, we deduce
\begin{equation}\label{E99}
\int_{B(0,1)^c} R^{-N_j} \partial_j^{N_j} \widehat{K_0} \Big(
\frac{u}{R} \Big) e^{i (\sigma - \frac{y}{R}). u} du \underset{R \to +
 \infty}{\to} \int_{B(0,1)^c} \partial_j^{N_j} \widehat{R_{1,1}} (u)
e^{i \sigma. u} du.
\end{equation}
Likewise, equations \eqref{E94} and \eqref{E95} give
$$\forall u \in \S^{N-1}, \Big| u_j R^{-k} \partial_j^k \widehat{K_0} \Big( \frac{u}{R} \Big) e^{i (\sigma - \frac{y}{R}). u} \Big| \leq A,$$
for every $k \in \{ N,\ldots,N_j-1 \}$, and
$$\forall u \in \S^{N-1}, \Big| u_j R^{1-N} \partial_j^{N-1}
\widehat{K_0} \Big( \frac{u}{R} \Big) \Big| \leq A,$$
while the last term in the right-hand side of equation \eqref{E93} verifies
$$\forall u \in B(0,1), \Big| R^{-N} \partial_j^N \widehat{K_0} \Big( \frac{u}{R} \Big) (e^{i (\sigma - \frac{y}{R}). u} - 1) \Big| \leq \frac{A}{|u|^{N-1}} \Big|\sigma - \frac{y}{R} \Big| \leq
\frac{A}{|u|^{N-1}},$$
provided $R \geq 2 |y|$. Hence, by assertion \eqref{E98}, and the dominated convergence theorem,
\begin{align*}
& \sum_{k=N}^{N_j-1} \Big( -i \Big( \sigma_j - \frac{y_j}{R} \Big)
\Big)^{N - k - 1} \int_{\S^{N-1}} \frac{\partial_j^k \widehat{K_0}
 (\frac{u}{R})}{R^k} e^{i (\sigma - \frac{y}{R}). u} u_j du \\ + &
\int_{\S^{N-1}} \frac{\partial_j^{N-1} \widehat{K_0}
 (\frac{u}{R})}{R^{N-1}} u_j du + \int_{B(0,1)} \frac{\partial_j^N
 \widehat{K_0} (\frac{u}{R})}{R^{N-1}} (e^{i (\sigma -
 \frac{y}{R}). u} - 1) du \\ \underset{R \to +
 \infty}{\to} & \sum_{k=N}^{N_j-1} (-i \sigma_j)^{N - k - 1}
\int_{\S^{N-1}} \partial_j^k \widehat{R_{1,1}} (u) e^{i \sigma. u} u_j du + \int_{\S^{N-1}} \partial_j^{N-1} \widehat{R_{1,1}}(u) u_j du \\ + & \int_{B(0,1)} \partial_j^N \widehat{R_{1,1}}(u) (e^{i
 \sigma. u} - 1) du.
\end{align*}
Thus, it follows from equations \eqref{E93} and \eqref{E99} that
\begin{equation}\label{E100}
\begin{split}
R^N K_0(R\sigma - y) \underset{R \to + \infty}{\to} & \frac{i^N}{(2
 \pi \sigma_j)^N} \bigg( (-i \sigma_j)^{N - N_j} \int_{B(0,1)^c}
\partial_j^{N_j} \widehat{R_{1,1}} (u) e^{i \sigma. u} du \\ + &
\sum_{k=N}^{N_j-1} (-i \sigma_j)^{N - k - 1} \int_{\S^{N-1}}
\partial_j^k \widehat{R_{1,1}} (u) e^{i \sigma. u} u_j du \\ + &
\int_{\S^{N-1}} \partial_j^{N-1} \widehat{R_{1,1}}(u) u_j du +
\int_{B(0,1)} \partial_j^N \widehat{R_{1,1}}(u) (e^{i \sigma. u} - 1) du.
\end{split}
\end{equation}
Then, in order to obtain assertion \eqref{E35}, we integrate by parts the first term in the right-hand side of equation \eqref{E100}. Indeed, the partial derivative $\partial_j^k \widehat{R_{1,1}}$ belongs to $L^1(B(0,1)^c)$ for every $k \geq N+1$, so, after several integrations by parts, equation \eqref{E100} becomes
\begin{align*}
R^N K_0(R\sigma - y) \underset{R \to + \infty}{\to} & \frac{i^N}{(2 \pi \sigma_j)^N} \bigg( \frac{i}{\sigma_j} \int_{B(0,1)^c} \partial_j^{N+1} \widehat{R_{1,1}} (u) e^{i \sigma. u} du + \frac{i}{\sigma_j} \int_{\S^{N-1}} e^{i \sigma. u} u_j \\ & \partial_j^N \widehat{R_{1,1}} (u) du + \int_{\S^{N-1}} \partial_j^{N-1} \widehat{R_{1,1}}(u) u_j du + \int_{B(0,1)} \partial_j^N \widehat{R_{1,1}}(u) \\ & (e^{i \sigma. u} - 1) du \bigg).
\end{align*}
Then, assertion \eqref{E35} holds by formula \eqref{E39}, which
completes the proof of Theorem \ref{T5}.
\end{proof}

\subsection{Rigorous derivation of the convolution equations}

We now aim at rigorously deriving the convolution equations \eqref{E10}, \eqref{E11} and \eqref{E42}. Indeed, our analysis of the asymptotics of the solitary waves relies on these equations.

In the introduction, we already proved that equations \eqref{E10} and \eqref{E11} hold almost everywhere. However, we will also consider the gradient of equation \eqref{E11}, whose derivation is rather more difficult. In order to give it a rigorous sense, we establish equation \eqref{E42}, which holds for smooth functions $f$ with sufficient decay at infinity. In particular, the function $v^{p+1}$ satisfy such assumptions by Theorems \ref{T7} and \ref{T8}, which gives a rigorous sense to the gradient of equation \eqref{E11}.

Let us now establish Lemma \ref{L3}.

\def\proof{\par{\it Proof of Lemma \ref{L3}}. \ignorespaces}
\begin{proof}
Let $k \in \{ 1,\ldots,N \}$ and let us consider the function $h_k$ given by
\begin{equation}\label{E101}
\begin{split}
\forall x \in \R^N, h_k(x) = & i \int_{B(0,1)^c} K_k(y) f(x-y) dy + i
\int_{B(0,1)} K_k(y) (f(x-y) - f(x)) dy \\ + & \bigg( \int_{\S^{N-1}}
K_0(y) y_k dy \bigg) f(x).
\end{split}
\end{equation}
Our proof splits into two steps. The first one states the continuity of $g$, $h_1$, $\ldots$, $h_N$ on $\R^N$, while the second one establishes that the partial derivative $\partial_k g$ of $g$ in the sense of distributions is equal to $h_k$. Then, we conclude that the function $g$ is of class $C^1$ on $\R^N$ and that its first order partial derivatives are given by formula \eqref{E42}.

\setcounter{step}{0}
\begin{step}\label{S31}
Continuity of the functions $g$, $h_1$, $\ldots$, $h_N$.
\end{step}

The function $f$ is continuous on $\R^N$, so, by assumption $(i)$, it belongs to $L^q(\R^N)$ for every $q \geq 1$. It follows from Young's inequalities that the function $g$ is well-defined in $L^q(\R^N)$ for every $1 < q < \frac{2N - 1}{2N - 3}$. In particular, it is given for almost every $x \in \R^N$ by
\begin{equation}\label{E102}
g(x) = \int_{\R^N} K_0(y) f(x-y) dy.
\end{equation}
Moreover, the kernel $K_0$ belongs to $L^1(B(0,1))$ by Corollary \ref{C1}. Therefore, by continuity of $f$, assumption $(i)$ and a standard application of the dominated convergence theorem, the function $g_1:x \mapsto \int_{B(0,1)} K_0(y) f(x-y) dy$ is continuous on $\R^N$. Likewise, the function $g_2:x \mapsto \int_{B(0,1)^c} K_0(y) f(x-y) dy$ is continuous on $\R^N$. Indeed, consider $x_0 \in \R^N$ and compute by assumption $(i)$ and Theorem \ref{T3},
$$\forall x \in B(x_0,1), \forall y \in B(0,1)^c, |K_0(y) f(x-y)| \leq \frac{A}{|y|^N (1 + |x-y|^{N(p+1)})}.$$
This gives
$$\forall 1 < |y| < |x_0| + 1, |K_0(y) f(x-y)| \leq A_{x_0},$$
and
$$\forall y \in B(0,|x_0| + 1)^c, |K_0(y) f(x-y)| \leq \frac{A}{|y|^N (1 + (|y| - |x_0| - 1)^{N(p+1)})}.$$
Thus, by a standard application of the dominated convergence theorem, $g_2$ is continuous at the point $x_0$. Hence, it is continuous on $\R^N$, as well as the function $g = g_1 + g_2$.

On the other hand, by Theorem \ref{T3}, the kernel $K_k$ belongs to $L^1(B(0,1)^c)$. Thus, by continuity of $f$, assumption $(i)$ and a standard application of the dominated convergence theorem, the function $h_k^1: x \mapsto i \int_{B(0,1)^c} K_k(y) f(x-y) dy$ is continuous on $\R^N$. Likewise, by Proposition \ref{P5}, the kernel $K_0$ belongs to $C^0(\R^N \setminus \{ 0 \})$, so, by continuity of $f$, the function $h_k^2: x \mapsto ( \int_{\S^{N-1}} K_0(y) y_k dy) f(x)$ is continuous on $\R^N$. Finally, assumption $(ii)$ yields
$$\forall y \in B(0,1), |K_k(y) (f(x-y) - f(x))| \leq \| \nabla f
\|_{L^\infty(\R^N)} \sum_{j=1}^N |y_j K_k(y)|.$$
Since the functions $y \mapsto y_j K_k(y)$ belong to $L^1(B(0,1)$ by
Corollary \ref{C1}, the function $h_k^3: x \mapsto i \int_{B(0,1)}
K_k(y) (f(x-y) - f(x)) dy$ is also continuous on $\R^N$. Hence, the
function $h_k = h_k^1 + h_k^2 + h_k^3$ is continuous on $\R^N$.

\begin{step}\label{S32}
First order partial derivatives of $g$ in the sense of distributions.
\end{step}

Let us now consider some test function $\phi \in C^\infty_c(\R^N)$. By definition \eqref{E102}, we compute
\begin{equation}\label{E103}
<\partial_k g, \phi> = - \int_{\R^N} f(y) \Lambda_{\phi}(y) dy,
\end{equation}
where
\begin{equation}\label{E104}
\forall y \in \R^N, \Lambda_{\phi}(y) = \int_{\R^N} K_0(x-y)
\partial_k \phi(x) dx.
\end{equation}
Then, let us fix $y \in \R^N$. Since $\phi$ belongs to $C^\infty_c(\R^N)$, there is some $R > |y| + 2$ such that
$$\Lambda_{\phi}(y) = \int_{B(y,1)} K_0(x-y) \partial_k (\phi(x) -
\phi(y)) dx + \int_{B(0,R) \setminus B(y,1)} K_0(x-y) \partial_k
\phi(x) dx.$$
However, the kernels $K_0$ and $K_k$ belong to $C^0(\R^N \setminus \{ 0 \})$ by Proposition \ref{P5}. Therefore, since the kernel $K_k$ is equal to $-i \partial_k K_0$, the kernel $K_0$ is of class $C^1$ on $\R^N \setminus \{ 0 \}$ such that
$$\forall z \in \R^N \setminus \{ 0 \}, \partial_k K_0(z) = i
K_k(z).$$
Hence, by integrating by parts,
\begin{equation}\label{E105}
\begin{split}
\Lambda_{\phi}(y) = & \int_{B(y,1)} K_0(x-y) \partial_k (\phi(x) -
\phi(y)) dx - i \int_{B(0,R) \setminus B(y,1)} K_k(x-y) \phi(x) dx \\ - & \int_{S(y,1)} K_0(x-y) (x_k - y_k) \phi(x) dx.
\end{split}
\end{equation}
Moreover, by Corollary \ref{C1}, the kernel $K_0$ belongs to
$L^1(B(0,1))$, so
\begin{equation}\label{E106}
\int_{B(y,1)} K_0(x-y) \partial_k (\phi(x) - \phi(y)) dx =
\underset{\eps \to 0}{\lim} \int_{\eps < |x-y| < 1} K_0(x-y)
\partial_k (\phi(x) - \phi(y)) dx.
\end{equation}
Hence, since $K_0$ is of class $C^1$ on $\R^N \setminus \{ 0 \}$, we
compute for every $\eps \in ]0,1[$,
\begin{equation}\label{E107}
\begin{split}
& \int_{\eps < |x-y| < 1} K_0(x-y) \partial_k (\phi(x) - \phi(y)) dx \\ = & - i \int_{\eps < |x-y| < 1} K_k(x-y) (\phi(x) - \phi(y)) dx + \int_{S(y,1)} K_0(x-y) (\phi(x) - \phi(y)) \\ & (x_k - y_k) dx - \int_{S(y,\eps)} K_0(x-y) (\phi(x) - \phi(y)) \frac{x_k - y_k}{\eps} dx.
\end{split}
\end{equation}
On one hand, $\phi$ belongs to $C^\infty_c(\R^N)$, which gives
$$\forall x \in B(y,1), \Big| K_k(x-y) (\phi(x) - \phi(y)) \Big| \leq A \sum_{l=1}^N |(x_l - y_l) K_k(x-y)|,$$
so, since the functions $x \mapsto x_l K_k(x)$ belong to $L^1(B(0,1))$
by Corollary \ref{C1},
\begin{equation}\label{E108}
\int_{\eps < |x-y| < 1} K_k(x-y) (\phi(x) - \phi(y)) dx
\underset{\eps \to 0}{\to} \int_{B(y,1)} K_k(x-y) (\phi(x) - \phi(y))
dx.
\end{equation}
On the other hand, the kernel $K_0$ belongs to $L^1(B(0,1))$. Therefore, there exists a sequence of positive real numbers $\delta_n$ which tends to $0$ when $n \to + \infty$, and which satisfies
$$\forall n \in \N, \exists \eps_n \in ]0,\delta_n[, \delta_n
\int_{S(0,\eps_n)} |K_0(z)| dz \leq \int_{B(0,\delta_n)} |K_0(z)|
dz.$$
Thus, since $K_0$ and $\phi$ belong respectively to $L^1(B(0,1))$ and
$C^\infty_c(\R^N)$,
\begin{align*}
\bigg| \frac{1}{\eps_n} \int_{S(y,\eps_n)} K_0(x-y) (x_k -
y_k) (\phi(x) - \phi(y)) dx \bigg| \leq & \frac{A}{\eps_n}
\int_{S(y,\eps_n)} |K_0(x-y)| |x - y|^2 dx \\ \leq & A
\int_{B(0,\delta_n)} |K_0(z)| dz \underset{n \to + \infty}{\to} 0.
\end{align*}
Therefore, by equations \eqref{E106}, \eqref{E107} and \eqref{E108},
\begin{align*}
\int_{B(y,1)} K_0(x-y) \partial_k (\phi(x) - \phi(y)) dx = & - i \int_{B(y,1)} K_k(x-y) (\phi(x) - \phi(y)) dx \\ + & \int_{S(y,1)} K_0(x-y) (x_k - y_k) (\phi(x) - \phi(y)) dx,
\end{align*}
which gives by equation \eqref{E105},
\begin{equation}\label{E109}
\begin{split}
\Lambda_{\phi}(y) = - & i \int_{B(y,1)} K_k(x-y) (\phi(x) - \phi(y)) dx - \bigg( \int_{S(y,1)} K_0(x-y) (x_k - y_k) dx \bigg) \\ & \phi(y) - i \int_{B(0,R) \setminus B(y,1)} K_k(x-y) \phi(x) dx.
\end{split}
\end{equation}
On the other hand, by definition \eqref{E104}, $\Lambda_\phi$ belongs to $L^q(\R^N)$ for $1 < q \leq + \infty$. Indeed, if the distance between $y$ and the support of $\phi$ is more than $1$, we compute by Theorem \ref{T3},
$$|\Lambda_\phi(y)| \leq A \int_{{\rm Supp}(\phi)}
\frac{|\phi(x)|}{|x-y|^N} dx \leq \frac{A}{{\rm d}(y, {\rm
 supp}(\phi))^N} \| \phi \|_{L^1(\R^N)},$$
while if this distance is less than $1$,
$$|\Lambda_\phi(y)| \leq A \| \phi \|_{L^\infty(\R^N)} \| K_0
\|_{L^1({\rm Supp}(\phi))},$$
by Corollary \ref{C1}. Thus, since the support of $\phi$ is compact, $\Lambda_\phi$ does
belong to $L^q(\R^N)$ for $1 < q \leq + \infty$. However, the function $f$ belongs to $L^q(\R^N)$ for $1 \leq q \leq + \infty$ by assumption
$(i)$. Therefore, it follows from equations \eqref{E103}, \eqref{E104}
and \eqref{E109}, and Fubini's theorem that
\begin{align*}
< \partial_k g, \phi> = & \int_{\R^N} f(y) \bigg( i \int_{B(y,1)}
K_k(x-y) (\phi(x) - \phi(y)) dx \\ + & \bigg( \int_{S(y,1)} (x_k -
y_k) K_0(x-y) dx \bigg) \phi(y) + i \int_{B(y,1)^c} K_k(x-y) \phi(x)
dx \bigg) dy.
\end{align*}
It now remains to make the changes of variables $z = x-y$ and $t = y+z$ to get
\begin{align*}
< \partial_k g, \phi> = & \int_{\R^N} \bigg( i \int_{B(0,1)} K_k(z)
(f(t-z) - f(t)) dz + \Big( \int_{S(0,1)} K_0(z) z_k dz \Big) f(t) \\ + & i \int_{B(0,1)^c} K_k(z) f(t-z) dz \bigg) \phi(t) dx,
\end{align*}
by Fubini's theorem, which gives
$$< \partial_k g, \phi> = <h_k, \phi>,$$
by definition \eqref{E101}. Therefore, the partial derivative $\partial_k g$ of $g$ in the sense of distributions is equal to $h_k$. Since $g$, $h_1$, $\ldots$, $h_N$ are continuous functions by Step \ref{S31}, $g$ is of class $C^1$ on $\R^N$ with partial derivatives given by formula \eqref{E42}. This completes the proof of Lemma \ref{L3}.
\end{proof}

\section{Asymptotics of the solitary waves}

This second part is mainly devoted to the proof of Theorem \ref{T1}. We first compute the optimal algebraic decay of the solitary waves stated in Theorem \ref{T7} by the standard argument mentioned in the introduction. In the second section, we complete the proof of Theorem \ref{T1} by deducing the asymptotics of a solitary wave from Propositions \ref{P7} and \ref{P8}. Finally, the last section is devoted to the proof of Theorem \ref{T2} which links the asymptotics of a solitary wave to its energy and its action in the case of the standard Kadomtsev-Petviashvili equation. In particular, the proof of Theorem \ref{T2} relies on Lemma \ref{L4}, which also implies the non-existence of non-trivial solitary waves when $p \geq \frac{4}{2N - 3}$ (see Corollary \ref{C2} for more details).

\subsection{Algebraic decay of the solitary waves}

In this section, we derive the algebraic decay of the solitary waves
stated in Theorem \ref{T7}. As mentioned in the introduction, Theorem \ref{T7} follows from a standard inductive argument which links the algebraic decay of the solitary waves to the algebraic decay of the associated kernels. We first determine some small algebraic decay for the solitary waves. It follows from Proposition \ref{P6} which gives some integral algebraic decay for $\nabla v$ and $\partial_1^2 v$. Then, we inductively improve the algebraic decay of $v$ and $\nabla v$ by using the superlinearity of equations \eqref{E11} and \eqref{E42}. This is possible as long as the rate of decay is less important than the rate of decay of $K_0$ and $K_k$. Thus, the solitary waves decay at least as fast as $K_0$, while their gradient decays at least as fast as $K_k$. This leads to Theorem \ref{T7} whose proof follows below.

\def\proof{\par{\it Proof of Theorem \ref{T7}}. \ignorespaces}
\begin{proof}
We split the proof into five steps. In the first one, we use Proposition \ref{P6} to infer some small algebraic decay for $v$ and $\nabla v$.

\setcounter{step}{0}
\begin{step}\label{S51}
There exists some $\alpha_0 > 0$ such that $v$ and $\nabla v$ belong to $M^\infty_{\beta}(\R^N)$ for $\beta \in [0, \alpha_0]$.
\end{step}

Step \ref{S51} results from equations \eqref{E10} and \eqref{E11}. Indeed, equation \eqref{E10} gives for every $\beta \geq 0$, and almost every $x \in \R^N$,
\begin{equation}\label{E110}
\begin{split}
|x|^\beta |v(x)| \leq & A \bigg( \int_{\R^N} |H_0(x-y)|
|y|^\beta |v(y)|^p |\partial_1 v(y)| dy \\ + & \int_{\R^N} |x-y|^\beta
|H_0(x-y)| |v(y)|^p |\partial_1 v(y)| dy \bigg).
\end{split}
\end{equation}
However, if $0 \leq \beta \leq N-1$, by Theorem \ref{T3} and Corollary \ref{C1}, we compute for every $1 < q < \frac{2N-1}{2N-4}$,
\begin{align*}
& \int_{\R^N} |x-y|^\beta |H_0(x-y)| |v(y)|^p |\partial_1 v(y)| dy \\
\leq & A \bigg( \int_{B(x,1)^c} |x-y|^{\beta - N + 1} |v(y)|^p
|\partial_1 v(y)| dy + \int_{B(x,1)} |H_0(x-y)| |v(y)|^p |\partial_1
v(y)| dy \bigg) \\ \leq & A \bigg( \int_{B(x,1)^c} |v(y)|^p
|\partial_1 v(y)| dy + \| H_0 \|_{L^q(B(0,1))} \| v^p \partial_1 v
 \|_{L^{q'}(B(x,1))} \bigg).
\end{align*}
Therefore, by Theorem \ref{T8}, there exists some $A \geq 0$ such that
\begin{equation}\label{E111}
\forall 0 \leq \beta \leq N-1, \int_{\R^N} |x-y|^\beta |H_0(x-y)|
|v(y)|^p |\partial_1 v(y)| dy \leq A.
\end{equation}
On the other hand, Theorem \ref{T3} and Corollary \ref{C1} yield for
every $1 < q < \frac{2N-1}{2N-4}$ and $r > \frac{N}{N-1}$,
\begin{align*}
& \int_{\R^N} |H_0(x-y)| |y|^\beta |v(y)|^p |\partial_1 v(y)| dy \\ \leq & A \bigg( \bigg( \int_{B(0,1)^c} |y|^{(1-N)r} dy \bigg)^\frac{1}{r} \bigg( \int_{B(x,1)^c} |y|^{\beta r'} |v(y)|^{p r'} |\partial_1 v(y)|^{r'} dy \bigg)^\frac{1}{r'} \\ + & \| H_0 \|_{L^q(B(0,1))} \bigg( \int_{B(x,1)} |y|^{\beta r'} |v(y)|^{p q'} |\partial_1 v(y)|^{q'} dy \bigg)^\frac{1}{q'} \bigg).
\end{align*}
However, by Theorem \ref{T8}, the function $v$ is bounded on $\R^N$,
so, by assertion \eqref{E45}, for every $s > 1$ and $0 < \beta < \min \{ 1, \frac{2}{s} \}$,
$$\int_{\R^N} |y|^{\beta s} |v(y)|^{p s} |\partial_1 v(y)|^s dy \leq A \bigg( \int_{\R^N} |y|^2 \partial_1 v(y)^2 dy \bigg)^\frac{\beta s}{2} \| \partial_1 v \|_{L^{\frac{2(1 - \beta) s}{2 - \beta s}}(\R^N)}^{(1 - \beta) s} \leq A.$$
Therefore, by Theorem \ref{T8} in the case $\beta = 0$, we obtain for $0 \leq \beta < \min \{ 1, \frac{2}{r'}, \frac{2}{q'} \}$,
\begin{equation}\label{E112}
\int_{\R^N} |H_0(x-y)| |y|^\beta |v(y)|^p |\partial_1 v(y)| dy \leq
A.
\end{equation}
Thus, by equations \eqref{E110}, \eqref{E111} and \eqref{E112}, there
is some $\alpha_1 > 0$ such that the function $v$ belongs to
$M^\infty_\beta(\R^N)$ for every $0 \leq \beta \leq \alpha_1$.

On the other hand, by Theorem \ref{T3} and Corollary \ref{C1}, the
kernel $K_0$ belongs to $L^q(\R^N)$ for $1 < q < \frac{2N-1}{2N-3}$,
while the function $v^p \nabla v$ belongs to $L^1(\R^N)$ by Theorem
\ref{T8}. Therefore, we can derive from equation \eqref{E11} the
following equation,
$$\nabla v = K_0 * (v^p \nabla v),$$
which holds in $L^q(\R^N)$ for $1 < q < \frac{2N-1}{2N-3}$. In
particular, it yields for every $\beta \geq 0$ and almost every $x \in
\R^N$,
\begin{equation}\label{E113}
\begin{split}
|x|^\beta |\nabla v(x)| \leq & A \bigg( \int_{\R^N} |K_0(x-y)|
|y|^\beta |v(y)|^p |\nabla v(y)| dy \\ + & \int_{\R^N} |x-y|^\beta
|K_0(x-y)| |v(y)|^p |\nabla v(y)| dy \bigg).
\end{split}
\end{equation}
However, if $0 \leq \beta \leq N$, Theorem \ref{T3} and Corollary
\ref{C1} give for every $1 < q < \frac{2N-1}{2N-3}$,
\begin{align*}
& \int_{\R^N} |x-y|^\beta |K_0(x-y)| |v(y)|^p |\nabla v(y)| dy \\ \leq
& A \bigg( \int_{B(x,1)^c} |x-y|^{\beta - N} |v(y)|^p |\nabla v(y)| dy
+ \| K_0 \|_{L^q(B(0,1))} \| v^p \nabla v \|_{L^{q'}(B(x,1))} \bigg).
\end{align*}
Hence, by Theorem \ref{T8}, there is some $A \geq 0$ such that
\begin{equation}\label{E114}
\forall 0 \leq \beta \leq N, \int_{\R^N} |x-y|^\beta |K_0(x-y)|
|v(y)|^p |\nabla v(y)| dy \leq A.
\end{equation}
On the other hand, Theorem \ref{T3} and Corollary \ref{C1} yield for
every $1 < q < \frac{2N-1}{2N-3}$ and $r > 1$,
\begin{align*}
& \int_{\R^N} |K_0(x-y)| |y|^\beta |v(y)|^p |\nabla v(y)| dy \\ \leq &
A \bigg( \bigg( \int_{B(0,1)^c} |y|^{- N r} dy \bigg)^\frac{1}{r}
\bigg( \int_{B(x,1)^c} |y|^{\beta r'} |v(y)|^{p r'} |\nabla v(y)|^{r'}
dy \bigg)^\frac{1}{r'} \\ + & \| K_0 \|_{L^q(B(0,1))} \bigg(
\int_{B(x,1)} |y|^{\beta r'} |v(y)|^{p q'} |\nabla v(y)|^{q'} dy
\bigg)^\frac{1}{q'} \bigg).
\end{align*}
However, by Theorem \ref{T8}, the function $v$ is bounded on $\R^N$,
so, by assertion \eqref{E45}, for every $s > 1$ and $0 < \beta < \min \{ 1, \frac{2}{s} \}$,
$$\int_{\R^N} |y|^{\beta s} |v(y)|^{p s} |\nabla v(y)|^s dy \leq A
\bigg( \int_{\R^N} |y|^2 |\nabla v(y)|^2 dy \bigg)^\frac{\beta s}{2}
\| \nabla v \|_{L^\frac{2(1 - \beta) s}{2 - \beta s}(\R^N)}^{(1 -
 \beta) s} \leq A.$$
Hence, (by Theorem \ref{T8} in the case $\beta = 0$), we get for every $0 \leq \beta < \min \{ 1, \frac{2}{r'}, \frac{2}{q'} \}$,
\begin{equation}\label{E115}
\int_{\R^N} |K_0(x-y)| |y|^\beta |v(y)|^p |\nabla v(y)| dy \leq A.
\end{equation}
Finally, by equations \eqref{E113}, \eqref{E114} and \eqref{E115}, there is also some $\alpha_2 > 0$ such that $\nabla v$ belongs to $M^\infty_\beta(\R^N)$ for every $0 \leq \beta \leq \alpha_2$. Then, it only remains to set $\alpha_0 = \min \{ \alpha_1, \alpha_2 \}$ to complete the proof of Step \ref{S51}.

\begin{remark}
By Step \ref{S51}, the function $v$ is lipschitzian on
$\R^N$. Indeed, $\nabla v$ is bounded on $\R^N$.
\end{remark}

We now improve the algebraic decay of $v$ by applying the inductive
argument mentioned in the introduction to equation \eqref{E11}.

\begin{step}\label{S52}
Let us consider $\alpha > 0$ and assume that the function $v$ belongs to $M^\infty_\beta(\R^N)$ for every $\beta \in [0, \alpha]$. Then, it belongs to $M^\infty_\beta(\R^N)$ for every $\beta \in [0, N] \cap [0, (p+1) \alpha[$.
\end{step}

Indeed, equation \eqref{E11} yields
\begin{equation}\label{E116}
|x|^\beta |v(x)| \leq A \bigg( \int_{\R^N} |x-y|^\beta |K_0(x-y)|
|v(y)|^{p + 1} dy + \int_{\R^N} |K_0(x-y)| |y|^\beta |v(y)|^{p +
 1} dy \bigg).
\end{equation}
for every $\beta > 0$ and for almost every $x \in \R^N$. However, by Theorem \ref{T3} and Corollary \ref{C1},
\begin{align*}
\int_{\R^N} |x-y|^\beta |K_0(x-y)| |v(y)|^{p + 1} dy \leq & A \bigg(
\int_{B(x,1)^c} |x-y|^{\beta - N} |v(y)|^{p + 1} dy \\ + &
\int_{B(x,1)} |K_0(x-y)| |v(y)|^{p + 1} dy \bigg),
\end{align*}
so that
\begin{equation}\label{E117}
\begin{split}
\int_{\R^N} |x-y|^\beta |K_0(x-y)| |v(y)|^{p + 1} dy \leq & A \Big( \| v \|^{p+1}_{L^{p+1}(\R^N)} + \| K_0 \|_{L^1(B(0,1))} \| v \|^{p+1}_{L^\infty(\R^N)} \Big) \\ \leq & A,
\end{split}
\end{equation}
for every $0 \leq \beta \leq N$, by Theorem \ref{T8}. Likewise, Theorem \ref{T3} and Corollary \ref{C1} give
\begin{align*}
\int_{\R^N} |K_0(x-y)| |y|^\beta |v(y)|^{p + 1} dy & \leq A \bigg( \bigg( \int_{B(x,1)^c} |y|^{\beta q'} |v(y)|^{(p + 1) q'} dy \bigg)^\frac{1}{q'} \\ & \bigg( \int_{B(0,1)^c} \frac{dz}{|z|^{N q}} \bigg)^\frac{1}{q} + \| K_0 \|_{L^1(B(0,1))} \| v \|^{p +  1}_{M^\infty_\frac{\beta}{p+1}(\R^N)} \bigg),
\end{align*}
for every $q > 1$. However, by the assumption of Step \ref{S52}, there is some $q > 1$ such that for every $\beta \in [0, (p+1) \alpha[$, the function $y \mapsto |y|^\beta v(y)^{p+1}$ belongs to $L^{q'}(\R^N)$. Hence, we obtain
\begin{equation}\label{E118}
\int_{\R^N} |K_0(x-y)| |y|^\beta |v(y)|^{p + 1} dy \leq A \Big( \| |.|^\beta v^{p + 1} \|_{L^{q'}(\R^N)} + \| v \|^{p + 1}_{M^\infty_\frac{\beta}{p+1}(\R^N)} \Big) \leq A.
\end{equation}
Thus, by equations \eqref{E116}, \eqref{E117} and \eqref{E118}, the
function $v$ belongs to $M^\infty_\beta(\R^N)$ for $\beta \in [0,N]
\cap [0, (p+1) \alpha[$, which is the desired result.

Then, we deduce the rate of decay of $v$ given by Theorem \ref{T7}.

\begin{step}\label{S53}
The function $v$ belongs to $M^\infty_N(\R^N)$.
\end{step}

Since $p+1 > 1$, the geometric sequence given by $u_0 = \alpha_0$ and $u_{n+1} = (p+1) u_n$ tends to $+ \infty$ when $n$ tends to $+ \infty$. Thus, by a straightforward inductive argument, it follows from Steps \ref{S51} and \ref{S52} that the function $v$ belongs to $M^\infty_N(\R^N)$.

We now turn to the algebraic decay of $\nabla v$. In particular, we improve the rate of decay given by Step \ref{S51} by applying the inductive argument of the introduction to equation \eqref{E42}.

\begin{step}\label{S54}
Let us consider $\alpha > 0$ and assume that the function $\nabla v$ belongs to $M^\infty_\beta(\R^N)$ for every $\beta \in [0, \alpha]$. Then, it belongs to $M^\infty_\beta(\R^N)$ for every $\beta \in [0, \min \{ N+1, (p+1) N, pN + \alpha \}]$.
\end{step}

Indeed, by Theorem \ref{T8}, the function $v^{p+1}$ is bounded and
continuous on $\R^N$. Moreover, by Steps \ref{S51} and \ref{S53}, it
belongs to $M^\infty_{(p+1)N}(\R^N)$, and its gradient, to
$L^\infty(\R^N)$. Therefore, by Lemma \ref{L3}, the following equality holds for every $k \in \{ 1, \ldots, N \}$ and $x \in \R^N$,
\begin{align*}
\partial_k v(x) = & i \int_{B(0,1)^c} K_k(y) v(x-y)^{p+1} dy + i
\int_{B(0,1)} K_k(y) (v(x-y)^{p+1} - v(x)^{p+1}) dy \\ + & \bigg(
\int_{\S^{N-1}} y_k K_0(y) dy \bigg) v(x)^{p+1}.
\end{align*}
In particular, by Theorem \ref{T3} and Corollary \ref{C1}, this yields for every $\beta > 0$,
\begin{equation}\label{E119}
\begin{split}
|x|^\beta |\partial_k v(x)| \leq & |x|^\beta \int_{B(0,1)^c} |K_k(y)|
|v(x-y)|^{p+1} dy + A |x|^\beta |v(x)|^{p+1} \\ + & |x|^\beta
\int_{B(0,1)} |K_k(y)| |v(x-y)^{p+1} - v(x)^{p+1}| dy.
\end{split}
\end{equation}
However, Step \ref{S53} yields for every $\beta \in [0, (p+1) N]$,
\begin{equation}\label{E120}
|x|^\beta |v(x)|^{p+1} \leq \| v
\|^{p+1}_{M^\infty_{\frac{\beta}{p+1}}(\R^N)} \leq A.
\end{equation}
On the other hand, by Theorem \ref{T3},
\begin{align*}
|x|^\beta \int_{B(0,1)^c} |K_k(y)| |v(x-y)|^{p+1} dy \leq & A
\bigg( \int_{B(0,1)^c} |y|^{\beta - N - 1} |v(x-y)|^{p+1} dy \\ + &
\int_{B(0,1)^c} |y|^{-N-1} |x-y|^\beta |v(x-y)|^{p+1} dy \bigg),
\end{align*}
so,
\begin{equation}\label{E121}
|x|^\beta \int_{B(0,1)^c} |K_k(y)| |v(x-y)|^{p+1} dy \leq A \Big( \| v \|^{p+1}_{L^{p+1}(\R^N)} + \| v \|^{p + 1}_{M^\infty_\frac{\beta}{p + 1}(\R^N)} \Big) \leq A,
\end{equation}
for every $\beta \in [0, \min\{ N+1, (p+1) N \}]$, by Theorem \ref{T8} and Step \ref{S53}. Finally, by Theorem \ref{T8} and Step \ref{S51}, the function $v$ is continuous on $\R^N$, while its gradient $\nabla v$ is bounded on $\R^N$. Hence, we compute for every $(x,y) \in (\R^N)^2$,
\begin{align*}
\Big| v(x-y)^{p+1} - v(x)^{p+1} \Big| \leq & (p+1) \int_0^1 |v(x-ty)|^p |\nabla v(x-ty).y| dt \\ \leq & A \sum_{l=1}^N |y_l| \int_0^1 \frac{dt}{1 + |x-ty|^{p N + \alpha}} \\ \leq & \frac{A}{1 + |x|^{pN + \alpha}} \sum_{l=1}^N |y_l|,
\end{align*}
which gives by Theorem \ref{T8}, Step \ref{S53} and the assumption of Step \ref{S54},
$$|x|^\beta \int_{B(0,1)} |K_k(y)| |v(x-y)^{p+1} - v(x)^{p+1}| dy \leq
\frac{A |x|^\beta}{1 + |x|^{pN + \alpha}} \int_{B(0,1)} \sum_{l=1}^N
|y_l| |K_k(y)| dy.$$
Then, Corollary \ref{C1} yields
\begin{equation}\label{E122}
|x|^\beta \int_{B(0,1)} |K_k(y)| |v(x-y)^{p+1} - v(x)^{p+1}| dy \leq A,
\end{equation}
for every $\beta \in [0, pN + \alpha]$. Thus, by equations \eqref{E119}, \eqref{E120}, \eqref{E121} and \eqref{E122}, the function $\nabla v$ belongs to $M^\infty_\beta(\R^N)$ for $0 \leq \beta \leq \min \{ N+1, (p+1)N, pN + \alpha \}$, which concludes the proof of Step \ref{S54}.

In conclusion, we infer the rate of decay of $\nabla v$ given by
Theorem \ref{T7}.

\begin{step}\label{S55}
The function $\nabla v$ belongs to the space $M^\infty_{\min \{ N+1,
 (p+1)N \}}(\R^N)$.
\end{step}

Indeed, the arithmetic sequence given by $u_0 = \alpha_0$ and $u_{n+1} = u_n + pN$ tends to $+ \infty$ when $n$ tends to $+ \infty$. Thus, by Steps \ref{S51} and \ref{S54}, the function $\nabla v$ belongs to $M^\infty_{ \min \{ N+1, (p+1)N \}}(\R^N)$. This concludes the proofs of Step \ref{S55} and of Theorem \ref{T7}.
\end{proof}

\subsection{Asymptotics of the solitary waves}

In this section, we complete the proof of Theorem \ref{T1}. Indeed, in the previous section, we proved Theorem \ref{T7} which describes the algebraic decay of the solitary waves and of their gradient. In order to show Theorem \ref{T1}, it remains to compute the first order asymptotics of any solitary wave $v$, i.e. the limit of the function $x \mapsto |x|^N v(x)$ when $|x|$ tends to $+ \infty$. This is the goal of Propositions \ref{P7} and \ref{P8}. In Proposition \ref{P7}, we compute the pointwise limit when $R$ tends to $+ \infty$ of the functions $v_R$ defined by formula \eqref{E46}. Our argument results from Theorems \ref{T3}, \ref{T5} and \ref{T7}, Corollary \ref{C1} and a standard application of the dominated convergence theorem. Then, in Proposition \ref{P8}, we deduce from Theorem \ref{T7} and a standard application of Ascoli-Arzela's theorem that this convergence is uniform in the case $p \geq \frac{1}{N}$. Finally, Theorem \ref{T1} follows from Theorem \ref{T7} and from Propositions \ref{P7} and \ref{P8}. However, let us first write the proof of Proposition \ref{P7}.

\def\proof{\par{\it Proof of Proposition \ref{P7}}. \ignorespaces}
\begin{proof}
Let $\sigma \in \S^{N-1}$. By formula \eqref{E46}, we compute for
every $R > 0$,
\begin{equation}\label{E123}
\begin{split}
v_R(\sigma) = & \frac{1}{p+1} \bigg( \int_{B(R\sigma, \frac{R}{2})}
R^N K_0(R \sigma - y) v(y)^{p+1} dy \\ + & \int_{B(R\sigma,
 \frac{R}{2})^c} R^N K_0(R \sigma - y) v(y)^{p+1} dy \bigg).
\end{split}
\end{equation}
However, Theorem \ref{T5} gives
$$R^N K_0(R\sigma - y) v(y)^{p+1} \underset{R \to + \infty}{\to}
\frac{\Gamma(\frac{N}{2})}{2 \pi^\frac{N}{2}} (1 - N \sigma_1^2)
v(y)^{p+1},$$ 
for every $y \in \R^N$, while
$$\Big| R^N K_0(R\sigma - y) v(y)^{p+1} 1_{|R\sigma - y| \geq \frac{R}{2}} \Big| \leq A \frac{R^N}{(R \sigma - y)^N (1 + |y|^{N(p+1)})} \leq \frac{A}{1 + |y|^{N(p+1)}},$$
by Theorems \ref{T3} and \ref{T7}. Therefore, the dominated convergence theorem yields
\begin{equation}\label{E124}
\int_{B(R\sigma, \frac{R}{2})^c} R^N K_0(R \sigma - y) v(y)^{p+1} dy
\underset{R \to + \infty}{\to} \frac{\Gamma(\frac{N}{2})}{2
 \pi^\frac{N}{2}} (1 - N \sigma_1^2) \int_{\R^N} v(y)^{p+1} dy.
\end{equation}
On the other hand,
\begin{align*}
\bigg| \int_{B(R\sigma, \frac{R}{2})} R^N K_0(R \sigma - y)
v(y)^{p+1} dy \bigg| \leq & A \int_{B(R\sigma, \frac{R}{2})}
\frac{R^N}{|y|^{N(p+1)}} |K_0(R \sigma - y)| dy \\ \leq &
\frac{A}{R^{N p}} \bigg( \int_{B(0,1)} |K_0(z)| dz + \int_{1 < |z| <
 \frac{R}{2}} \frac{dz}{|z|^N} \bigg),
\end{align*}
by Theorem \ref{T7}, so that
$$\bigg| \int_{B(R\sigma, \frac{R}{2})} R^N K_0(R \sigma - y)
v(y)^{p+1} dy \bigg| \leq \frac{A}{R^{N p}} (1 + \ln(R)) \underset{R
 \to + \infty}{\to} 0,$$
by Theorem \ref{T3} and Corollary \ref{C1}. Thus, equations \eqref{E123} and \eqref{E124} lead to
$$v_R(\sigma) \underset{R \to + \infty}{\to}
\frac{\Gamma(\frac{N}{2})}{2 \pi^\frac{N}{2} (p + 1)} (1 - N \sigma_1^2) \int_{\R^N} v(y)^{p+1} dy,$$
which is exactly assertion \eqref{E47}.
\end{proof}

Then, we establish the uniformity of the pointwise limit computed above in the case $p \geq \frac{1}{N}$. This follows from Proposition \ref{P8} whose proof is mentioned below.

\def\proof{\par{\it Proof of Proposition \ref{P8}}. \ignorespaces}
\begin{proof}
Let us assume by contradiction that $(v_R)_{R > 0}$ does not converge uniformly to $v_\infty$ when $R$ tends to $+ \infty$. Then, there is some $\eps > 0$ and a sequence of positive real numbers $(R_n)_{n \in  \N}$ tending to $+ \infty$, such that
\begin{equation}\label{E125}
\forall n \in \N, \| v_{R_n} - v_\infty \|_{L^\infty(\S^{N-1})} \geq
\eps.
\end{equation}
However, since $p \geq \frac{1}{N}$, we deduce from Theorem \ref{T7}
\footnote{Here, the notation $\nabla^{\S^{N-1}}$ denotes the gradient on the sphere $\S^{N-1}$ immersed in $\R^N$. More precisely, if we consider some index $i \in \{ 1, \ldots, N \}$ and some function $f \in C^\infty(\S^{N-1},\C)$, the notation $\partial_i^{\S^{N-1}}$ is defined by
$$\forall x \in \S^{N-1}, \partial_i^{\S^{N-1}} f(x) = \underset{t \to 0}{\lim} \frac{f(\frac{x+te_i}{|x+te_i|}) - f(x)}{t},$$
where $(e_1,\ldots,e_n)$ is the canonical basis of $\R^N$. Then, the gradient $\nabla^{\S^{N-1}} f$ is defined by $\nabla^{\S^{N-1}} f = (\partial_1^{\S^{N-1}} f, \ldots, \partial_N^{\S^{N-1}} f)$.},
\begin{equation}\label{E126}
\begin{split}
\| v_{R_n} \|_{L^\infty(\S^{N-1})} \leq & A, \\ \| \nabla^{\S^{N-1}}
v_{R_n} \|_{L^\infty(\S^{N-1})} \leq & A R_n^{N+1} \| \nabla v(R_n.)
\|_{L^\infty(\S^{N-1})} \leq \frac{A}{R_n^{\min\{ 1, pN \} -1}} \leq
A.
\end{split}
\end{equation}
Therefore, by equation \eqref{E126} and Ascoli-Arzela's theorem, up to a subsequence, $(v_{R_n})_{n \in \N}$ converges in
$L^\infty(\S^{N-1})$. By Proposition \ref{P8}, its limit is necessary equal to $v_\infty$, which leads to a contradiction with assertion \eqref{E125}. Thus, $(v_R)_{R > 0}$ uniformly converges to $v_\infty$ when $R$ tends to $+ \infty$, which is the desired result.
\end{proof}

Finally, we end the proof of Theorem \ref{T1} by invoking Theorem \ref{T7}, and Propositions \ref{P7} and \ref{P8}.

\def\proof{\par{\it Proof of Theorem \ref{T1}}. \ignorespaces}
\begin{proof}
Indeed, the function $x \mapsto |x|^N v(x)$ is bounded on $\R^N$ by
Theorem \ref{T7}. Moreover, by Proposition \ref{P7}, assertion
\eqref{E8} holds for every $\sigma \in \S^{N-1}$, and by Proposition
\ref{P8}, the convergence given by assertion \eqref{E8} is actually
uniform when $\frac{1}{N} \leq p < \frac{4}{2N - 3}$. This concludes
the proof of Theorem \ref{T1}.
\end{proof}

\subsection{Link between the asymptotics and the energy of the
 solitary waves for the standard Kadomtsev-Petviashvili equation}

This last section deals with the standard Kadomtsev-Petviashvili equation. Here, we link the asymptotics of a solitary wave to its energy and its action. As mentioned in the introduction, this link results from the standard Pohozaev identities of Lemma \ref{L4}, which were derived by A. de Bouard and J.-C. Saut in \cite{deBoSau1}.

\def\proof{\par{\it Proof of Theorem \ref{T2}}. \ignorespaces}
\begin{proof}
Let us consider $k \in \{ 2, \ldots, N \}$ and $0 < p < \frac{4}{2N-3}$. We deduce from formulae \eqref{E49}, \eqref{E50} and \eqref{E51} after some algebraic computations that
\begin{align}
\Big( 4 + p(3 - 2N) \Big) \int_{\R^N} v(x)^{p+2} dx = & 2(p+1)(p+2)
\int_{\R^N} v(x)^2 dx, \label{E127} \\ \Big( 4 + p(3 - 2N) \Big)
\int_{\R^N} \partial_1 v(x)^2 dx = & pN \int_{\R^N} v(x)^2 dx,
\label{E128} \\ \Big( 4 + p(3 - 2N) \Big) \int_{\R^N} v_k(x)^2 dx = & p \int_{\R^N} v(x)^2 dx. \label{E129}
\end{align}
Assume now that $p=1$, and $N=2$ or $N=3$. Then, formula \eqref{E7} becomes
\begin{equation}\label{E130}
\forall \sigma \in \S^{N-1}, v_\infty(\sigma) =
\frac{\Gamma(\frac{N}{2})}{4 \pi^\frac{N}{2}} (1 - N \sigma_1^2)
\int_{\R^N} v(x)^2 dx.
\end{equation}
However, by formulae \eqref{E3}, \eqref{E127}, \eqref{E128} and
\eqref{E129}, the energy of $v$ is equal to
\begin{equation}\label{E131}
E(v) = \frac{2N-5}{2(7-2N)} \int_{\R^N} v(x)^2 dx,
\end{equation}
so, by equation \eqref{E130},
$$\forall \sigma \in \S^{N-1}, v_\infty(\sigma) = \frac{(7 - 2N)
 \Gamma(\frac{N}{2})}{2 (2N-5) \pi^\frac{N}{2}} (1 - N \sigma_1^2)
E(v),$$
which is formula \eqref{E9}. Likewise, by equations \eqref{E2} and \eqref{E131}, the action of $v$ is given by
$$S(v) = \frac{1}{7 - 2N} \int_{\R^N} v(x)^2 dx,$$
so, by equation \eqref{E130},
$$\forall \sigma \in \S^{N-1}, v_\infty(\sigma) = \frac{(7 - 2N)
 \Gamma(\frac{N}{2})}{4 \pi^\frac{N}{2}} (1 - N \sigma_1^2)
S(v),$$
which completes the proof of formula \eqref{E9}.
\end{proof}

Finally, for the sake of completeness, we complete this subsection by another straightforward consequence of the identities of Lemma \ref{L4}: the non-existence of non-trivial solitary-wave solutions of equation \eqref{E1} in $Y$ for $p \geq \frac{4}{2N - 3}$ in any dimension $N \geq 4$. 

\def\proof{\par{\it Proof of Corollary \ref{C2}}. \ignorespaces}
\begin{proof}
Indeed, formula \eqref{E128} gives for every $p > 0$,
$$\Big( 4 + p(3 - 2N) \Big) \int_{\R^N} \partial_1 v(x)^2 dx = Np \int_{\R^N}
v(x)^2 dx.$$
Therefore, if $p \geq \frac{4}{2N-3}$,
$$\int_{\R^N} v(x)^2 dx \leq 0.$$
Thus, $v$ is identically equal to $0$, which concludes the proof of
Corollary \ref{C2}.
\end{proof}

\section*{Acknowledgements}

The author is thankful to F. B\'ethuel for its careful attention to the preparation of this paper, and to P. G\'erard, M. Maris, J.-C. Saut and N. Tzvetkov for interesting and helpful discussions. He is also grateful to the referee who provided many improvements to the article.

\bibliographystyle{plain}
\bibliography{Bibliogr}

\end{document}